\magnification=\magstephalf
\parindent = 0 pt
\input amstex
\documentstyle{amsppt}
\refstyle{A}
\widestnumber\key{CDW}
\def\pf{{\blacksquare}}

\def\fh{\frak{h}}
\def\fg{\frak{g}}

\topmatter
\title On dynamical Poisson groupoids I\endtitle
\author Luen Chau Li and
Serge Parmentier\endauthor
\address{L.-C. Li, Department of Mathematics, Pennsylvania State University,
University Park, PA 16802. USA}\endaddress
\email luenli\@math.psu.edu\endemail
\address S. Parmentier, Institut Girard Desargues
(UMR 5028 du CNRS),  Universit\'e Lyon 1, 43 Blvd du 11 Novembre 1918, 69622
Villeurbanne cedex. France\endaddress
\email serge\@desargues.univ-lyon1.fr\endemail
\date {September 17, 2002}\enddate
\keywords Poisson groupoids, duality, dynamical $r-$
matrices, symplectic double groupoids\endkeywords

\abstract We address the question of duality for the dynamical Poisson
groupoids of Etingof and Varchenko over a contractible base.  We also
give an explicit description for the coboundary case associated with the
solutions
of (CDYBE) on simple Lie algebras as classified by the same
authors.  Our approach is based on
the study of a class of Poisson structures on trivial Lie groupoids within
the category of biequivariant Poisson manifolds. In the former case,
it is shown that the dual Poisson groupoid of such a dynamical Poisson
groupoid is isomorphic to a Poisson groupoid (with trivial Lie groupoid
structure) within this category.  In the latter case, we find that
the dual Poisson groupoid is also of dynamical type modulo Poisson
groupoid isomorphisms.  For the coboundary dynamical Poisson groupoids
associated with constant $r-$ matrices, we give an explicit construction
of the corresponding symplectic double groupoids. In this case, the
symplectic leaves of the dynamical Poisson groupoid are shown to be the
orbits of a Poisson Lie group action.

\endabstract
\endtopmatter

\document
{\bf 1. Introduction.}
\bigskip
The classical dynamical Yang-Baxter equation (CDYBE)
 was introduced by
Felder in \cite{F} as a consistency condition for the
Knizhnik-Zamolodchikov-Bernard equations.
The geometric meaning of (CDYBE) was subsequently unraveled by Etingof and
Varchenko in the fundamental paper [EV].
While the solutions of the classical Yang-Baxter equation are related to
Poisson 
Lie groups \cite{D},
the authors in [EV] showed that an appropriate geometrical setting for
(CDYBE)
is that of a special class of
Poisson groupoids (as defined in \cite{W1}), the so-called
coboundary dynamical Poisson groupoids.
Given a Lie group $G$, a Lie subgroup $H\subset G$, an $Ad_H^{*}$ invariant
open set $U\subset \fh^*$ (here $\fh^*$ is the dual of $\fh=Lie(H))$,  and a
solution of (CDYBE), Etingof and Varchenko constructed a Poisson bracket on
the set
$X=U\times G\times U$ compatible with its trivial Lie groupoid structure.
This Poisson bracket
intertwines a left
and a right inclusion
of the restricted symplectic cotangent $H\times U$ into $X$ together with a
Sklyanin-like term on
$G.$  In addition, the authors in \cite{EV} identified an appropriate
abstract
context in which to
view these objects as  the category of $H$-bi-equivariant Poisson
manifolds ${\Cal C}_{U}.$
\smallskip
It is classical  that the study of Poisson Lie groups relies in an
essential way on duality and the construction of doubles \cite{D, STS,
LW1}.
For Poisson groupoids,  the notion of duality was introduced by
Weinstein in \cite{W1},
and  was developed by MacKenzie and Xu in \cite{MX1, MX2}.
In the same paper \cite{W1}, Weinstein also introduced the notion of
symplectic double groupoids (see also \cite{M2}), and described a program
for
showing that, at least locally, Poisson groupoids in duality arise as the
base of a symplectic double groupoid.
\smallskip
In order to state our objectives and results, let us begin by recalling
that a symplectic groupoid is a pair $(\Gamma, \Pi)$, consisting of a Lie
groupoid $\Gamma$ together with a non-degenerate Poisson structure
$\Pi$, in such a way that the graph of the multiplication map is a
Lagrangian submanifold of $\Gamma\times\Gamma\times {\overline \Gamma}$
\cite{W2,K}.
It is a classical fact that Poisson structures can be understood at
least locally by the notion of symplectic groupoids.  On the other hand,
double groupoids are intrinsically complicated objects introduced by
Ehresmann \cite{E} in the 1960's and have found usage in category theory
\cite{E}, homotopy theory \cite{BH}, differential geometry \cite{P}, and
Poisson groups \cite{M3, LW2}.  By definition, a double Lie groupoid is
a quadruple  $({\Cal  S}; {\Cal H}, {\Cal  V}, B)$ where ${\Cal  H}$ and
${\Cal  V}$ are Lie groupoids over
$B$, and 
${\Cal  S}$ is equipped with two Lie groupoid structures, a horizontal
structure with base ${\Cal  V}$, and a vertical structure with base ${\Cal
H}$, such that the structure maps of each groupoid structure on ${\Cal  S}$
are morphisms with respect to the other. Finally, a symplectic double
groupoid is a double Lie groupoid $({\Cal S}; {\Cal H}, {\Cal V}, B)$
in which $\Cal S$ is equipped with a symplectic structure such that
both $\Cal S\rightrightarrows\Cal V$ and $\Cal S\rightrightarrows\Cal H$
are symplectic groupoids.  Note that for the case of Poisson Lie groups,
the program in \cite{W1} which we mentioned above has been carried out
globally in \cite{LW2}.  Thus a Poisson Lie group and its dual are the
bases of a symplectic groupoid.
\smallskip
This work is the first part of a series to understand the geometry of
dynamical Poisson groupoids.  Our goal here is three-fold.  First
of all, for a general dynamical Poisson groupoid $X=U\times G\times U$
(not necessarily of coboundary type), we would
like to characterize certain  properties of its (global) dual Poisson
groupoid  in a simple nontrivial case in which its existence is
guaranteed.  
In this connection, we should point out that in contrast to
(finite dimensional) Lie algebras,
not all Lie algebroids can be integrated to Lie groupoids \cite{AM}.
For (finite dimensional) general Lie algebroids, the necessary and
sufficient 
condition for  integrability  was only obtained quite recently
in  \cite{CF}.  Thus we work at the outset with
the class ${\Cal C}_\ast$ of Poisson groupoids $X=U\times G\times U$
(with the trivial Lie groupoid structure) which admits a (base preserving)
Poisson groupoid
morphism $I:H\times U\longrightarrow X$, where $U$ is $Ad_H^*$ invariant
and contractible.  If $X$ is a dynamical Poisson groupoid in
${\Cal C}_\ast$, the corresponding Lie algebroid dual $A(X)^*$ must be
transitive, i.e., the anchor map is a surjective submersion.  Consequently,
we can invoke a general theorem of MacKenzie \cite{M1}, according to
which $A(X)^* \simeq TU\oplus (U\times \fg')$, where the latter is
the trivial Lie algebroid over $U$ and $\fg'$ is a typical fiber of
the adjoint bundle of $A(X)^*$.  As a result, $A(X)^*$ integrates
to a unique global Lie groupoid $X^*$ isomorphic to the trivial Lie
groupoid $U \times G'\times U$, where $G'$ is the connected and
simply connected Lie group with $Lie(G')=\fg'$.  Thus the existence of
the dual Poisson groupoid is not an issue.  Our main result
in this direction (Theorem 3.2.4) is the following: if $X$ is a dynamical
Poisson groupoid in ${\Cal C}_\ast$, then its dual Poisson groupoid $X^*$
is isomorphic to a Poisson groupoid $(U\times G'\times U, \{\, ,\,
\}_{U\times
G'\times U})$ in ${\Cal C}_\ast$.  In particular, the Poisson
structure $\{\, , \,\}_{U\times G'\times U}$ is uniquely determined by a
(unique) Poisson groupoid morphism
$I':H\times U\longrightarrow U\times G'\times U$
and a unique groupoid 1-cocycle $P'$ on $U\times G'\times U$.  The proof
of this theorem consists of two steps: in the first step, we establish
the existence of the unique Poisson groupoid morphism $I'$; while the
second step involves a careful analysis of the form of the Poisson
bracket for a Poisson groupoid in ${\Cal C}_\ast$(Theorem 2.2.5).
As a corollary of
Theorem 3.2.4, we obtain via Poisson reduction a reduced duality diagram
for the Poisson quotients $G/H \times U$ and $G'/H \times U$ and for
the vertex Lie algebras $\fg$ and $\fg'$.  In the special
case when $\fh^*=0$, this duality diagram is just the well-known diagram of
Drinfeld for Poisson Lie groups.
\smallskip
In [EV], extending Belavin and Drinfeld's classic paper \cite{BD},
Etingof and Var-
chenko obtained a classification
of solutions of (CDYBE) for pairs $(\fg, \fh)$ of Lie algebras,
where $\fg$ is simple and $\fh\subset \fg$ is a Cartan subalgebra. These
solutions of
(CDYBE) are parametrized by subsets $S$ of a simple system of roots
$\Delta^s$ and closed
meromorphic two-forms on $\fh^*.$
\smallskip 
Our second objective is to give an explicit study of duality
for the coboundary
dynamical Poisson groupoids associated with this class of dynamical
$r-$ matrices.
 Note that in this
case, the base $U$ (where the $r-$ matrix is analytic) is neither
contractible
nor simply-connected. We proceed in two steps. To start with,  we construct
(see Theorem 4.4) an explicit
trivialization of the Lie algebroid dual $A(X)^*$ of the (full) coboundary
Poisson 
groupoid $X=U\times G\times U.$
This, in particular, establishes the integrability of $A(X)^*$ as a
Lie algebroid.  Then an argument similar to that of Theorem
3.2.4 
applied to any connected and simply connected open subset
$U'$ of $U$ shows that the dual Poisson groupoid of $U'\times G\times U'$ is
isomorphic to a dynamical Poisson
groupoid $U'\times G'\times U'.$
Here, the vertex Lie group $G'$ is a semi-direct product
$L_S\ltimes I_S$ where $L_S\subset G$ is the Levi factor and $I_S$ is
a normal Lie subgroup containing the product $N^+_S\times N^-_S$ of
unipotent radicals.  More importantly, the Poisson bracket is uniquely
determined by the value of a Lie groupoid $1$-cocycle
$P': U'\times G'\times U'\longrightarrow L({\fg'}^*, \fg')$ whose
partial derivatives are explicitly given in terms of the Lie-Yamaguti
data of the reductive pair $(\fg,\fh).$
\smallskip
Our final objective in this paper is to understand how to construct
symplectic double groupoids for the coboundary dynamical case in the
special instance where the $r-$ matrix is constant.  For this class of coboundary
dynamical Poisson groupoids, the base is $\fh^*$ and so Theorem 3.2.4
applies.
However, from the point of view of constructing the symplectic double
groupoids, it is more natural (and considerably simpler) to work
directly with the dual Poisson groupoid whose Lie algebroid is
$T^*\fh^*\times \fg^*$.  Since we have a constant $r-$ matrix, the Lie
group $G$ equipped with the Sklyanin bracket is a Poisson Lie group
(for simplicity, we assume $G$ is complete) and as it turns out, the
dual Poisson groupoid of $X$ is given by $X^*=H\times \fh^*\times G^*$
with appropriate structure maps ($G^*$ is the dual Poisson group of $G$)
and the Poisson structure is a product structure (Theorem 5.1.4).
The construction of a symplectic double groupoid having $X$ and $X^*$ as
side groupoids proceeds via a number of steps (Proposition 5.2.3, Corollary
5.2.6, Corollary 5.2.8,Theorem 5.2.10 and 5.2.13).
First of all, we show $X$ and $X^*$
form a matched pair of Lie groupoids.  The upshot of this is that
$X$ and $X^*$ act on each other via groupoid actions and give rise
to a vacant double Lie groupoid $({\Cal S}_{vac}; X^*, X, \fh^*)$.
However, this is not the correct underlying double Lie groupoid of
the sought-for symplectic double groupoid (in contrast to the Poisson
group case).  In the second step of the construction, we extend the
Lie groupoids $X$ and $X^*$ to the product groupoids $X^*_e=X^*\times H
\rightrightarrows H\times \fh^*$ and $X_e=(H\times H)\times
X\rightrightarrows
H\times \fh^*$ ($H\times H \longrightarrow H$ is the coarse groupoid).  Then
we show that there is a left action of $X^*_e$ on $X$ and a right action
of $X_e$ on $X^*$.  The corresponding action groupoids
${\Cal S}\simeq X^*_e\ltimes X\rightrightarrows X$ and
${\Cal S}\simeq X^*\rtimes X_e\rightrightarrows X^*$
 then give the horizontal structure and the vertical
structure respectively of a nonvacant double Lie groupoid
$({\Cal S}; X^*, X, \fh^*)$ which has $({\Cal S}_{vac}; X^*, X, \fh^*)$
as a double Lie subgroupoid.  Finally, we show that
the double Lie groupoid $({\Cal S}; X^*, X, \fh^*)$ where $\Cal S$
is equipped with an appropriate symplectic structure is a desired
symplectic double groupoid.  We would like to point out that
the actions of the extended Lie groupoids on the unextended ones
obey a number of properties (Proposition 5.2.9) which are important
in showing that $({\Cal S}; X^*, X, \fh^*)$ is a double Lie groupoid.
The reader should contrast these properties with actions via `twisted
automorphisms' (Proposition 5.2.4, \cite{M3, LW1}).  As an
application/amplification of this result, we show the existence of
a natural Poisson Lie group structure on the set $H\times H\times G^*$
such that the symplectic leaves of $(X, \{\, , \,\}_{X})$ are the
orbits of a Poisson action of $H\times H\times G^*$ on $X$(Theorem 5.2.28).
Finally, we use this result to describe the symplectic leaves of a
natural Poisson quotient associated with $X$.
\smallskip
The paper is organized as follows.  In Section 2, we begin by giving some
background material which we recall here for the convenience of the reader.
The rest of Section 2 is devoted to the description of all Poisson
groupoids $(X, \{\, , \,\})$ which admit a Poisson groupoid morphism
$I:H\times U\longrightarrow X=U\times G\times U$, where $X$ is the trivial
Lie groupoid over a connected base $U$.  Section 3 is concerned with
Poisson groupoids in duality with dynamical Poisson groupoids over a
contractible base $U$.  It also treats duality diagrams for the
Poisson quotients mentioned earlier and for the vertex Lie algebras.
In Section 4, we consider the coboundary dynamical Poisson groupoids
associated with a class of solutions of (CDYBE) for pairs $(\fg,\fh)$
of Lie algebras, where $\fg$ is simple, and $\fh$ is a Cartan subalgebra
of $\fg$ \cite{EV}.  Here, we obtain a more refined description of the dual
Poisson groupoid.  Finally, Section 5 treats the coboundary dynamical
Poisson groupoids in the constant $r-$ matrix case
in detail.  We begin with an explicit description of the dual Poisson
groupoid $X^*$ whose Lie algebroid is $T^*\fh^*\times \fg^*$.  Then
we move on to the construction of a symplectic double groupoid having
$X$ and $X^*$ as side groupoids.  We conclude the paper by describing
the symplectic leaves of $(X, \{\, , \,\})$ as well as a Poisson quotient
associated with $X$.
\smallskip
We shall address the construction of symplectic double groupoids for the
general dynamical case, together with its relationship to other works (in
particular \cite{LWX}) in a sequel to
this paper. On the other hand, the links between duality and the recent
work \cite{KW}, as well as the relevance of coboundary dynamical Poisson
groupoids to integrable systems (see the papers \cite{HM, LX1, LX2} in this
connection)
will be considered in separate publications.
\bigskip
\bigskip
{\bf Acknowledgements.}  L.-C. Li would like to thank the members of 
Institut G. Desargues for hospitality and CNRS support (UMR 5028) during his visits to Universit\'e  Lyon 
1.
\bigskip
\bigskip
\parindent=0pt
{\bf 2. A class of biequivariant Poisson groupoids.}
\bigskip
{\bf 2.1. Preliminaries.}
\bigskip
In this preliminary subsection, we recall some of the basic concepts and
constructs
which we shall use in this paper (other results will be recalled when
needed). 
\bigskip
Let $\Gamma$ be a Lie groupoid over $B$ (see \cite{DSW,M1} for details),
with target and source
maps $\alpha, \beta: \Gamma\longrightarrow B$, and multplication map
$m:\Gamma\ast 
\Gamma\longrightarrow B$ defined on the set of composable pairs
$\Gamma\ast 
\Gamma := \{(x,y)\, \mid \, \beta (x)=\alpha (y)\}$. We shall denote
the  unit section by $\epsilon : B\longrightarrow \Gamma$, and the
inversion map by
$i:\Gamma\longrightarrow \Gamma.$
\medskip
{\bf Definition 2.1.1} \cite{W1} (Poisson groupoid.)
\smallskip
{\it A Lie groupoid  $\Gamma$ equipped with a Poisson structure $\Pi$ is
called a Poisson groupoid if
and 
only if  the graph of the multiplication map

$$Gr (m)\subset \Gamma\times \Gamma\times {\overline \Gamma}$$ is a
coisotropic 
submanifold, i.e. if and only if
$$\big( \Pi\oplus \Pi \oplus -\Pi\big) \,  (\omega, \omega')= 0,\quad
\forall \omega, 
\omega'\in (T (Gr (m)))^\perp\subset T^* (\Gamma\times
\Gamma\times\Gamma).$$

$\Gamma$ is called a symplectic groupoid if $\Pi$ is non degenerate with
$Gr(m)$ a Lagrangian submanifold.
\smallskip
In both cases, we shall say that the Poisson structure and the groupoid
structure are compatible.}
\bigskip
Let $G$ be a connected Lie group, $H\subset G$ a connected Lie subgroup with
respective Lie
algebras $\fg$ and $\fh$ and let $U\subset \fh^*$ be a connected $Ad_H^*-$
invariant
open subset.
In \cite{EV}, Etingof and Varchenko introduced the category ${\Cal
C}_{U}$
of biequivariant 
Poisson manifolds over $U$ as follows.
\smallskip
An object in ${\Cal C}_{U}$ is a Poisson manifold $(X, \{\, , \,
\}_X)$ equipped with commuting left Hamiltonian $H-$action $\phi^-$ and
right
Hamiltonian $H-$ action $\phi^+$ with $U$-valued $Ad_H^*$ equivariant
momentum maps
$j_\pm: X\longrightarrow
U$ satisfying the polarity condition
$$\{j_+^* \varphi, \,
j_-^* \psi\}_X= 0,\, {\hbox { for all }}\varphi, \psi\in C^\infty (U).$$
A morphism in ${\Cal C}_{U}$ between $(X, \{\, , \, \}_X)$ and
$(X', \{\, , \, \}_{X'})$ is an equivariant Poisson map $\sigma:
X\longrightarrow X'$ such that $j'_\pm\, \circ\, \sigma =j_\pm.$

\medskip
{\bf Definition 2.1.2}\cite{EV} (Poisson groupoid in ${\Cal C}_{U}$)
\smallskip
{\it We say that $X\in {\Cal C}_{U}$ is a Poisson groupoid
in ${\Cal C}_{U}$ if it is equipped with a compatible groupoid
structure over 
$U$ such that $\alpha= j_-,\, \beta= j_+.$}
\bigskip
{\bf Example 2.1.3} (The Hamiltonian unit)
\smallskip
The most basic (but not simplest) symplectic groupoid in ${\Cal
C}_{U}$ is the (restricted) Hamiltonian unit $H\times U$ equipped  with the
nondegenerate bracket
$$\{f,g\} (h,p)= - <D'g, \delta f> + <D'f, \delta g> - <p, [\delta f,
\delta g]>,$$

$(<D'f, Z>= {d\over dt}_{\mid_0} f( he^{tZ}, p),\quad <\delta f,
\lambda>= {d\over dt}_{\mid_0}  f(h, p+ t\lambda),\, Z\in \fh, \lambda\in
\fh^*)$
the $H- $ actions 
$$\phi^- _k (h,p)= (kh, p),\quad \phi^+ _k (h,p)= (hk, Ad^*_k p),$$
and the action groupoid structure
$$\eqalign {&\alpha_0 (h,p)= j_- (h,p)= Ad^*_{h^{-1}} p ,\quad \beta_0
(h,p)= 
j_+ 
(h,p)= p\cr
&(h, j_- (k,q))\cdot (k,q)= (hk, q),\quad \epsilon (q) = (1,q),\quad
i(h,p)= (h^{-1}, Ad^*_{h^{-1}} p).\cr}$$
If $U=\fh^*$, this is clearly isomorphic to the cotangent symplectic
groupoid $T^*H$ [W2]
under the trivialization map.
\medskip
We now recall a fundamental construction of [EV] which interprets dynamical
$r-$ matrices in terms of Poisson groupoids.

 Let $\iota: \fh\longrightarrow \fg$ be the Lie inclusion. We say
that a
smooth map $R:
U\longrightarrow L(\fg^*, \fg)$ (here and henceforth we denote by
$L(\fg^*,\fg)$ the set of linear maps from $\fg^*$ to $\fg$)  is a classical
dynamical $r-$ matrix
if and only if  it is skew symmetric
$$<R(q)(A), B>=- <A, R (q) B>,$$
and satisfies the classical dynamical Yang- Baxter condition
$$\eqalign {& dR (q) \iota^* A \, B - dR (q) \iota^* B \, A +
\iota 
d<R(q) A, B>\cr
& -[R (q) A, R(q) B] - R (q) ad^*_{R (q) (A)} B + R (q)
ad^*_{R (q) B} A= \chi (A,B),\cr}\eqno (2.1.1)$$
where $\chi : \fg^*\times \fg^* \longrightarrow  \fg$ is
$ad_{\fg} -$ invariant, i.e.
$$ad_X \, \chi (A,B) + \chi (ad^*_X A, B) + \chi (A, ad^* _X B)=
0,$$
for all $A, B\in \fg^*, X\in \fg$, and all $q\in U.$

\smallskip
The dynamical $r-$ matrix is said to be $ad^*_{\fh} -$ equivariant
if and only if
$$dR(q) ad^*_Z q + R (q) ad^*_{\iota (Z)} + ad_{\iota
(Z)} \, R (q) = 0, \eqno (2.1.2)$$
for all $Z\in \fh,$ and all $q\in U.$
\medskip
Note that if $\chi (A,B) = 0$ in Eqn. (2.1.1), the resulting equation is
called the classical dynamical Yang-Baxter equation (CDYBE) \cite{F}.  On
the other hand, if 
$\chi (A,B) = [T (A), T (B)]$ for some nonzero symmetric map $T:
\fg^*\longrightarrow \fg$ with $ad_X T +  T ad^*_X= 0$, the resulting
equation is called the modified dynamical Yang-Baxter equation (mDYBE).

\medskip

Let $X= U \times G\times U$. For $f\in C^\infty (X)$, define its
partial derivatives and left/right gradients (w.r.t. $G$) by
$$\eqalign {&<\delta_1 f, \lambda>= {d\over dt}_{\mid _0} f(p+t\lambda, x,
q),\quad <\delta_2 f, \lambda>= {d\over dt}_{\mid_0} f(p, x, q+t\lambda),\,
\lambda\in \fh^*\cr
&<Df, X>= {d\over dt}_{\mid_0} f(p, e^{tX}x,
q),\quad <D'f, X>= {d\over dt}_{\mid_0} f(p, xe^{tX}, q), \, X\in g.\cr}$$
\medskip
We shall equip $X$ with the trivial Lie groupoid structure over $U$
with
structure maps 
$$\eqalign {&\alpha (p,x,q)= p,\, \beta (p,x,q)=q,\, \epsilon (q) = (q,
1,q), \, 
i(p,x,q)= (q, x^{-1}, p)\cr
&\hskip 60 pt m ((p,x,q), (q,y, r))= (p, xy, r).\cr}\eqno (2.1.3)$$

\medskip

The following 
 theorem gives the Poisson groupoid analog of coboundary Poisson Lie groups
(in the context of trivial Lie groupoids over $U).$
\smallskip

{\bf Theorem 2.1.4}\cite{EV} {\it
\smallskip
(a) The formula $$\eqalign {&\{f,g\}_X (p,x,q)= <p, [\delta_1f,\,
\delta_1g]> -<q, [\delta_2f,\, \delta_2g]>\cr
&\hskip 50 pt - <\iota\delta_1 f,\, Dg> - <\iota\delta_2 f, D'g>\cr
&\hskip 50 pt + <\iota\delta_1 g,\, Df> + <\iota\delta_2 g , D'f>\cr
&\hskip 50 pt +<R(p) Df, Dg> - <R(q) D'f, D'g>\cr}$$
defines a Poisson bracket on $X$ if and only if
$R: U\longrightarrow L(\fg^*, \fg)$ is an $ad^*_{\fh} -$
equivariant dynamical $r -$ matrix.

(b) The trivial Lie groupoid $X$ equipped with the Poisson bracket $\{\, ,
\,
\}_X$ in (a) and the
Hamiltonian $H- $ actions
$$\phi^- _h (p,x,q)= (Ad^*_{h^{-1}} p, hx, q), \quad \phi^+_h
(p,x,q) = (p, xh, Ad^*_{h} q),$$
 is a Poisson groupoid in
${\Cal 
C}_{U}.$
\smallskip
We shall call $(X, \{\, ,\, \}_X)$ a coboundary dynamical
Poisson groupoid.}
\medskip
Note that the dynamical Poisson groupoid $X$ of Thm 2.1.5 admits a Poisson
groupoid embedding 
$$I: H\times U \longrightarrow X: (h,p)\mapsto (Ad^*_{h^{-1}} p,
h, p).\eqno (2.1.4)$$
where $H\times U$ is the Hamitonian unit. As we shall see in later
sections, this property turns out to play a
crucial role in the study of duality.
\bigskip 
We now recall the notion of a Lie algebroid (for more details see [DSW],
[M1]).
\medskip
{\bf Definition 2.1.5} {\it A Lie algebroid is a smooth vector bundle $q:
A\longrightarrow B$
equipped with a Lie bracket $[\, , \,
]_A$ on the set 
$\Gamma (A)$ of smooth sections of $A$ and a smooth base preserving bundle
map $a: 
A\longrightarrow TB$, called the anchor map, such that
$$
\align & a \, [\zeta, \, \eta]_A= [a(\zeta), a(\eta)]_{B} \\
& [ \zeta, \,f\,  \eta]_A=f\, [\zeta, \eta]_A+ a (\zeta) (f) \,\eta,\
\endalign
$$
for all $\zeta, \eta\in \Gamma (A)$ and all $f\in C^\infty (B).$}
\bigskip

The Lie algebroid of a smooth groupoid $\Gamma$ over $B$ is the vector
bundle $$A(\Gamma):=\big(Ker (T\alpha)\big)_{\mid _{\epsilon (B)}}$$
 over $B$ 
with anchor map $a$ given by the restriction of $T [\alpha, \beta]$ to
$A(\Gamma)$ (here $[\alpha, \beta] (z)$
= $(\alpha (z), \beta (z)),\, z\in
\Gamma$) 
and bracket of sections $[X, Y] (b):= [X^l , Y^l]_{\Gamma} (\epsilon (b))$
where
$$X^l: \Gamma\longrightarrow Ker (T\alpha)$$
is the unique left invariant vector field whose restriction to $\epsilon
(B)$ is $X.$
\medskip
Let $V$ be a vector space and let $\rho: \Gamma\longrightarrow Aut (V)$
be a smooth groupoid morphism where $Aut (V)$ is viewed as a groupoid over
its
unit element $I_V.$
\medskip
{\bf Definition 2.1.6} {\it A smooth map $\Sigma: \Gamma\longrightarrow V$
is 
called a groupoid $1-$cocycle iff
$$\Sigma (xy)= \Sigma (x) + \rho (x) \Sigma (y)$$
for all $(x,y)\in \Gamma\ast \Gamma.$
The induced map $\Sigma_\ast: A(\Gamma) \rightarrow V$ defined as
the restriction of $T \Sigma$ to $A(\Gamma)$ is called the
induced  Lie algebroid $1-$ cocycle.}
\bigskip
Finally, we recall the notion of an action of a Lie groupoid
$\Gamma\rightrightarrows B$ on a manifold
$S$ with moment map $f: S\longrightarrow B.$ (We use the terminology of
\cite{MW}.)
\smallskip
Let
$$
\align
&\Gamma \ast _f S=\{ (x, s)\in \Gamma\times S\, \mid\, \beta (x)= f(s)\}\\
&S \ast_f \Gamma=\{ (s,x)\in S\times \Gamma\, \mid\, f(s)= \alpha (x)\}\\
\endalign
$$
\medskip

{\bf Definition 2.1.7}
\smallskip
{\it (a) A left action of $\Gamma$ on $S$ with moment $f$ is a smooth map
$\phi^l: \Gamma \ast_f S\longrightarrow S: (x,s)\mapsto x\cdot s$ such that
$$f(x\cdot s)= \alpha (x),\quad y\cdot (x\cdot s)= (yx)\cdot s,\quad
\epsilon 
(f(t)) \cdot t=t,$$
for all $(y,x)\in \Gamma \ast \Gamma,\, (x,s)\in \Gamma \ast_f S, t\in S.$
\medskip
(b) A right action of $\Gamma$ on $S$ with moment $f$ is a smooth map
$\phi^r: S \ast_f\Gamma \longrightarrow S: (s,x)\mapsto s\cdot x$ such that
$$f(s\cdot x)= \beta (x),\quad (s\cdot x)\cdot y= s \cdot (xy),\quad t\cdot
\epsilon (f(t))= t,$$
for all $(x,y)\in \Gamma \ast \Gamma, (s,x)\in S \ast_f \Gamma, t\in S.$}

\bigskip 
\bigskip

{\bf 2.2. Trivial Lie groupoids in ${\Cal C}_{U}.$}
\bigskip

Our purpose in this subsection is to provide an explicit class of Poisson
brackets 
on trivial Lie groupoids
which extends the construction of  thm
2.1.4 (a), and is essential for our subsequent study of duality.
\medskip
We assume that the Lie subgroup $H\subset G$ is connected.
We begin with a general property.
\medskip
{\bf Proposition 2.2.1} {\it Let $Y$ be a Poisson groupoid over $U$
with source and target maps
$\alpha$ and $\beta$ and unit map $\epsilon.$  If there exists a (base
preserving) 
Poisson groupoid morphism
$$I: H\times U \rightarrow Y,$$
(here $H\times U$ is the Hamiltonian unit) then $Y$
belongs to ${\Cal C}_{U}.$}
\medskip
{ Proof.}\hskip 5 pt  It follows from a general property of Poisson
groupoids [W1] that
$$\{\alpha^*
\varphi, \beta^* \psi\}_Y= 0,\, \forall
\varphi,\psi\in C^\infty (U).$$ So it remains to show that $Y$ admits
two commuting 
Hamiltonian $H - $ actions $\phi^-,\, \phi^+$ with equivariant momentum maps
$\alpha$ and $\beta.$
\smallskip
Set, as in Def. 2.1.7,
$$
\align & (H\times U)\ast_\alpha Y= \{ (h,p,y)\mid \beta_0 (h,p)= \alpha
(y)\} \\
&Y\ast_\beta (H\times U)=\{(y,h,p)\mid \beta (y)= \alpha_0 (h,p)\}
\endalign
$$
Here $\alpha_0, \beta_0$ are as in Example 2.1.3.

The morphism $I$ induces a left (resp. right) groupoid action of  $H\times
U$ on $Y$ over $\alpha$ (resp. over $\beta$):
$$
\align &\phi^-: (H\times U)\ast_\alpha Y\longrightarrow Y\\
&(h,\alpha (y), y)\mapsto I(h,\alpha (y))\cdot y\\
&\phi^+ : Y\ast_\beta (H\times U)\longrightarrow Y\\
&(y,h,Ad^*_h \beta (y))\mapsto y\cdot I (h, Ad^*_h \beta
(y))\
\endalign
$$
which, upon the natural identifications
$$(H\times U)\ast_\alpha Y\simeq H\times Y, \,Y\ast_\beta (H\times
U)\simeq Y\times H,$$
induce a left and  a right action of $H$ on $Y$ (also denoted $\phi^\pm ).$
\medskip
We now show that $\phi^-$ is Hamiltonian with momentum map $\alpha.$ (The
verification for $\phi^+$ and $\beta$ is similar.)
\smallskip
Note that $\alpha$ is equivariant since $\alpha (\phi^- _k (y))= \alpha
(I(k,
\alpha (y))\cdot y)= \alpha
(I(k, \alpha (y)))=\alpha_0 (k,\alpha (y))$= $Ad^*_{k^{-1}} \alpha (y).$
\smallskip
Let $Z^- (y) = {d\over dt}_{\mid_0} \phi^- _{e^{tZ}} (y)$ be the
infinitesimal generator of the action corr. to $Z$. We want to show $Z^- $
coincides with the Hamiltonian vector field ${\widehat X}_{f_Z\circ
\beta}$ where $f_Z\in C^\infty (U)$ is defined by $f_Z (q)= <Z,q>,\,
\forall q\in U.$
\smallskip
Since $I(1,q)= \epsilon (q),$ we have

$$Z^- (y) = {d\over dt}_{\mid_0} I(e^{tZ}, \alpha (y))\cdot y= T_{\epsilon
\circ \alpha (y)} r_y \, T_{(1,\alpha (y))} I (Z,0),$$
which shows that $Z^- \in Ker T\beta.$ On the other hand,
$$\eqalign {Z^- (\epsilon\circ \alpha (y))&= {d\over dt}_{\mid_0} I(e^{tZ},
\alpha (y)) \cdot I (1, \alpha (y))\cr
&= {d\over dt}_{\mid_0} I(e^{tZ} , \alpha (y))\cr
&= T_{(1,\alpha (y))} I (Z,0),\cr}$$
thus $Z^- $ is right invariant. Since ${\widehat X}_{f_Z\circ \alpha}$ is
also right invariant \cite{X}, it suffices to show that both vector fields
coincide on $\epsilon (U).$ But from the Poisson property of $I$, we
have
$$\eqalign { {\widehat X}_{f_Z\circ \alpha}(\epsilon\circ \alpha (y))&=
\Pi_Y^\# (\epsilon\circ \alpha (y)) d(f_Z\circ \alpha)\cr
&= T_{(1, \alpha (y))} I\, \Pi_0^\# (1, \alpha (y)) d( f_Z\circ \alpha\circ
I)\cr
&= T_{(1,\alpha (y))}I\,\Pi_0^\# (1, \alpha (y)) d(f_Z\circ \alpha_0)\cr
& =  T_{(1,\alpha (y))}I\, \Pi_0^\# (1, \alpha (y)) (ad^*_Z \alpha (y),
Z)\cr
&= T_{(1, \alpha (y))}I\, (Z,0).\cr}$$
 Now, since $Z^-$ is Hamiltonian, its flow $\phi^-_{e^Z}$
at $t=1$ preserves the Poisson bracket of $Y$. Therefore, the connectedness
of $H$ implies that $\phi^-$ is Hamiltonian with momentum map $\alpha .$
Hence the claim.$\pf$
\bigskip
For the rest of this subsection, we let $X= U\times G\times U$
be 
the trivial Lie groupoid of section
2.1 (see eqn (2.1.3)).
We shall describe all pairs

$$(X,\{\, , \, \}),\quad I: H\times U\rightarrow X,$$
where $\{\, , \, \}$ is a Poisson bracket on $X$ compatible with its
groupoid structure and $I$ is a morphism of Poisson groupoids, where
$H\times U$ is the Hamiltonian unit.
\bigskip
Let $\rho: G\longrightarrow Aut(V)$ be a representation of $G$ on
the vector space $V.$ We are going to restrict ourselves to groupoid
$1 $ cocycles $P:X\longrightarrow V$ which satisfy
$$P(p,xy, q)= P(p, x, r) + \rho (x) P(r, y, q) {\hbox { for all }}
p,q,r\in U, x,y\in G.$$
\medskip
{\bf Proposition 2.2.2}  {\it $P$ is a $1-$ cocycle on $X$ iff
$$P(p,x,q)= -l(p) + \pi (x) + \rho (x) l(q),$$
for some smooth map $l: \fh^*\longrightarrow V$ with $l(q_0)= 0$ for some
$q_0\in U$, and a
group $1- $ 
cocycle $\pi: G\longrightarrow V.$}
\medskip
{ Proof.}   Clearly any such map is a $1-$ cocycle.
Conversely, if $P$ is a cocycle then $P(p,1,p)=P(p,1,q_0)+P(q_0,1,p)=0$ and
$\pi
(x)= P(q_0,x,q_0)$ is a group
cocycle. The claim then follows from $ P(p, x, q)= P(p, x,q_0)+ \rho (x)
P(q_0,
1, q)= 
P(p,1, q_0) + P(q_0,x,q_0) +\rho (x) P(q_0,1,q) = -P(q_0,1,p) + \pi (x)
+\rho (x)
P(q_0,1,q). \pf$
\bigskip
{\bf Proposition 2.2.3}
{\it Let $\Pi\in \bigwedge ^2 T X$ be a bivector field. Then the graph of
$m:X\ast X\longrightarrow X$ is $\Pi -$ coisotropic in $X\times X\times
{\overline 
X}$ iff
$$
\align &\Pi^\# (p,x,q) (Z_1, B, Z_2)= (-K(p) Z_1 -A^*
(p)
T^*_1 r_x B, \\
&\hskip 15 pt T_1 r_x A (p) Z_1 + T_1 l_x A (q) Z_2+ T_1 r_x
P(p,x,q) T^*_1 r_x B,\\
&\hskip 30 pt K(q) Z_2 - A^* (q) T^*_1 l_x B),\\
\endalign
$$
for some smooth maps $K: U\longrightarrow L(\fh,\fh^*),\, A:
U\longrightarrow
L(\fh,\fg)$, and  a  groupoid $1- $ cocycle $P:X\longrightarrow L(\fg^*,
\fg)$
for the adjoint action.  Here, $K$ and $P$ are pointwise skew-symmetric.}
\medskip
{ Proof.} The graph of the multiplication $m$ is
$$Gr (m)= \{\big( (p,x,q), (q,y,r), (p,xy,r)\big)\}\subset X\times
X\times {\overline X}.$$
Therefore, $\Omega\in\big( T_{\big( (p,x,q), (q,y,r), (p,xy,r)\big)}
Gr(m)\big)^\perp$ 
if and only if
$$\Omega=\big( (Z_1, \omega , Z_2), (-Z_2, T^*_y (r_{y^{-1}}\circ l_x)
\omega, Z_3) , (-Z_1, - T^*_{xy} r_{y^{-1}} \omega, -Z_3)\big),$$
for some $Z_1,Z_2,Z_3\in h$ and $\omega\in (T_x G)^*.$
\smallskip
One then verifies (see the appendix for the details) that
the $\Pi - $ coisotropy of $Gr(m)$:
 $$(\Pi \oplus \Pi  \oplus -\Pi) (\Omega, \Omega ')=
0,\quad \forall \Omega, \Omega'\in (T Gr(m))^\perp,$$
is equivalent to our assertion. $\pf$
\bigskip 
Now, a map  $I: H\times U\longrightarrow X$ is a (base preserving)
groupoid morphism 
iff $$I(k,q)= (Ad^*_{k^{-1}} q, \chi (k,q), q)\eqno (2.2.1)$$
 for some smooth map $\chi$ satisfying
$$\chi (hk, q)= \chi (h, Ad^*_{k^{-1}} q) \chi (k,q),\eqno (2.2.2)$$
in particular $\chi (1,q)=1$, and if $0\in U$, the map $H\longrightarrow G:
k\mapsto \chi (k,0)$ is a
group
morphism.
\smallskip
Note that Eqn. (2.2.2) says that $\chi: H\times U\longrightarrow G$ is a
groupoid morphism when $G$ is viewed as a groupoid over its unit element.
Applying the Lie functor to $\chi$ then provides an algebroid morphism
$$\eqalign {&A(\chi): U\times \fh\longrightarrow \fg\cr
& (q,Z)\mapsto T_{(1,q)} \chi (Z, ad^*_Z q),\cr}$$

which we shall henceforth denote as $(q,Z)\mapsto A_\chi (q) Z.$
The morphism property then says that for all $Z,Z'\in \fh, p\in U$,
$${\eqalign {A_\chi (p) [Z,Z']= dA_\chi &(p)\cdot ad^*_Z p\cdot Z' -
dA_\chi (p)\cdot ad^*_{Z'} p
\cdot Z\cr
&+ [A_\chi (p) Z, A_\chi (p) Z'].\cr}}\eqno (2.2.3)$$
\bigskip
{\bf Proposition 2.2.4} {\it If $(X, \Pi)$ is a Poisson
groupoid with $\Pi^\#$  expressed
as in Prop. 2.2.3 above, then the map $I$ is a Poisson map iff
\smallskip
(a) $\quad K(p) Z= ad^*_Z p,\quad \forall Z,Z'\in \fh.$
\smallskip
(b) $\quad A_\chi (p) = A(p)$
\smallskip
(c)  For all $\alpha, \beta\in \fg^*$,
\smallskip
  $ <\alpha, P (I (h,p))\beta> = <\lambda_\alpha, Z_\beta  >
-<\lambda_\beta,  Z_\alpha > - <p, [Z_\alpha, Z_\beta]>,$
\smallskip
where $\lambda_\alpha\in \fh^*$ and $ Z_\alpha\in \fh$ are defined by
$$\eqalign {&<\lambda_\alpha, Z>=<\alpha,  T_{(h,p)} (r_{(\chi (h,p)
)^{-1}}\circ \chi) (T_1l_h Z,0)>,\quad \forall Z\in \fh\cr
&<Z_\alpha, \lambda>= <\alpha, T_{(h,p)} (r_{(\chi (h,p)
)^{-1}}\circ \chi) (0,\lambda)>,\quad \forall \lambda\in \fh^*.\cr}$$}
\medskip
{Proof.} This is a spelled out form of $\{I^* f,I^* g\}_{H\times
U}= I^* \{f, g\}_X . \pf$
\bigskip
Assembling the above propositions, we can now formulate the main assertion
of
this 
subsection.
\bigskip
{\bf Theorem 2.2.5} {\it If $(X,\{\, ,\, \})$ is a
Poisson groupoid which admits
a Poisson groupoid morphism $I:H\times U\rightarrow X$ as in Eqn. (2.2.1),
then,
\smallskip
(a)  The Poisson bracket must be of
 the form
\smallskip
$$\eqalign {&\{f,g\} (p,x,q)= <p,[\delta_1 f, \delta_1 g]> -
<q,[\delta_2f,\delta_2g]>\cr
&\hskip 30 pt -<A_\chi (p)\delta_1f, Dg> -<A_\chi (q) \delta_2 f, D'g>\cr
&\hskip 30 pt +<A_\chi (p) \delta_1g,
Df>+<A_\chi (q)\delta_2 g, D'f>\cr
&\hskip 30 pt+ <Df, P(p,x,q) Dg>,\cr}$$
where the groupoid $1-$ cocycle $P$ satisfies Prop. 2.2.4 (c)
\medskip
(b) The Jacobi identity for $\{\, , \, \}$ is equivalent to the condition

$$\eqalign {& <\beta, [P \alpha,P \gamma]>- <\beta, DP \cdot P\alpha
\cdot\gamma> 
+<\beta, \delta_1 P \cdot (A_\chi(p)^* \alpha )\cdot \gamma>\cr
& + <\beta , \delta_2 P
\cdot (A_\chi(q)^* Ad^*_x \alpha)\cdot \gamma>
 + cp(\alpha, \beta, \gamma)= 0\cr}$$
for all $ \alpha,\beta,\gamma\in \fg^*,$
where $P$ stands for $P(p,x,q)$ and
$$
\align
&\delta_1 P (\lambda)=
{d\over dt} _{\mid 0} P(p+t\lambda, x, q),\, \delta_2P (\lambda)= {d\over
dt}_{\mid 0} P(p,x,q+t\lambda),\\
&\hskip 40 pt DP\cdot X= {d\over dt}_{\mid 0} P(p, e^{tX}x, q).\\
\endalign
$$
\medskip
(c) $(X,\{\,,\,\})$ belongs to ${\Cal C}_{U}$ with Hamiltonian $H- $
actions 
$$\eqalign {&\phi^- : H\times X\longrightarrow X\cr
& (h, p,x,q)\mapsto (Ad^*_{h^{-1}}
p, \chi (h,p) x, q)\cr
&\phi^+ : X\times H\longrightarrow X\cr
& (p,x,q,h)\mapsto (p, x\chi(h,
Ad^*_h q) , Ad^*_h q).\cr}$$}
\medskip
{ Proof. } (a) This is  simply a restatement of Prop. 2.2.3 and Prop. 2.2.4.

\smallskip
(b) We give 
the main steps (see the appendix for the details of the calculations):
\smallskip
 For the Jacobi identity, we use the shorthand notation
$J_{ijk},\, i,j,k\in \{1,2,*\}$ to stand for $ \{p_i^* \varphi, \{p_j^*
\varphi',\, p_k^* \psi\}\} + c.p.$,
 where, as an index, $* = G$, and $p_1$, $p_G$, $p_2$ are the projections
from
$X = U\times G \times U$ onto the first, second,
and third factor of $X$ respectively.
Thus for  example $J_{12*}= \{p_1^* \varphi, \{p_2^*
\varphi',\, p_G^* \psi\}_\} + c.p..$
\smallskip
Clearly, we have $J_{ijk}= 0$ for $i,j,k\in \{1,2\}$ and $J_{* 12}= 0$ and
these
do not impose any conditions. On the other hand, $J_{* 11}= 0
\Leftrightarrow J_{* 22}= 0 \Leftrightarrow A_\chi$
satisfies (2.2.3). Writing $P(p,x,q)= -l(p) + \pi (x) + Ad_x l(q)
Ad^*_x$ as in Prop.
2.2.2, we have
\smallskip

$$\eqalign { &\hskip 40 pt J_{* * 1}= 0\Leftrightarrow J_{**
2}= 0 
\Leftrightarrow\cr
& <\alpha, (dl(p) ad^*_Zp  + ad _{A_\chi (p) Z}
l(p)  + l(p) ad^*_{A_\chi (p)Z} + d\pi (1) A_\chi(p) Z) \beta>\cr
& =  + <dA_\chi(p) (A_\chi(p)^*
\alpha) Z, \beta>  - <dA_\chi (p) (A_\chi(p)^* \beta) Z, \alpha>,\quad
\forall \alpha,\beta\in \fg^*.\cr}\eqno (\star)$$
But the latter follows upon differentiating  the identity of Prop. 2.2.4
(c) at $(1,p)$. Indeed, that $lhs (\star)={d\over dt}_{\mid 0} <\alpha,
P\circ I
(e^{tZ},
p) \beta>$ is clear.
On the other hand, upon using
$$
\align
&<\lambda, Z_\alpha (1,p) >=<\alpha,  {d\over dt}_{\mid_0} \chi
(1,p+t\lambda)\, \chi (1,p)^{-1}>= 0\\
&<\lambda_\alpha (1,p), Z>= <\alpha, {d\over dt}_{\mid_0} \chi (e^{tZ}, p)
\chi (1,p)^{-1}>= <A_\chi (p)^* \alpha, Z>,\\
\endalign
$$
we have
$$
\align
&{d\over dt}_{\mid_0} \big( <\lambda_\alpha, Z_\beta> (e^{tZ},p)-
<\lambda_\beta, 
Z_\alpha> (e^{tZ},p) -<p, [Z_\alpha (e^{tZ},p),\, Z_\beta (e^{tZ},p)]>
\big)\\
&= <A_\chi(p)^* \alpha, {d\over dt}_{\mid_0} Z_\beta (e^{tZ}, p)>-<A_\chi
(p)^* 
\beta, {d\over dt}_{\mid_0} Z_\alpha (e^{tZ}, p)>\\
&=rhs (\star).\\
\endalign
$$ 
Thus, with our assumptions, the bracket $\{\,,\,\}$ satisfies the Jacobi
identity if and only
if $J_{***}= 0$ which is precisely (b).
\medskip
(c) This is Prop. 2.2.1: The morphism $I:H\times U\rightarrow X$
induces a left groupoid
action of $H\times U$ on $X$ over $\alpha_X$

$$\phi^- :(H\times U)\ast_{\alpha_X} X\longrightarrow X: ((h,p),
(p,x,q))\mapsto I(h,p)\cdot (p,x,q)$$
and a right groupoid action on $X$ over $\beta_X$
$$\phi^+ : X\ast_{\beta_X} (H\times U)\longrightarrow X: ((p,x,q), (k,
Ad^*_k q) )\mapsto (p,x,q)\cdot I (k, Ad^*_k q).$$
With the identifications $(H\times U)\ast_{\alpha_X} X\simeq H\times
X$ and $X\ast_{\beta_X} (H\times U)\simeq X\times H$, these actions are
the
ones given above. $\pf.$
\bigskip
We end this subsection with a definition and some remarks.
\bigskip
{\bf Definition 2.2.6}
\smallskip
{\it Following Etingof and Varchenko, we shall say that the Poisson
groupoid $(X,\{\, , \, \})$ of Theorem
2.2.5 is of  {\bf dynamical} type iff $I(k,q)= (Ad^*_{k^{-1}} q,
k, q).$  In this case, the corresponding Lie algebroid dual $A(X)^*$
(which is a Lie algebroid) will be called a dynamical Lie algebroid.}

\bigskip
{\bf Remarks 2.2.7}
\smallskip
(a) For $X$ of dynamical type, $\chi (h,p)= \chi (h,1)= h,$ thus
$A_\chi (p) Z=  Z.$ Therefore, the first six terms of the Poisson
bracket coincide with those of the coboundary case. The last term
however, which is given by
$$P(p,x,q)= -l(p) + \pi (x) + Ad_x l(q) Ad_x^*$$
differs from the coboundary case by the group $1-$ cocycle $\pi:
G\rightarrow L(\fg^*,
\fg).$
\smallskip
As we shall demonstrate in section 4 for a class of solutions of the
modified dynamical Yang-Baxter equation on simple Lie algebras, Poisson
groupoids of dynamical type
with $\pi \not=0$ arise naturally as Poisson groupoid duals (modulo Poisson
groupoid isomorphisms) of certain
coboundary dynamical groupoids. Note that the situation here is
analogous to that for Poisson Lie groups: Typically, the Poisson Lie group
dual of a 
Lie group equipped with the Sklyanin bracket is not of coboundary type.

\medskip
(b)  For $X$ of dynamical type, the identity of Prop. 2.2.4 (c) simplifies
to 
  $$P(Ad^*_{h^{-1}} p, h, p) = 0$$
 In other words $P$ vanishes
on the $H-$ orbit of $\epsilon (U)\subset
X.$ 
This condition is the natural extension of the $\fh-$ equivariance of the
dynamical $r-$ matrix which it reduces to when $\pi (x)=0.$
\bigskip
\bigskip 

{\bf 3. Duality.}
\bigskip
The purpose of this section is to characterize the Poisson groupoid dual
to a dynamical Poisson groupoid $X= U\times G\times U$ where $U\subset
\fh^*$ is an 
$Ad^*_H -$ invariant contractible open set, and study
some 
derived duality diagrams for Poisson quotients.

\bigskip
{\bf 3.1. Duality of Poisson groupoids.}
\bigskip
Following \cite{W1},\cite{M2}, \cite{MX2}, we begin by recalling the
notion of duality of Poisson groupoids and the
definition of the 
dual (when it exists) of a Poisson groupoid.
\medskip
Let $(Y,\{\, , \, \}_Y)$ be a Poisson groupoid over $B$ with target
and source maps 
$\alpha, \beta,$ and unit map $\epsilon.$ We use, as above, the sign
convention $\{f,g\}_Y=
<df,\, \Pi^\# dg>.$
\smallskip
Since $Y$ is Poisson, the set of $1-$ forms $\Omega^1 (Y)$ inherits a Lie
bracket from $C^\infty (Y)$ \cite{W1}, \cite {KSM}, given by
$$[\omega, \omega']= -L_{\Pi^\# \omega} \omega' + L_{\Pi^\# \omega'}
\omega - d<\omega, \Pi^\# \omega'>,\eqno (3.1.1)$$

and the map
$$\Omega^1 (Y) \longrightarrow \frak X (Y): \omega\mapsto -\Pi^\# \omega$$
is a morphism of Lie algebras, where $\frak X (Y)$ is the set of vector
fields
on $Y$ with ordinary Lie bracket. Therefore, $T^* Y$ is a Lie algebroid
over $Y.$
\medskip
Now, it follows from a general result in \cite{W1} that the unit submanifold
$\epsilon (B)$ of the Poisson groupoid $Y$ is
coisotropic in $Y$, hence its conormal bundle
$$N^*\big(\epsilon (B)\big):= \bigcup_{q\in B} \big(
T_{\epsilon (q)} \epsilon
(B)\big)^\perp \subset T^* Y_{\mid \epsilon (B)}$$
inherits a Lie algebroid structure:
the bracket of two sections $\theta_1, \theta_2: B \longrightarrow
N^*(\epsilon (B))$ is
$$[\theta_1, \theta_2]_{N^*} (q)= [{\overline \theta_1},
{\overline\theta_2}] \epsilon (q)\eqno
(3.1.2)$$
for arbitrary ${\overline \theta_1}, {\overline \theta_2}\in Ker
\,\epsilon^*$ 
subject to
${\overline \theta_1} \circ \epsilon = \theta_1, {\overline \theta_2}\circ
\epsilon = \theta_2,$
while the anchor map $a_*: N^*\big(\epsilon
(B)\big)\rightarrow TB$ is given by
the restriction of $-\Pi^\#$ to $N^*\big(\epsilon (B)\big).$
\smallskip
Since we have a natural identification $N^*\big(\epsilon (B)\big)\simeq
A(Y)^*$, we shall therefore always take $A(Y)^*$ with the induced Lie 
algebroid structure and
the pair $(A(Y), A(Y)^*)$ will be called the tangent Lie bialgebroid of
$(Y, \{\, ,\, \}).$ For the precise definition of Lie bialgebroids see
[MX1]; note however that $(A,A^*)$ is a Lie bialgebroid if and
only if $(A^*, A)$ is.
\bigskip
{\bf Definition 3.1.1} \cite{M2}
{\it We shall say that two Poisson groupoids $Y$ and $Y'$
over the same base are in duality if and only if
the Lie bialgebroids $(A(Y), A(Y)^*)$ and $(A(Y')^*, A(Y')^-)$ are
isomorphic. Here, $A(Y')^-$ is obtained from $A(Y)$ by
changing the sign of both anchor and bracket of sections.}
\bigskip 
Note that if the Lie algebroid $A(Y)^*$ is integrable, then there exists 
(by Lie $I$)
a unique source-simply connected Lie groupoid $Y^*$ integrating
$A(Y)^*$.  In this case, it follows
 from a general theorem of MacKenzie and Xu \cite{MX2} that the
latter 
may be
equipped with a unique Poisson bracket $\{\, , \, \}_{Y^*}$ compatible
with its groupoid structure such that $(Y^*,\{\, , \, \}_{Y^*})$
has tangent Lie bialgebroid $(A(Y)^*, A(Y)).$
\smallskip
The Poisson groupoids $Y$ and $Y^*$ are Poisson groupoids in
duality and $(Y^*, \{\, , \, \}_{Y^*})$ is called the
dual of $(Y, \{\, , \, \}_Y).$
\medskip
The following example is already in \cite{W1}:
\medskip
{\bf Example 3.1.2}
The Poisson groupoid dual to the Hamiltonian unit $H\times
U$ of Example 2.1.3
 is the coarse groupoid $U\times U$ with Poisson bracket
$$\{f, g\}_{U\times U} (p,q)= <p,[\delta_1 f , \delta_1 g]>-
<q, [\delta_2 f, \delta_2 g]>.$$

Note that $U\times U$ belongs to ${\Cal C}_{U}$ with the
Hamiltonian $H-$ actions
$$\phi^- _k (p,q)= (Ad^*_{k^{-1}} p, q),\quad \phi^+ _k
(p,q)= (p, Ad^*_k q).$$
The associated tangent Lie bialgebroid is given by
$$\big(A(H\times U), A(U\times U)\big)= \big(U\times
\fh, 
U\times \fh^*\big)$$
and the respective Lie brackets on
 smooth sections are as follows
$$\eqalign { & [Z_1,Z_2] (q)= [Z_1(q), Z_2(q)]_h + dZ_2 (q) ad^* _{Z_1(q)}
q- 
dZ_1
(q) ad^*_{Z_2(q)} q, \quad Z_1,Z_2: U\rightarrow \fh\cr
&[X_1, X_2] (q)= dX_2 (q) X_1 (q) - dX_1 (q) X_2 (q),\quad
X_1, X_2: U \rightarrow \fh^*.\cr}$$
\bigskip
We now recall a special (and simplest) instance of Lie bialgebroid
morphisms. 
Let $(A,A^*)$ and $(A', (A')^*)$ be two Lie bialgebroids over
$B$ with bundle projections $q: A\rightarrow B, \, q_*:
A^*\rightarrow 
B$, anchors $a: A\rightarrow TB,\, a_*: A^* \rightarrow TB$ , and
similarly for $(A', (A')^*).$
\bigskip
{\bf Definition 3.1.3}

{\it  A bundle map
$$
\align & \phi: A\longrightarrow A' \\
&\hskip 10 pt q \searrow\hskip 10 pt \swarrow q' \\
&\hskip 25 pt B \
\endalign
$$
is called a Lie bialgebroid morphism if and only if
both $\phi$ and $\phi^*$ are Lie algebroid morphisms, i.e.
\medskip
$$\eqalign {&a'\circ \phi = a, \quad \phi [X, Y]_A= [\phi (X),
\phi (Y)]_{A'}\cr
&a_*\circ \phi^*= a'_*,\quad \phi^* [\alpha',
\beta']_{(A')^*}= [\phi^* (\alpha'),
\phi^* (\beta')]_{A^*}\cr
& {\hbox {  for all sections }} X,Y:B\rightarrow A,\, \alpha', \beta':
B\rightarrow (A')^*.\cr}$$}
\medskip
The following property \cite{MX1}  is basic.
\medskip
{\bf Proposition 3.1.4} \hskip 5 pt {\it Let $Y, Y'$ be Poisson groupoids
over 
$B$  with tangent Lie bialgebroids $(A,
A^*)$ and $(A', (A')^*)$.
If $\mu: Y\longrightarrow Y'$ is a base preserving morphism of Poisson
groupoids, then 
$$ A(\mu): A\longrightarrow A'$$
is a (base preserving) morphism of Lie bialgebroids.}

\bigskip
\bigskip
{\bf 3.2. The dual of a dynamical Poisson groupoid}
\bigskip

Throughout this subsection, we shall equip $X= U\times G\times U$
with the trivial Lie groupoid structure and we assume that the subgroup
$H\subset G$ is connected and simply connected. We begin with a
description of
the tangent Lie bialgebroid of a
Poisson groupoid $(X, \{\, , \,\})$, a special instance of which is
described 
in
\cite{BKS}. 
\smallskip
{\bf Proposition 3.2.1}
{\it Let $(X,\{\, , \, \})$ be a Poisson groupoid  with
Poisson bracket as in Proposition 2.2.3.}
{\it Then the Lie bialgebroid tangent to $X$ is (isomorphic to) the pair
$\big(
U\times 
\fh^*\times \fg, \, U\times \fh\times \fg^*\big)$ with
anchor maps}
$$
\align &a: U\times \fh^*\times \fg\rightarrow U\times
\fh^*: (q, 
\lambda, X)\mapsto (q,\lambda) \\
&a_*: U\times \fh\times \fg^*\rightarrow U\times \fh^*:
(q, Z, B)\mapsto (q, -K(q)Z + A^* (q) B)\
\endalign
$$
{\it and  Lie brackets of sections $[\, , \, ]$, $[\, , \, ]_*$ given by
the following expressions}
$$
\align &[(\lambda, X), (\lambda', X')](q)= (d\lambda' \cdot
\lambda  - 
d\lambda \cdot \lambda', d X'\cdot \lambda - d X\cdot \lambda'
+ [X(q),X'(q)]_g) \\
& \\
&<[(Z,B), (Z',B')]_*, (\Lambda, Y)> (q)= \\
&<-dZ'(K(q)Z-A^*(q)\, B) + dZ (K(q) Z' -A^*(q)\, B'), \Lambda>\\
&-<Z', dK (q) (\Lambda) Z>\\
&-<B, \delta_1 P (\Lambda) B'> - <B', dA (q) (\Lambda) \, Z> + <B, dA (q)
(\Lambda) Z'>\\
&+<dB (K(q) Z' - A^* (q) B') - dB' (K(q) Z - A^* (q) B), Y>\\
&+ <ad^*_ {A(q) Z} B' - ad^*_{A(q) Z'} B, Y>\\
&+<B, \partial P (Y) B'>,\\
\endalign
$$

{\it where $\lambda, \lambda': U \rightarrow \fh^*$;$\, X,X':
U\rightarrow \fg$; $\,Z,Z': U\rightarrow \fh$;$\,
B, B':U\rightarrow \fg^*$ are smooth maps, $\, \Lambda\in \fh^*,
Y\in \fg$,  all differentials and sections are evaluated at $q$, and the
partial 
derivatives of the groupoid $1- $ cocycle $P$ (see Thm. 2.2.5) are evaluated
at $(q,1,q).$ }
\bigskip
{ Proof: } The Lie algebroid $A(X)$ is well known (see, for example, [M1]).
Although for the dynamical coboundary case an algebraic description of
the Lie algebroid dual $A(X)^*$ was
given in \cite{BKS}, it was not derived there from the Poisson groupoid
using
Weinstein's 
coisotropic calculus. So we shall briefly indicate the steps of the
calculation.
\smallskip
The unit section of $X$ is given by $\epsilon :U\rightarrow U\times G
\times U: q\mapsto (q,1,q)$. Therefore, $\gamma\in N(\epsilon (U))_q$
if and only if $\gamma= (-Z, B, Z)$, for some  $Z\in \fh$ and $B\in \fg^*.$
\smallskip
Let $\alpha,\alpha': U \rightarrow N(\epsilon (U))$ be two sections
written as 
$\alpha (q)= (-Z(q), B(q), Z(q))$ and $\alpha' (q)= (-Z'(q), B'(q),
Z'(q))$. Set $\omega (p,x,q)= (-Z(q), T^*_x l_{x^{-1}} B(q), Z(q))$,
and similarly for $\omega'$. By (3.1.2),  it suffices
\smallskip
 1) to compute
$$[\alpha, \alpha']_{N(\epsilon (U))} (q)= [\omega, \omega'](q,1,q),$$
where the {\sl rhs} is given by (3.1.1) with Hamiltonian operator
$$
\align &\Pi^\# (p,x,q) (Z_1, B, Z_2)= (-K(p) Z_1 -A^*
(p)
T^*_1 r_x B, \\
&\hskip 15 pt T_1 r_x A (p) Z_1 + T_1 l_x A (q) Z_2+ T_1 r_x
P(p,x,q) T^*_1 r_x B,\\
&\hskip 30 pt K(q) Z_2 - A^* (q) T^*_1 l_x B),\\
\endalign
$$
and 
\smallskip
2) to choose an identification of $N(\epsilon (U))$ with $U\times
\fh\times \fg^*.$
\smallskip
The computation for (1)  is rather standard (although somewhat lengthy) and
may be performed with
the help of
$$<L_{\Pi^\# \omega} \omega', (\Lambda,  X^l, \Lambda') >=
d<\omega', (\Lambda,  X^l, \Lambda') > \cdot \Pi^\#\omega +
<\omega', [(\Lambda,X^l, \Lambda '), \Pi^\# \omega]>,$$
where $\Lambda, \Lambda'\in \fh^*, X\in \fg,$ and  $X^l$ is the left
invariant vector field on $G$ with $X^l (1)= X.$
\smallskip
As for 2) the natural identification to make is given by $\iota_-:
N(\epsilon (U))_q\rightarrow \fh\times \fg^*:
(-Z, B, 
Z)\mapsto (Z,B)$. Setting $[(Z,B), (Z',B')]_* (q)= \iota_- [\omega,
\omega'] (q,1,q)$ then gives the stated formula.$\pf$
\bigskip
{\bf Remark 3.2.2}
The coboundary dynamical case in \cite{BKS} corresponds to the choice:
$\chi 
(h,q)= h$, and the groupoid $1-$ cocycle
$$P(p,x,q)= -R 
(p) + Ad_x 
R (q) Ad^*_x.$$

Note that $A(q)\,Z= A_\chi (q) \,Z=\iota\, Z$ is constant in this case,
while the induced algebroid $1 - $ cocycle
defined as
$$\eqalign {&P_* : U\times \fh^*\times \fg\longrightarrow
L(\fg^*, 
\fg)\cr
& (q,\Lambda, Y)\mapsto (\partial P\cdot Y + \delta_2 P\cdot
\Lambda)(q,1,q)= (\partial P\cdot Y - \delta_1 P\cdot \Lambda)(q,1,q) \cr}$$
which appears in the bracket of Prop. 3.2.1 is given by
$$P_* (q,\Lambda, Y)= dR(q) \Lambda +
R(q) ad^*_Y + ad_Y R (q) .$$
\bigskip
\bigskip
The proposition which follows is a special instance of the functorial
relationship between Poisson groupoids and Lie bialgebroids.
\bigskip
{\bf Proposition 3.2.3}
{\it Let $(X,\{\, , \, \})$ be a Poisson groupoid
with Poisson bracket as in Prop. 2.2.3 and let
$$\eqalign {&I: H\times U\longrightarrow U\times G\times
U\cr
&(h,p)\mapsto (Ad^*_{h^{-1}} p, \chi (h,p), p)\cr}$$
be a groupoid morphism. If $A(I)$ is a morphism of Lie bialgebroids, then
$I$ is also a
Poisson map.  Hence $I$ is a Poisson groupoid morphism.}

\medskip
Proof. We have to show that the three conditions of Prop. 2.2.4 are
satisfied. From the definition of $I$, it is clear that
the induced morphism $A(I): U\times \fh\longrightarrow U\times
\fg\times \fh^*$
is given by 
$$(q,Z)\mapsto (q, A_\chi (q) (Z), ad^*_Z q),$$ so its
dual map $A(I)^*: U\times g^*\times \fh\longrightarrow U\times
\fh^*$ is of the form
$$ (q, B, Z) \mapsto (q, A_\chi^*(q) B - ad^*_Z q),$$
where $A_\chi (q) (Z)= A(\chi) (q,Z).$ Thus $A(I)^*$ preserves the
anchor maps if and only if
$$K(q) (Z)= ad^*_Z q {\hbox { and }} A^* (q)= A_\chi^* (q),$$
while a direct calculation shows that
$$A(I)^* [(Z,B), (Z',B')]_* =  [A(I)^* (Z,B),
A(I)^*(Z',B')]_{U\times \fh^*}\eqno (\star)$$
is equivalent to
$$
\align
&<B, \partial P (A_\chi ({\bold Z})) B'- \delta_1P (ad^*_{\bold Z}
q) B'>\\
&\hskip 5 pt = <B', dA_\chi (A_\chi^* (B)) {\bold Z}>- <B, dA_\chi
(A_\chi^* (B')) {\bold Z}>.\hskip 40 pt (\star\star)\\
\endalign
$$ 
Clearly, the anchor conditions are precisely the conditions (a) and (b) of
Prop.2.2.4. Therefore, it remains to
verify condition (c). To this end, let $\rho: G\longrightarrow Aut( L(\fg^*,
\fg))$ be the adjoint action and set
$$
\align 
&\rho_\chi: H\times U\longrightarrow Aut (L(\fg^*, \fg))\\
&\hskip 15 pt (h,p)\mapsto \rho (\chi (h,p)).\\
\endalign
$$
Then $\rho_\chi$ is a groupoid representation. We begin by showing that both
sides of condition 2.2.4 (c) are groupoid $1-$cocycles for $\rho_\chi$ (see
Def. 2.1.6).
\smallskip 
 That the left hand side
$P\circ I: H\times U\rightarrow L(\fg^*, \fg)$
is a groupoid $1-$cocycle for $\rho_\chi$ is immediate from the fact
that $P: U\times G\times U\rightarrow L(\fg^*, \fg)$ is such a
cocycle for $\rho.$
\smallskip 
As for the right hand side, we have to show that the map
$\Sigma: H\times U\rightarrow L(\fg^*, \fg)$
defined by
$$<\alpha, \Sigma (h,p) \beta>= <\lambda_\alpha, Z_\beta>- <\lambda_\beta,
Z_\alpha> - <p, [Z_\alpha, Z_\beta]>$$
satisfies
$$\Sigma (hk,p)= \Sigma (h,Ad^*_{k^{-1}} p)+ Ad_{\chi (h,
Ad^*_{k^{-1}} p)} \Sigma (k,p) Ad^*_{\chi (h, Ad^*_{k^{-1}} p)}.$$
But this follows by a direct calculation which makes successive use of the
following three identities
$$
\align
<Z_\alpha (hk,p), \lambda>&= <Z_\alpha (h, Ad^*_{k^{-1}} p),
Ad^*_{k^{-1}} \lambda> + <Z_{Ad^*_{\chi (h, Ad^*_{k^{-1}}p)}
\alpha}(k,p),\lambda>\\
<\lambda_\alpha (hk,p), Z>&= <Ad^*_k \lambda_\alpha (h, Ad^*_{k^{-1}} p),
Z>\\
<\lambda_\alpha (hk,p),Z>&=<\lambda_{Ad^*_{\chi (h,
Ad^*_{k^{-1}}p)}\alpha} (k,p), Z>+ <p, [Ad_{k^{-1}} Z_\alpha (h,
Ad^*_{k^{-1}} p), Z]>\\
\endalign
$$
Now, since $H$ is simply connected, by Prop. 7.3 of \cite{X} the two
groupoid
$1-$cocycles $P\circ I$ and $\Sigma$ coincide iff their
induced algebroid  cocycles are the same. But the latter is equivalent to
$(\star\star)$ above. This concludes the proof.$\pf$
\bigskip
\bigskip
Recall that a Lie algebroid $A$ over a (connected) base $B$ is said to be
transitive iff its anchor map $a: A\rightarrow TB$ is a surjective
submersion.  In this case
the Kernel $Ker \, a$ of $a$ is a Lie algebra bundle \cite{M1} called the
adjoint bundle of $A$ whose
fibers are called the vertex (or isotropy) Lie algebras $\frak k$ of $A$.
If $A$ is a transitive Lie algebroid over a contractible base $B$, it is
shown in \cite{M1} that $A$
isomorphic to the trivial Lie algebroid $TB\oplus (B\times {\frak k})$
(Whitney
sum), where $\frak k$ is the typical fiber; in particular $A$ integrates to
a
global Lie groupoid isomorphic to
$B\times K\times B$ where $K$ is the connected and simply connected Lie
group
with $Lie (K)= {\frak k}.$
\medskip

With these facts, we immediately obtain a description of the dual
of a dynamical Poisson groupoid.
\smallskip
 Let $X$ be a dynamical Poisson groupoid as in Def. 2.2.6 over the
 contractible base $U$ with embedding of the Hamiltonian unit given by
 $$I(h,q)= (Ad^*_{h^{-1}} q, h, q).$$
 Let $\iota :\fh\rightarrow
\fg$ be the inclusion map.
\medskip
{\bf Theorem 3.2.4} (Duality)
\smallskip
{\it The dual Poisson groupoid $X^*$ of
$X$ is isomorphic to the  Poisson groupoid $
(U\times 
G'\times U, \{\, , \, \}_{U \times G'\times U})$
where $G'$ is the connected and simply connected Lie group whose Lie algebra
is the vector space $\frak k:=\{(Z,A)\in \fh\times \fg^*\, \mid \, ad^*_Z
(q_0)= \iota^*A\}$ 
for some $q_0\in U$, equipped with the Lie bracket
$$
\align &[(Z, A), (Z',A')]= \big(-[Z,Z'] -
<A,\delta_1 
P ( \cdot ) A'>, \\
&\hskip 15 pt  ad^*_{\iota Z} A' - ad^*_{\iota Z'} A + <A,
\partial P ( \cdot ) \, A'>\big),\\
\endalign
$$
and the Poisson bracket $\{\, , \,\}_{U\times
G'\times U}$ is given by Theorem 2.2.5 for a (unique) Poisson groupoid
morphism $I':H\times U\rightarrow U\times
G'\times U$ and a (unique)  groupoid $1 - $ cocycle
$$P': U\times G'\times U\longrightarrow L({\frak k}^*, {\frak k}) $$
for the adjoint action of ${\frak k}.$}

\medskip
{ Proof. } For the first part,  observe that the anchor map of $A(X)^*$
$$a_* (q, Z, A)= (q, -ad^*_Z q +\iota^* A)$$
is a surjective submersion since $\iota$ is injective. Thus
$A(X)^*$  is 
transitive and therefore, by MacKenzie's theorem, it is isomorphic to the
trivial Lie algebroid $A'=\fh^*\times \fh^*\times (Ker a_*)_{q_0} .$
Now the 
fiber $(Ker \, a_*)_{q_0}$ is the vector space $\frak k$
equipped with the Lie bracket  given by the restriction of the bracket of
sections of 
$A(X)^*$ of Prop. 3.2.1 (with $K(q) Z= ad^*_Z q$ and $A(q)Z=\iota Z$). Hence
the claim.
\smallskip
For the second part, let
$$\tau: A(X)^*\longrightarrow A'$$
be the (base preserving) trivializing isomorphism of MacKenzie's theorem,
and denote by
$$T:X^*\longrightarrow U\times G'\times U$$
 the unique groupoid isomorphism such that $A(T)= \tau.$
We may thus transport the Poisson groupoid structure of $X^*$ to
$U\times 
G'\times U$ by setting
$$\{f,g\}_{U\times G'\times U}= \{f\circ T,\, g\circ T\}_{X^*}
\circ T^{-1}.$$
We now show that there is a (base preserving) Poisson groupoid
morphism 
$$I': H\times U\longrightarrow U\times G'\times U,$$
where $H\times U$ is the Hamiltonian unit.
 \smallskip
Consider the Poisson groupoid morphism (this is the anchor map of $X$)
$$J: U\times G\times U\longrightarrow U\times U:
(p,x,q)\mapsto (p,q)$$ where $U\times U$ is the coarse groupoid
of Example 3.1.2.  Its  induced Lie bialgebroid morphism
$$A(J): U\times \fg\times \fh^*\longrightarrow U\times \fh^*$$
is of course just the anchor map $a$ of $A(X)$.
By the lifting property of Lie algebroid morphisms, and Prop.
3.2.3 above, the dual morphism
$$A(J)^*=a^*: U\times \fh\longrightarrow U\times \fg^*\times \fh$$
may be lifted uniquely to a (base preserving) Poisson groupoid morphism
$$J^*: H\times
U\longrightarrow X^*.$$
Thus $I'=T\circ J^*$ is the sought-for Poisson groupoid morphism.  The
uniqueness of  $P'$ now follows
from the uniqueness of the Poisson structure of a (suitably simply
connected) Poisson groupoid with prescribed tangent Lie bialgebroid
\cite{MX2}.  Hence the claim.
$\pf$
\bigskip 
{\bf Caveat} 
\smallskip
We shall see in section 5 that, for $\fh\not=0$, even when $X$ is coboundary
with {\bf constant } $r- $ matrix, the vertex group $G'$ is different from
the Poisson Lie group dual to $G$ equipped with the Sklyanin
bracket $\{\, , \, \}_{(R,-R)}.$
\bigskip
\bigskip

We close this subsection with a description of natural Poisson quotients
associated with Thm 3.2.4. We now assume that the contractible set $U$
contains $0.$
\bigskip
Let $X$ be as in Thm 2.2.5 with the map $H\longrightarrow G: h\mapsto \chi
(h,0)$ one to one. Consider the restriction of the left
Hamiltonian action
$$\phi^- : H\times X\longrightarrow X: (h,(p,x,q))\mapsto
(Ad^*_{h^{-1}} p, \chi (h,p) x, q)$$
to $\alpha^{-1} (0)= \{0\}\times G\times U$:
$$\phi^- :H\times G\times U\longrightarrow G\times U:
(h,x,q)\mapsto (\chi (h,0) x, q).$$
Let $\pi : G\times U \rightarrow G/H\times U: (x,p)\mapsto
({\overline x}, p)$ be the canonical projection.
\bigskip
{\bf Proposition 3.2.5}  (Hamiltonian reduction)
\smallskip
{\it The Poisson bracket $\{\, , \, \}_{red.}$ of the reduced space
$\alpha^{-1} (0)/H\simeq G/H\times
U$ vanishes at $({\overline 1}, 0).$ Its linearization at
$({\overline 1}, 0)$ coincides with the
vertex Lie algebra of $A(X)^*$ at $0.$ }
\bigskip
{ Proof.} 
We have to calculate the Poisson bracket $\{f, g\}_X (0, x, q)$ of two
 functions $f,g\in C^\infty (X)$ whose restriction to $\{0\}\times G\times
 U$ is $H - $ invariant.
\smallskip
Since $H$ is connected, the restriction of $f$ to $\{0\}\times G\times
U$ is $H - $ invariant if
and only 
if 
$$A_\chi^* (0) Df (0,x,q)= 0 {\hbox { for all }} x\in G, q\in U.$$
Thus (see Thm 2.2.5)
$$\eqalign {&\{f, g\}_X (0, x, q)= - <q, [\delta_2 f, \delta_2 g]> -
<A_\chi (q) \delta_2 f, D'g>\cr
&\hskip 20 pt + <A_\chi (q)\delta_2 g, D'f> + <Df, P(0,x,q) Dg>.\cr}$$
Now, $P(0,1,0)=0$ and $A_\chi^* (0) D'f (0,1,0)= A_\chi^* (0)
Df(0,1,0)= 0$.  Therefore,
$$\{f, g\}_X (0,1,0)= 0.$$
Set $Z= \delta_2 f , Z'=\delta_2 g , A= Df, A'= Dg$, all evaluated at
$(0,1,0).$ A direct calculation then gives
$$
\align
 &d \{f,g\}_X (0, 1, 0) (0, Y, \lambda)= -<\lambda, [Z,Z']> -
<dA_\chi (0)(\lambda)\, Z, A'>\\
&\hskip 10 pt +<dA_\chi (0)
(\lambda) Z', A>+ <ad^*_{A_\chi (0)Z} A', Y>\\
&\hskip 5 pt - <ad^*_{A_\chi (0) Z'} A, Y>
+<A,(\partial P ( X)+ \delta_2 P (\lambda)) A'>.\\
\endalign
$$
To conclude, observe that this coincides with the restriction of the Lie
bracket of Prop. 3.2.1 (with $K(q) Z= ad^*_Z q$ and $A=A_\chi$) to the
kernel $Ker \, a_*$ at $0. \pf$
\medskip
Prop. 3.2.5 provides in some sense an indirect Poisson integration of the
vertex Lie
algebra ${(Ker \, a_*)} _{0}$ of $A(X)^*$ by the natural
quotient space 
$G/H\times U.$
For $X$ dynamical, combining Prop. 3.2.5 with Thm 3.2.4 then gives a reduced
vertex
diagram reminiscent of the Poisson Lie group duality of Drinfeld.
\medskip
Let $X=U\times G\times U$ be a dynamical Poisson groupoid as in Def.
2.2.6 with dual Poisson groupoid $X^*\simeq
U\times  G'\times U.$ Assume that the map  $H\rightarrow G':h\mapsto
\chi' (h,0)$ is one to
one.  Denote the units of $G$ and $G'$ by $1$ and
equip both spaces $G/H\times U ,\, G'/H\times
U$ with the Poisson brackets obtained via Poisson
reduction.
\medskip
Let $\fh^\perp\subset g^*$ be the annihilator of $\fh$.
\medskip
{\bf Theorem  3.2.6} (Reduced duality diagram)
\smallskip
{\it We have the diagram
$$
\align 
&\hskip 20 pt X^*\hskip 180 pt X \\
&\hskip 40 pt red.\searrow\hskip 100  pt \swarrow red.\\
&\hskip 50 pt G'/H\times U\hskip 50 pt G/H\times U\\
&\hskip 100 pt \searrow \hskip 20 pt \swarrow\\
&\hskip 100 pt lin.\, at \, ({\overline 1}, 0)\\
&\hskip 100 pt  \swarrow\hskip 20 pt  \searrow \\
&A(X)^* \supset ({Ker\, a_*})_{0}= \fh\times \fh^\perp \hskip 40 pt
g= ({Ker\, 
a})_0 \subset A(X)\\
\endalign
$$
\smallskip
In case $H$ is reduced to its unit, this diagram reduces to
that of Drinfeld's duality for Poisson Lie groups. }$\pf$
\bigskip
\bigskip
{\bf 4. An explicit case study of duality}
\bigskip 
In \cite{EV}, Etingof and Varchenko obtained, among other things, a
classification of solutions of the
(CDYBE) for pairs $(\fg, \fh)$ of Lie algebras, where $\fg$ is simple,
and $\fh \subset \fg$ is a Cartan subalgebra.
Our purpose in this section is to give an explicit study of duality
for the corresponding class of coboundary dynamical Poisson groupoids.
\smallskip
We begin by recalling the general form of these dynamical $r-$ matrices.
\smallskip
First, let us fix some notation.  Let $\fg$ be a complex simple Lie
algebra with Killing form $(\, , \,)$, $\fh \subset \fg$ a Cartan
subalgebra, 
and
$\fg=\fh \oplus \sum_{\alpha \in\Delta} \fg_{\alpha}$
the root space decomposition.
We let $\Delta^s$ be a fixed simple system of roots and denote by
$\Delta^\pm$ 
the corresponding positive/negative
system.  For any positive root $\alpha \in \Delta^+$, we  choose
root vectors $e_{\alpha} \in \fg_{\alpha}$ and
$e_{-\alpha} \in \fg_{-\alpha}$ which are dual with respect to
$(\, , \,)$ so that $[e_{\alpha} , e_{-\alpha}]=h_{\alpha}$.  We
also fix an orthonormal basis $(x_i)_{1\le i\le {rank (\fg)}}$ of $\fh$.
Lastly, for a subset of simple roots $\Gamma\subset \Delta^s$, we shall
denote the 
root span of $\Gamma$ by
$<\Gamma>\subset \Delta$ and set ${\overline
\Gamma}^\pm= \Delta^\pm \setminus <\Gamma>^\pm.$
\smallskip
For any subset $\Gamma \subset \Delta^s$,  we give
 the $ ad_\fh-$ invariant solutions of $(mDYBE)$
(see eqs. (2.1.1), (2.1.2)) associated with the triple
$(\fg, \fh, \Gamma)$
as (cf.\cite{EV}):
\smallskip

$$R(q)B= \sum_{i,j} C_{ij} (q)<x_j,B> x_i + \sum_{\alpha\in \Delta}
\phi_\alpha (q)<e_{-\alpha},B> e_{\alpha} \eqno (4.1)$$
where 
\smallskip
$$\eqalign {&\phi_\alpha (q) = {1\over 2}\,\, {\hbox { for }} \alpha\in
{\overline 
\Gamma}^+,\quad\phi_\alpha (q) = - {1\over 2} \,\,{\hbox { for }} \alpha\in
{\overline 
\Gamma}^-\cr
&\phi_\alpha (q)= {1\over 2} coth ({(\alpha, q-\mu)\over 2})\,\, {\hbox {
for }} \alpha\in <\Gamma>,\cr}$$
and where $\sum_{i,j} C_{ij} dq^i\otimes dq^j$ is any closed
meromorphic $2-$ form on $\fh^*$ and $\mu\in \fh^*$ is arbitrary.
\medskip
We shall denote by $U$ the domain of analyticity of $R$ and let $G$
be the connected and simply-connected Lie group with $Lie(G)=\fg.$
Note that $U$ is trivially $Ad_H^*-$ invariant as $H$ is abelian, hence
we can consider the coboundary dynamical Poisson groupoid
$X=U \times G\times U$ associated with $R$.  Our immediate goal is
to construct an explicit trivialization of the dynamical Lie
algebroid $A(X)^* \simeq U \times \fh^* \times \fg^*$.
Note that, as $U$ is not contractible, this is not guaranteed
by MacKenzie's theorem.
In what follows, we shall make the identification $\fg^* \simeq \fg$
using the Killing form $(\, , \,)$.  Then we have
$\fh^\perp \simeq \frak n:=\sum_{\alpha \in \Delta}\fg_{\alpha}$,
and the Lie bracket between the sections of the dynamical Lie algebroid
takes the form 
$$
\align 
[(Z, B), & (Z', B')]_*(q) \\
& = (dZ'(q)i^*B(q)- dZ(q)i^*B'(q)
  - (dR (q)(.)B(q), B'(q)),  \\
&\hskip 15pt  -dB(q)i^*B'(q)+dB'(q)i^*B(q)   \\
&\hskip 15pt  -[R(q)(B(q))+Z(q), B'(q)] - [B(q), R(q)(B'(q))+Z'(q)]),
\qquad
\thetag{4.2}
\endalign
$$
\smallskip
We shall begin our construction with a description of the vertex
Lie algebra ${\Cal V}_q= {(Ker\,
a_*)}_q=\fh\times \fh^\perp\simeq \fh\times {\frak n}$ of ${A^*}$ at $q\in
{U}.$
To do so, let us introduce the following Lie subalgebras of $\fg$
associated with $\Gamma \subset \Delta^s$:
\smallskip
 $$\eqalign {&\fh_\Gamma:= <(h_\gamma)_{\gamma\in \Gamma}>_{\bold
C},\cr
 &\fh_\Gamma^\perp := {\hbox { the orthogonal complement of }} \fh_\Gamma
 {\hbox { in }} \fh
 {\hbox { w.r.t. }} (\, , \,)_{\mid_ {\fh\times \fh}},\cr
 &{\frak l}_\Gamma:= \fh_\Gamma \oplus <(e_\alpha)_{\alpha\in
<\Gamma>}>_{\bold C}
 {\hbox { the Levi factor }},\cr
  &{\overline 
{\frak n}}^\pm_\Gamma:= <(e_\alpha)_{\alpha\in {\overline
\Gamma}^\pm}>_{\bold C}
{\hbox { 
the nilpotent radicals }}.\cr}$$
Clearly $[\fh_\Gamma, {\overline
{\frak n}}^\pm_\Gamma]\subset {\overline
{\frak n}}^\pm_\Gamma$ and $[{\frak l}_\Gamma, {\overline
{\frak n}}^\pm_\Gamma]\subset {\overline
{\frak n}}^\pm_\Gamma.$
\medskip
If $\fg_1, \fg_2$ are two Lie algebras, we denote by $\fg_1\ominus \fg_2$
the
vector space $\fg_1\oplus \fg_2$ equipped with the Lie bracket $[x_1+x_2,
y_1+y_2]= [x_1, y_1] - [x_2, y_2].$
 Let ${\frak I}_\Gamma= \fh_\Gamma ^\perp \ltimes \big({\overline
{\frak n}}^+_\Gamma \ominus {\overline
{\frak n}}^-_\Gamma\big)$ be the semidirect product Lie algebra where
$\fh_\Gamma^\perp$ acts on each summand of the anti-direct sum by the
adjoint
action of $\fg.$ Set
$$\fg':={\frak l}_\Gamma \ltimes {\frak I}_\Gamma$$
where ${\frak l}_\Gamma$ also acts on ${\frak I}_\Gamma$ by the adjoint
action of $\fg$.
\bigskip
{\bf Proposition 4.1} {\sl  Let $q\in {U}$. Then the map $\psi (q) : {\Cal
V}_q\longrightarrow \fg'$ defined by
$$\eqalign { &\psi (q) (0,e_\alpha) = {-1\over 2 sinh \big({ (\alpha,
q-\mu)\over
2}\big)}  
e_\alpha\quad {\hbox { for }} \alpha\in <\Gamma>,\cr
&\psi (q) (0,e_\alpha)= -e^{\mp {1\over 2} (\alpha, q-\mu)} e_\alpha\, \quad
{\hbox { 
for }} \alpha\in  {\overline \Gamma}^\pm,\cr
&\psi (q) (Z,0) = -Z \quad {\hbox { for all }} Z\in \fh,\cr}$$
is an isomorphism of Lie algebras.}
\medskip
{ Proof.} The Lie bracket of ${\Cal V}_q= \fh\times {\frak n}$
can be calculated from Eqn. (4.2) and
 we have
 $$
\align
 & [(Z, n), (Z', n')]_*
= \big(-(dR (q)(.)n,\, n'), -[R(q)n+Z, n']\\
&\hskip 40 pt -[n, R(q)n'+Z']\big).\\
\endalign
$$
Writing $n=\sum_{\alpha\in \Delta} n^\alpha e_\alpha$ and similarly for
$n'$, we have
$$
\align
(dR (q) \, ({\bold \Lambda})\, n,\, n')&= \sum_{\alpha\in <\Gamma>}
d\phi_\alpha (q)\, ({\bold \Lambda})\, n^\alpha {n'}^{-\alpha}\\
&=
\sum_{\alpha\in <\Gamma>} ({1\over 4} -\phi_\alpha (q)^2) (\alpha, {\bold
\Lambda}) \, n^\alpha {n'}^{-\alpha}.\\
\endalign
$$
A direct calculation then gives the following Lie bracket relations
:
$$\eqalign { &[(0, e_\alpha),
(0, e_\beta)]_*=(0, -(\phi_\alpha (q) + \phi_\beta (q)) [e_\alpha, e_\beta])
{
\hbox { 
for  }}\alpha\in <\Gamma>, \beta\in \Delta,\, \alpha+\beta\not=0,\cr
 &[(0, e_\alpha),(0, e_{-\alpha})]_*=(- ({1\over 4} -
\phi_\alpha (q)^2) [e_\alpha, e_{-\alpha}], 0)\quad {\hbox { for }}\alpha\in
<\Gamma>,\cr
&[(0, e_\alpha), (0, e_\beta)]_*= (0,0)\quad {\hbox { for }}\alpha\in
{\overline
\Gamma}^+, \beta\in {\overline \Gamma}^-, \cr
 &[(0, e_\alpha),(0, e_\beta)]_*=
(0, -[e_\alpha, e_\beta])\quad {\hbox { for }} \alpha,\beta\in {\overline
\Gamma}^+\cr
&[(0, e_\alpha),(0, e_\beta)]_*=
(0,+[e_\alpha, e_\beta])\quad {\hbox { for }} \alpha,\beta\in {\overline
\Gamma}^-, \cr
& [(Z,0),(0, e_\alpha)]_*=(0, -\alpha (Z)
e_\alpha)\quad {\hbox { for all }} \alpha\in \Delta, Z\in \fh\cr
&[(Z,0),(Z',0)]_*=(0, 0),\quad Z,Z'\in \fh.\cr}$$
\smallskip
where the bracket $[\, , \, ]$ on the r.h.s.  is that of $\fg.$
After rescaling the basis of ${\frak n}$ by setting
$$ E_\alpha (q)= 2 sinh \big( {(\alpha, q-\mu)\over 2}\big)
e_\alpha, \, \alpha\in <\Gamma>; \quad E_\alpha (q)= e^{\pm {1\over 2}
(\alpha, q-\mu)} e_\alpha, \, \alpha\in {\overline \Gamma} ^\pm,$$
the above relations yield
$$\eqalign { &[(0, E_\alpha (q)),
(0, E_\beta (q))]_*=(0,- N_{\alpha,\beta} E_{\alpha+\beta} (q) ){
\hbox { 
for  }}\alpha\in <\Gamma>, \beta\in \Delta,\, \alpha+\beta\not=0,\cr
 &[(0, E_\alpha (q)),(0, E_{-\alpha}(q))]_*=(- [e_\alpha, -e_{-\alpha}],
0)\quad {\hbox { for }}\alpha\in
<\Gamma>,\cr
&[(0, E_\alpha (q)),(0, E_\beta(q))]_*= (0,0)\quad {\hbox { for }}\alpha\in
{\overline
\Gamma}^+, \beta\in {\overline \Gamma}^-, \cr
 &[(0, E_\alpha (q)),(0, E_\beta (q))]_*=
(0,-N_{\alpha, \beta} E_{\alpha+\beta}(q))\quad {\hbox { for }}
\alpha,\beta\in
{\overline
\Gamma}^+\cr
&[(0, E_\alpha (q)),(0, E_\beta(q))]_*=
(0, + N_{\alpha,  \beta} E_{\alpha+\beta}(q))\quad {\hbox { for }}
\alpha,\beta\in
{\overline
\Gamma}^-, \cr
& [(Z,0), (0, E_\alpha (q))]_*= (0,-\alpha (Z)
E_\alpha (q))\quad {\hbox { for all }} \alpha\in \Delta, Z\in \fh,\cr
&[(Z,0), (Z',0)]_*=(0, 0)\quad Z,Z'\in \fh ,\cr}$$
where $N_{\alpha, \beta}$ are the structure constants of $\fg.$
\smallskip
We shall check the first Lie bracket above  with $\beta\in {\overline
\Gamma}^-$; the others are similar.
\smallskip
For a root $\gamma\in \Delta,\,$ set $x_\gamma=(\gamma, q-\mu).$ We have
$\phi_\alpha (q)+ \phi_\beta (q) ={1\over 2} coth ({x_\alpha\over
2})-{1\over 
2}= {1\over {e^{x_\alpha}-1}}$,
$E_\alpha (q)= e^{x_\alpha\over 2} (1- e^{-x_\alpha}) e_\alpha,$ and
$E_\beta (q)= e^{-x_\beta\over 2} e_\beta.$
Thus,
$$[(0, E_\alpha (q)), (0, E_\beta (q))]_*=(0,- {e^{x_\alpha\over 2}
(1-e^{-x_\alpha}) e^{-x_\beta\over 2} \over (e^{x_\alpha}-1)} N_{\alpha,
\beta} e_{\alpha+\beta}).$$
Now, if $\alpha+\beta$ is a root, then it belongs to ${\overline
{\Gamma}}^-$;
thus $e_{\alpha+\beta}= e^{x_\alpha +x_\beta \over 2}
E_{\alpha+\beta} (q)$, and this immediately gives the assertion.
\smallskip

 The Lie bracket relations above show that the structure constants of
${\Cal V}_q$ in
the basis
$((x_i,0)\, 1 \le i\le rank(\fg)$;
$(0, E_\alpha (q)), \alpha\in
\Delta)$ are opposite to those of $\fg'.$ Therefore the map $\psi: {\Cal
V}_q\longrightarrow \fg '$ defined by
$(Z,0)\mapsto -Z$ and $(0, E_\alpha (q))\mapsto -e_\alpha$ is an
isomorphism of Lie algebras. Hence
the claim.$\pf$
\medskip
{\bf Corollary 4.2}  {\it The map $\widetilde{\psi}:U \times
\fg'\longrightarrow Ker\, a_*
:(q,\xi)\longrightarrow (q,-\Pi_{\fh}\xi, \psi(q)^{-1}(\Pi_{\frak n}\xi))$
is an isomorphism between the trivial Lie algebra bundle $U\times \fg'$
and the adjoint bundle $Ker\, a_*$.  Here, $\Pi_{\fh}$ and $\Pi_{\frak n}$
are
the projections relative to the direct sum decomposition
$\fg=\fh\oplus\frak n.\quad \pf$}
\bigskip
Let us briefly comment on the vertex isomorphism of Prop. 4.1.
If $\Gamma\subset \Delta^s$ is
the empty set, the vertex Lie algebra
${\Cal V}_q,\, q\in U,$ of ${A(X)^*}$ is isomorphic to $\fh\ltimes ({\frak
n}
^+\ominus {\frak 
n}^-)$, which is reminiscent of (although not identical to) the  Lie
algebra dual of $\fg$ equipped with the standard constant  $r-$ matrix
(see also example
5.1.8).
If $\Gamma= \Delta^s$, we have ${\Cal V}_q\simeq\fg'={\frak l}_\Delta= \fg.$
For a general subset $\Gamma\subset \Delta^s$ the vertex Lie algebra ${\Cal
V}_q\simeq\fg'$ 
is seen to naturally intertwine the Levi factor ${\frak l}_\Gamma$ with the
summand ${\frak I}_\Gamma$ which is again
reminiscent of the Lie algebra dual of $\fg$  equipped with the
standard constant $r-$ matrix.
\bigskip
Our next step is to construct a flat connection
$\theta_{*}: TU \simeq U\times \fh^*\longrightarrow U\times\fh\times\fg^*$
satifying the condition
$[\,\theta_{*}(\lambda),
\widetilde{\psi}(\xi)\,]_{*}=\widetilde{\psi}(d\xi\cdot\lambda)$
for $\lambda: U\longrightarrow \fh^*$ and $\xi: U\longrightarrow\fg'.$
To simplify notation,  we shall identify the elements
$(Z,n)\in \fh\times {\frak n}\simeq {\Cal V}_q$ of the vertex Lie algebra
with $Z+n\in \fg$ from now onwards.
\medskip
 Let $C^\#: {U}\rightarrow L(\fh^*, \fh)$ be the map defined by
$C^\# (q) \lambda= \sum_{i,j} C_{ij}(q) \lambda (x_j) x_i.$
We shall seek $\theta_*$ in the form
$\theta_{*}(q, \lambda)=(q, f(q)\lambda, \lambda)$, where
$f: U\longrightarrow L(\fh^*, \fh)$.  By definition,  $\theta_*$ is a flat
connection if and only if $\theta_{*} [\,\lambda, \lambda'\,]=
[\, \theta_{*}(\lambda), \theta_{*}(\lambda')\,]_*$ for
$\lambda, \lambda': U \longrightarrow \fh^*$  By using Eqn. (4.2),
a straightforward calculation shows that this is equivalent to the
following two conditions:
$$\eqalign {& (1)\, \, \, df(q) (\lambda'(q))\, \lambda(q)-
df(q) ( \lambda(q))\, \lambda'(q)
= - (dR(q) (\, . \,)\, \lambda(q), \lambda'(q))\cr
& (2)\, \, \, [R(q) (\lambda(q)), \lambda'(q)]+ [\lambda(q), R(q)
(\lambda'(q))]
=0}$$
for $q \in U$

\smallskip
On the other hand, the condition $[\, \theta_{*}(\lambda),\,
\widetilde{\psi}
(\xi)\,]_*
=\widetilde{\psi}(d \xi\cdot\lambda)$ is equivalent to

$$\eqalign {& (3)\, \, \,  (dR(q)  (\, .\,)\, \lambda(q), n)= 0\cr
& (4)\, \, \,  d \psi^{-1}(q)  \lambda(q) (n) -
[f (\lambda)(q), \psi^{-1}(q)(n)]\cr
& -  \big( [R(q)  (\lambda(q)), \psi^{-1}(q)(n)]+ [\lambda(q), R(q)
(\psi^{-1}(q) (n))]\big)=0,\cr}$$
for $n \in {\frak n}$ and $q \in U$.
\smallskip
From the properties of $R$, (2) and (3) are immediately seen to hold.
We now examine condition (4).  Set
$\psi(q) (e_\alpha)= \psi_\alpha (q) e_\alpha$, for all
$\alpha\in \Delta,$ then
$d \psi^{-1} (q) \lambda(q) (e_\alpha)=\phi_\alpha (q) \psi_\alpha (q)^{-1}
(\lambda (q), \alpha) e_\alpha.$
 Meanwhile, it is easy to check that
$$\eqalign {&[f(\lambda)(q), \psi^{-1}(q) e_\alpha]= \psi_\alpha (q)^{-1}
\alpha (f(\lambda (q)))
e_\alpha\cr
&\big( [R(q) (\lambda(q)), \psi^{-1}(q) e_\alpha]+ [\lambda (q), R(q)
(\psi^{-1}(q) e_\alpha)]\big)\cr
& =\big( \psi_\alpha(q)^{-1}
( \alpha (C^\#(q) (\lambda (q)))+
(\alpha, \lambda (q)) \phi_\alpha (q)\big)e_\alpha.\cr}$$
\smallskip
Therefore, condition (4) is equivalent to
$$\alpha \big(f(\lambda) (q) + C^\# (q) (\lambda (q))\big) =0, {\hbox {
for all }}
\alpha\in \Delta,$$
that is to $f= - C^\#.$
Finally, inserting $f=-C^\#$ into condition (1) shows that it is
trivially satisfied as it is equivalent to the closedness of the
$2-$ form $\sum_{i,j} C_{ij} dq^i \otimes dq^j.$ Hence
we have
\medskip
{\bf Proposition 4.3}  {\it The map
$$
\align
&\theta_{*}: TU \simeq U \times \fh^*
\longrightarrow U\times \fh \times \fg^*\\
&\quad (q, \lambda)\mapsto
\big(q, -C^\#(q)\lambda, \lambda)\\
\endalign
$$
is a flat connection on $U\times \fh\times \fg^*$ satisfying
$[\, \theta_{*}(\lambda),\, \widetilde{\psi}(\xi)\,]_*
=\widetilde{\psi}(d \xi\cdot\lambda)$
for $\lambda: U\longrightarrow \fh^*$ and $\xi: U\longrightarrow\fg'.$}
\medskip

{\bf Theorem 4.4} (Trivialization)
{\sl Let $A':={U} \times \fh^*\times \fg'$ be
the trivial Lie algebroid over U (see Prop. 3.2.1), then
the (bijective)  bundle map
$$\eqalign {&\sigma:A' \longrightarrow U\times \fh\times \fg^*\cr
&(q, \lambda, \xi)\mapsto \theta_{*}(q, \lambda)+\widetilde{\psi}(q, \xi)
\cr}$$
is an isomorphism of Lie algebroids.  Its inverse is given by
$\tau(q, Z, \lambda + n)= \mathbreak (q, \lambda, -C^\#(q)\lambda - Z + \psi
(q)n).$}
\medskip
{ Proof.} This is clear from the properties of $\theta_*$ and the
fact that $\widetilde{\psi}$ is an isomorphism of Lie algebra bundles.$\pf$
\bigskip
Note that the theorem implies, in particular, that the dynamical Lie
algebroid $A(X)^*\simeq U\times \fh\times \fg^*$ is integrable.  In what
follows, we let $U'$ be a connected and simply-connected open subset of
$U$ and we consider the coboundary Poisson groupoid
$X(U')=U' \times G\times U'$ associated with $R$.  We also
let $G'$ be the connected and simply
connected Lie group with $Lie (G')= \fg'$ and denote by
$$ T: X(U')^* \longrightarrow X'={U}'\times G'\times {U}'$$
 the unique (base preserving) Lie groupoid isomorphism such that $A(T)=
 \tau_{\mid_ {U'\times \fh\times \fg^*}}.$
\smallskip 
 If we define the Poisson bracket on $X'$ by
$$\{f,g\}_{X'}= \{f\circ T, g\circ T\}_{X(U')^*} \circ T^{-1},$$
then $(X', \{\, , \, \}_{X'})$ and
$(X(U') , \{\, , \, \}_{X(U')})$ are Poisson groupoids in
duality (see Def 3.1.1)
and $T$ is an isomorphism of Poisson groupoids.
\medskip
The following theorem characterizes the Poisson groupoid $(X',\{\, ,\,
\}_{X'}).$ 
 
\medskip
Let $j: \fh\longrightarrow \fg': Z\mapsto Z$ be the inclusion.
\medskip
{\bf Theorem 4.5}
{\it  The Poisson groupoid $(X', \{\, , \, \}_{X'})$ is of dynamical type
 with Poisson bracket
$$\eqalign {&\{f,g\}_{X'}(p,u,q)= <p,[\delta_1 f, \delta_1 g]> -
<q,[\delta_2f,\delta_2g]>\cr
&\hskip 30 pt -<j \delta_1f, Dg> -<j \delta_2 f, D'g>\cr
&\hskip 30 pt +<j\delta_1g,
Df>+<j\delta_2 g, D'f>\cr
&\hskip 30 pt+ <Df, P'(p,u,q) Dg>,\cr}$$
where  $P': {U}'\times G'\times {U}'\longrightarrow L({\fg'}^*,
\fg')$
is the unique skew symmetric groupoid cocycle whose tangent cocycle
$P'_* (q, {\bold \Lambda}, {\bold Z}+{\bold n}):=-\delta_1 P'\,(
{\bold
\Lambda })+ \partial P' \, ({\bold Z} + {\bold n})$ is given by
$$\eqalign {& (n, \delta_1 P'\, ({\bold \Lambda})\, n')= (\Pi_{\fh} [\psi^*
n,
\psi^* n'], {\bold \Lambda})\cr
& (\lambda, \delta_1 P'\,( {\bold \Lambda}) \,n')= 0\cr
& (\lambda, \delta_1 P' \,({\bold \Lambda}) \,\lambda')= (dC^\#(\lambda')
\,\lambda - dC^\# (\lambda) \,\lambda', {\bold \Lambda})\cr
&  (n, \partial P'\,( {\bold n})\, n')=- (\Pi_{\frak n} [\psi^* n,\psi^*
n'],
\psi^{-1} {\bold n})\cr
&  (\lambda, \partial P' \,({\bold n} )\,n')= - ( d({\psi^*}^{-1}
)(\lambda)
\,\psi^* n', {\bold n})- ([C^\# \lambda, n'], {\bold n})\cr
&  (\lambda, \partial P' \,({\bold n})\, \lambda') = 0\cr
&  \,\,\partial P'\, ({\bold Z})= 0.\cr}$$
Here, $\Pi_{\fh}$ and $\Pi_{\frak n}$ are the projections relative to
the direct sum decomposition $\fg=\fh\oplus {\frak n}$,
 $n,n',{\bold n}\in {\frak n}, {\bold Z}\in \fh$, and $\lambda,\lambda',
{\bold
\Lambda}\in 
\fh^*\simeq \fh,$ and the differentials of $P'$ are taken at $(q,1,q).$}
\medskip
Proof. For the sake of clarity, we shall begin by
repeating the argument of Theorem 3.2.4 here.
Consider the Poisson groupoid morphism
$$J: X\longrightarrow U'\times U': (p,h,q)\mapsto (p,q),$$
with induced Lie bialgebroid morphism
$A(J).$   Applying Prop. 3.2.3
 to the dual morphism
 $$A(J)^*:U'\times \fh\longrightarrow U'\times \fh\times \fg^* $$
we infer the existence of a (base preserving)
Poisson groupoid morphism
$$J^*: H\times U'\longrightarrow X^*$$
and so of a morphism
$$I'=T\circ J^*:H\times {U}'\longrightarrow {U}'\times G'\times {U}'.$$
Now, (see Eqn (2.2.1)) $I'$ is necessarily of the form
$$(h,p)\mapsto (p, \chi' (h,p), p),$$
for some groupoid morphism $\chi': H\times {U}'\longrightarrow G'$
with tangent map $A (\chi '): U'\times \fh\longrightarrow \fg'.$
Therefore, the  Poisson bracket $\{\, , \, \}_{X'}$ is given by Thm
2.2.5 for  
$A_{\chi'}$ and some groupoid $1-$ cocycle $P': X'\rightarrow
L({{\fg}'}^*, {\fg}').$
\smallskip
Denote  the Lie algebroid of ${U}'\times G'\times {U}'$ by $A'$
and 
let $(A', [\, ,\, ]', a'; {A'}^*, [\, ,\, ]'_*, a'_*)$ be the
Lie bialgebroid structure of Prop. 3.2.1.  For ${A'}^*$, we have
$K=0$ since $\fh$ is abelian, and
$A=A_{\chi '}.$
The Poisson bracket $\{\, , \, \}_{X'}$ is now uniquely determined by the
duality
requirement (see Def. 3.1.1) that
the trivialization map $\tau$ of Thm. 4.4 be a Lie bialgebroid isomorphism
from
$({A(X)^*}, A(X)^-)$ to $(A',{A'}^*)$, that is, by the condition that
the map
$$
\align
&\tau^*: {A'}^*=U'\times \fh\times {\fg '}^*\longrightarrow
A(X)^-=U'\times \times \fh^*\times \fg\\
&\quad (q, Z, \lambda+n)\mapsto
\big(q,-\lambda, C^\# (\lambda) + Z + \psi^* (n)\big)\\
\endalign
$$
satisfies
$$ - a_{A(X)} \, \tau^* = a'_*,\quad \tau^* [(Z, \lambda+n),
(Z', 
\lambda'+n')]'_*= -[\tau^* (Z, \lambda+n), \tau^* (Z',\lambda
'+n')]_{A(X)}.$$
Now (see Prop. 3.2.1) the anchor condition is equivalent to $A_{\chi'} (q)
(Z) = j (Z)= Z$ so the
groupoid morphism $\chi': H\times {U}' \longrightarrow G'$ is just the
inclusion of $H$ into $G'.$ Therefore (see Def. 2.2.6) $X'$ is of
dynamical type.
On the other hand, a direct calculation shows that the bracket
condition holds if and only if
$\delta_1 P'$ and $ \partial P'$ satisfy the equations given above. Hence
the claim. $\pf$
\bigskip
{\bf Remarks 4.6}
\smallskip
(a) The relationship between $P'$ and $P'_*$ is as
follows. Fix $q_0\in U'$ and write $P'$ as
$$P'(p,u,q)= - l (p) + \pi (u) + Ad_u l(q) Ad_u^*$$
for some map $l: {U}'\rightarrow L({\fg'}^*, \fg')$ with $l(q_0)= 0$ and
some group cocycle $\pi: G'\rightarrow L({\fg'}^*, \fg').$ We have
$$P'_* (q, {\bold \Lambda} , X')= d\pi (1) \, X' + ad_{X'} l(q) + l(q)
ad_{X'}^* + dl (q) {\bold \Lambda}.$$
Therefore (this is a special case of a result of \cite{X})
$$dl(q)( {\bold \Lambda})= P'_* (q, {\bold \Lambda}, 0),\quad d\pi (u)
T_1l_u X'=Ad_u d\pi (1) X' Ad_u^*= Ad_u P'_* (q_0, 0, X')\,
Ad_u^*.$$
\smallskip
(b) Writing out the equations of Thm 4.5 for $\delta_1 P'$ using the basis
$(e_\beta)$ 
of $\frak n$ and integrating yields
$$\eqalign {& l(q) (e_\beta)= (\phi_\beta (q) - \phi_\beta (q_0)) e_\beta,\,
{\hbox {
if }} \beta\in <\Gamma>\cr
&l(q) (e_\beta) = e^{(\beta, \mu)} (e^{-(\beta,q_0)}- e^{-(\beta, q)})
e_\beta,\, {\hbox {
if }} \beta\in {\overline \Gamma}^+\cr
&l(q) (e_\beta) = e^{-(\beta, \mu)} (e^{(\beta, q)}- e^{(\beta, q_0)})
e_\beta,\, {\hbox {
if }} \beta \in {\overline \Gamma}^- \cr
&l(q) (\lambda) = (C^\# (q) -C^\# (q_0)) (\lambda)\, {\hbox { for all }}
\lambda\in \fh^*\cr}$$
On the other hand, the remaining equations evaluated at $q=q_0$ give
$$\eqalign {&<\lambda, d\pi (1) ({\bold n}) e_\beta>= - (\phi_\beta (q_0)
(\lambda, \beta) + (C^\# (q_0) \lambda, \beta)) <e_\beta , {\bold n}>\cr
&<e_\alpha, d\pi (1) ({\bold n}) e_\beta>=- {\psi_{-\alpha} (q_0)
\psi_{-\beta}
(q_0)\over 
\psi_{-(\alpha+\beta)} (q_0) } \, <[e_\alpha, e_\beta], {\bold n}>\cr
&<\lambda, d\pi (1) ({\bold n}) \lambda'>= 0,\quad d\pi (1) {\bold Z}=
0\cr}$$
for all $\alpha, \beta\in \Delta, \lambda, \lambda'\in \fh^*, {\bold Z}+
{\bold n} \in \fg' $, where we have set $\psi (e_\alpha) = \psi_\alpha (q)
e_\alpha.$
\smallskip
The latter equations allow, in principle, for an explicit expression of $\pi
(u)$ but 
we shall postpone this integration to a future publication as it will not be
needed in the rest of this paper.
\bigskip
\bigskip 
{\bf 5. Coboundary dynamical Poisson groupoids - the constant $r-$ matrix
case.}
\bigskip
The purpose of this section is two-fold.  In Section 5.1, we give a
construction of the dual $X^*$ of the
coboundary dynamical Poisson groupoid $X= \fh^*\times G\times \fh^*$
(of Theorem
2.1.4) 
for the constant $r-$ matrix case, i.e., for the case where $R$ is a
constant
map from $\fh^*$ to $L(\fg^*, \fg)$.
As the reader will see,the construction involves the use of Poisson Lie
group theory. More specifically, the Poisson Lie group $G$ equipped with the
Sklyanin bracket admits an extension to a bigger Poisson Lie group whose
dual 
is critical in the construction.
In Section 5.2, we construct a
symplectic 
double groupoid which has $X$ and $X^*$ as its side groupoids.  This leads,
in
particular, to a description of the symplectic leaves of $X$ as orbits
of a Poisson Lie group action.
We shall discuss the non-constant $r-$ matrix case  in a forthcoming
publication.
\bigskip
{\bf 5.1. The dual Poisson groupoid. }
\bigskip
Let $\iota: \fh\longrightarrow \fg$ be the inclusion map. We assume here
that
the Lie groups $G$ and $H$  are connected and
simply
connected.
Let $R: \fg^*\longrightarrow \fg$ be a
skew-symmetric constant $r-$ matrix which satisfies
Eqn. (2.1.1) and Eqn.(2.1.2). Recall that the group $G$ equipped with the 
Sklyanin
bracket
$$\{f,g\}_G (x)= <R (Df), Dg>- <R( D'f), D'g> \eqno {(5.1.1)}$$
is a Poisson Lie group with tangent Lie bialgebra
$(\fg, [\, , \, ]; \fg^*, [\, , \, ]_\ast)$ where
$$[A, B]_\ast= 
ad^*_{R(A)} B - ad^*_{R (B)} A.\eqno (5.1.2)$$
\medskip
{\bf Lemma 5.1.1} {\it $H$ is a trivial Poisson Lie subgroup of $G$.}
\medskip
Proof. 
Since $R$ is $\fh -$ equivariant, we have

$$\eqalign {<[A,B]_\ast, \iota Z>=& <ad^*_{R(A)} B, \iota Z> -
<ad^*_{R (B)}
A,\iota Z>\cr
=&<B, [R (A),\iota Z]>+ <A, [\iota Z, R (B)]>\cr
=&<B, [R (A),\iota  Z]> -<A, R (ad^*_{\iota Z} B)>= 0,\cr}$$
therefore $[\fg^*, \fg^*]_\ast \subset \fh^\perp$. In particular, $
\fh^\perp$ is an ideal in $\fg^*$ and the connected Lie subgroup
 $H\subset G$ is a Poisson Lie subgroup with tangent
Lie bialgebra $(\fh,\fh^*)$ defined by
$$[\iota^* A, \iota^* B]_{\fh^*}:=\iota^* [A,B]_\ast.$$
 Hence the Lie bracket of $\fh ^*$ is identically zero and $H\subset
G$ is a trivial Poisson Lie subgroup. $\pf$
\bigskip
Let $(G^*, \{\, , \, \}_\ast)$ be Drinfeld's Poisson Lie group dual
to $(G,\, \{\, , \, \}_G), $
and let (\cite{STS}, \cite{LW})
$$\varphi^+ : G^*\times
{\overline G}\longrightarrow G^*,\quad \varphi^-: {\overline {G^*}}\times G
\longrightarrow
G \eqno (5.1.3)$$
be the right and left dressing actions. Recall that $\varphi^+$ and
$\varphi^-$ are 
Poisson Lie group actions and that $G$ and $G^*$ act on each
other by twisted
automorphisms 
$$\varphi^+ _x (uv)= \varphi^+_{\varphi^- _v (x)} (u) \, \varphi^+ _x (v),
\quad \varphi^- _{u} (xy)= \varphi^- _u (x) \, \varphi^- _{\varphi^+ _x
(u) } (y).\eqno (5.1.4)$$
They are related to (in fact defined by) the Poisson brackets of $G$ and
$G^*$ by the 
formulae  
$$\eqalign {\{\phi, \psi\}_\ast (u) &= - d\phi (u)
\lambda^+ (T_1 l_u^* d\psi (u)) (u)\cr
\{f, g\}_G (x)& = df (x)  \lambda^- (T_1r_x^* dg (x)) (x),\cr} \eqno
(5.1.5)$$ 
where $\lambda^+ (X) (u), \lambda^- (A) (x)$ are the infinitesimal
generators of $\varphi^+ $ and $\varphi^-.$
\smallskip
 It follows from Lemma 5.1.1 that the
restriction of the  right dressing action to $H$ induces a left Hamiltonian
action 
$$\phi^l: H\times G^* \longrightarrow G^*: (h,u)\mapsto \varphi^+
_{h^{-1}} (u).\eqno (5.1.6)$$
Moreover, since $G^*$ acts trivially on $H\subset G$ (i.e. $\varphi^-_u
(h)= h, \, u\in G^*,\, h\in H$), by (5.1.4) for each $h\in H,\,\phi^l _h $
is an automorphism of $G.$
\medskip
{\bf Caveat}
Note that our conventions differ from those of \cite{LW1}. Indeed
$\varphi^+$
is their left dressing action made right, while $\varphi^-$ is their right
dressing action made left.
\smallskip
  Recall also that the dressing vector fields $\lambda^+$ and
$\lambda^-$ may fail to globally integrate to define $\varphi^+$ and
$\varphi^-.$  In this subsection, we need only assume that $\phi^l$ is
globally defined.
\medskip
By construction, the map $\iota^* : \fg^*\longrightarrow
\fh^*$ is a 
morphism of Lie algebras.
Let 
$$I^*: G^* \longrightarrow \fh^*$$ be the (unique) morphism of Lie
groups 
integrating $\iota^*.$
\medskip
{\bf Lemma 5.1.2}  {\it $I^*$ is an $Ad^*_H -$
equivariant momentum map for the action $\phi^l$.}
\medskip
Proof. For $Z\in \fh$, let $j_Z\in C^\infty (G^*, \fh^*)$ be defined by
$j_Z(u)= <I^*
(u), Z>.$ We have to show that the Hamiltonian vector field ${\widehat
{X}} _{j_Z}$ coincides with the infinitesimal generator $-\lambda^+ (\iota
Z)$ of the action $\phi^l .$
But 
$$\eqalign { {d\over dt}_{\mid_0} j_Z (u e^{tA})&= {d\over dt}_{\mid_0}
<I^* 
(ue^{tA}), Z>\cr
&= {d\over dt}_{\mid_0} <I^* (u) + I^* (e^{tA}), Z>= <\iota^* A,
Z>.\cr}$$
Therefore 
$$\eqalign { {\widehat {X}} _{j_Z} (u)&= -\lambda^+ (T_1l_u^* dj_Z
(u))(u)\cr
&= -\lambda^+ (T_1l_u^* T^*_u l_{u^{-1}} \iota Z) (u)= -\lambda^+
(\iota 
Z) (u).\cr}$$
It remains to show 
$Ad^*_{h^{-1}}\, I^* (u) = I^* (\varphi^+
_{h^{-1}} (u)).$
Since both sides are group morphisms from $G^*$ to $\fh^*$
and $G^*$ is connected, it is enough to check
that the induced Lie morphisms are equal. Now
$${d\over dt}_{\mid_0}  I^* (\varphi^+ _{h^{-1}} (e^{tA}))={d\over
dt}_{\mid_0} Ad^*_{h^{-1}} I^* (e^{tA})$$
reads as
$$\iota^* T_1 \varphi^+ _{h^{-1}} A= Ad^*_{h^{-1}} \iota^* A.$$
But this equality follows from $T_1\varphi^+_{h^{-1}}= Ad^*_{h^{-1}}$
and the $Ad^*_H -$ equivariance of $\iota^*.$ Hence the claim.
\smallskip
(Note that one may also use  functoriality applied to the bialgebra
morphism $\iota.$) $\pf$
 \bigskip
{\bf Proposition 5.1.3} {\it
\smallskip
(a) The set $G\times \fh^*$ equipped with the multiplication
$(x,p)(y,q)= (xy, p+q)$ and  the Poisson bracket
$$\eqalign {&\{f,g\} (x,p)= <p, [\delta f , \delta g] + <R (Df), Dg>-
<R (D'f), 
D'g>\cr
& + <(Df-D'f), \iota (\delta g)> - <(Dg-D'g), \iota (\delta f)>\cr}$$
is a Poisson Lie group.
\smallskip
(b) The Drinfeld Poisson Lie group dual of $(G\times \fh ^*, \{\, , \,
\})$ 
is the set
$H\times G^*$  equipped with the semi-direct multiplication
$$(h,u)(k,v)= (hk, u\varphi^+
_{h^{-1}} (v))$$ and the Poisson bracket
$$\{\phi, \psi\}_\ast (h,u)= -\partial_* \phi \lambda^+ (T_1^*
l_u \partial_* \psi) (u),$$
where $\partial_* \phi $ is the partial derivative w.r.t. $G^*.$}
\medskip
Proof. (a)  This may be checked by a standard calculation which makes use
of Lemma 5.1.1, so we shall leave out the
details of the verification.
\smallskip
(b)  The tangent Lie bialgebra of $G\times \fh^*$ is given by
$(\fg\oplus \fh^*, [\, ,\, ]_\oplus; \fh\ltimes\fg^*, [\, , \, ]')$ where
$$[Z+A, Z'+A']'= [Z,Z'] - ad^*_{\iota (Z)} A' + ad^*_{\iota Z'}
A + [A,A']_\ast.$$
Therefore the dual group is $H\ltimes G^*$ as stated, while the
multiplicativity
of the Poisson bracket $\{\, , \, \}_*$ follows from the Hamiltonian
property (5.1.6) and the multiplicativity of the Poisson bracket of $G^*.$
Hence the assertion.
\smallskip
Note that $\fh\ltimes \fh^\perp \subset \fh\ltimes \fg^*$ is a Lie
subalgebra which is isomorphic to the vertex Lie algebra of Thm.
3.2.4.$\quad\pf$
\bigskip
Let $X= \fh^*\times G\times \fh^*$ be the dynamical Poisson groupoid
of Thm 2.1.4 with constant $r - $ matrix taken to be $-R$.  By the proof of
Thm. 3.2.4 and
Prop. 2.2.1, the dual groupoid
$X^*$ belongs 
to ${\Cal C}_{\fh^*}$. In the theorem below we shall give the
explicit ${\Cal C}_{\fh^*}$ structure of $X^*.$
\medskip
 If $f\in C^\infty 
(H\times \fh^*\times G^*)$ we define $\delta f$ and $D'f$ by the formulae
$$<\delta f, \lambda>= {d\over dt}_{\mid_0} f(h, p+t\lambda, u),\, <D'f, Z>=
{d\over dt}_{\mid_0} f(he^{tZ}, p, u),$$
 and we denote by $\partial_* f$  the partial derivative w.r.t. $G^*.$
\medskip
{\bf Theorem 5.1.4} (Dual Poisson groupoid (second form))
\smallskip
{\it Let $X$ be as above.
\smallskip
 (a) The 
set $\Gamma:= H\times \fh^*\times G^*$ together with
  the product Poisson bracket
$$\eqalign {&\{f,g\}_{\Gamma} (h,p,u)= -<D'g, \delta_1 f> + <D'f, \delta_1
g>\cr
&\hskip 50 pt - <p,
[\delta_1 f, \delta_1 g]>  -\partial_* f \lambda^+ (T_1l_u^* \partial_*
g)(u),\cr}$$ 
the commuting Hamiltonian actions of $H$
$$\phi^- _k (h,p,u) = (kh, p, \varphi^+ _{k^{-1}} u),\quad \phi^+_k (h,p,u)=
(hk, Ad^*_k \, p , u)$$
with equivariant momentum maps
$$j_- (h,p,u) = Ad^*_{h^{-1}} p + I^* (u),\quad j_+ (h,p,u)= p$$
($I^*$ is as in Lemma 5.1.2),
and the groupoid structure
$$\eqalign {&\alpha= j_-,\, \beta = j_+, \, \epsilon (q)= (1,q,1)\cr
&(h, j_-(k,q,v), u)\cdot (k,q, v)= (hk, q, u \varphi^+_{h^{-1}} v)\cr
&i(h,p,u)= (h^{-1}, j_-(h,p,u), \varphi^+ _h (u^{-1}))\cr}$$
is a Poisson groupoid in ${\Cal C}_{\fh^*}.$
\smallskip
(b) The Poisson groupoid  $\Gamma$ of (a) is the Poisson
groupoid dual $X^*$ of $X.$}
\bigskip
Proof. (a) That the actions $\phi^\pm$ are Hamiltonian with equivariant
momentum 
maps $j_\pm$ follows from Example 2.1.3,
the Hamiltonian property of the action $\phi^l$ (Eqn. (5.1.6)), and 
Lemma 5.1.2. On
the other hand, an easy verification, using
$\varphi^+_{h^{-1}} (uv)=
\varphi^+_{h^{-1}} (u) \, \varphi^+ _{h^{-1}} (v)$ and the $Ad^*_H - $
equivariance of $I^*$, shows that the groupoid axioms (for these
axioms see e.g.  \cite{W1})
are satisfied.  
The lengthy check that the graph of the multiplication
$$Gr (m)\subset \Gamma\times \Gamma\times {\overline {\Gamma}}$$
is a coisotropic submanifold is postponed to the appendix.
\smallskip
(b)  We have to show that the Lie
bialgebroid tangent to $\Gamma$ is isomorphic to the Lie bialgebroid $
(A(X)^*, A(X)) 
\simeq (\fh^*\times \fh\times \fg^*,
\fh^*\times \fh^*\times \fg)$ of
Prop. 3.2.1 (with constant $r-$ matrix $-R$). We shall only sketch the main
steps.
\smallskip
(i) The isomorphism $A(\Gamma)\simeq
A(X)^*.$ We have to compute (see the end of section 2.1) the value on
$\epsilon (\fh^*)$ of the Lie bracket of two left invariant
sections
$$X^l , {X'} ^l: \Gamma: \longrightarrow Ker\, T\,\alpha \subset T \Gamma.$$
We have
$$A(\Gamma)_q:=Ker\, T_{(1,q,1)} \alpha =\{ (Z, ad^*_Z q -\iota^* A,
A)\, \mid \, 
Z\in \fh, A\in \fg^*\}\simeq \fh\times \fg^*,$$
where the identification $\simeq$ is by dropping the middle term.
Now the left invariant vector field $X^l$ whose restriction to $\epsilon
(\fh^*)$ is 
$$X: \fh^*\longrightarrow A(\Gamma):
p\mapsto (q, Z(q), A(q))$$
is given by 
$$X^l (h,q,u)=(T_1 l_h\, Z(q), ad^*_{Z(q)} q -\iota^* A(q), T_1l_u
Ad^*_{h^{-1}} \, A(q)).$$
A lengthy calculation then shows that
$$[X, X'] (q):= [X^l , {X'}^l]_{\Gamma}(1,q, 1)$$
is given, after the identification $\simeq$, by
$$\eqalign { &[X, X'](q)=(dZ' (ad^*_Z q -\iota^* A)- dZ
(ad^*_{Z'} q -\iota^* A') + [Z(q), Z'(q)],\cr
&\hskip 90 pt +dA' (ad^* _Z q -\iota^* A) - dA (ad^*_{Z'} q-
\iota^* A')\cr
&\hskip 110 pt -ad^*_{\iota Z} A' + ad^* _{\iota Z'} A + [A(q),
A'(q)]_\ast),\cr}$$
where $(\fg^*, [\, ,\, ]_\ast)$ is as in Eqn.(5.1.2) and
all maps and differentials are evaluated at $q\in \fh^*.$ Thus the
bracket indeed coincides, up to sign,  with the one given
in Prop.
3.2.1 for $A(q)= \iota,\, K(q) Z= ad^*_Z q,$ and $P$ as in
Remark 3.2.2 with $r-$ matrix $-R.$
Now, applying the Lie functor to the morphism
$$[\alpha, \beta]: H\times \fh^* \times G^* \rightarrow \fh^*\times
\fh^*: 
(h,q,u)\mapsto (\alpha (h,q,u), \beta (h,q,u))$$
shows that the anchor $a_* (q, Z, A)= (q, ad^*_Z q -\iota^* A).$
\smallskip
(ii) The isomorphism ${ A(H\times \fh^* \times
G^*)}^*\simeq A(X).$ The unit section  is $\epsilon: \fh^*
\longrightarrow \Gamma:q\mapsto (1,q,1).$ Therefore
$\gamma\in N^* (\epsilon (\fh^*))_q$ if and only if $\gamma= (\lambda,
0, X)$ for some $\lambda\in \fh^*$ and $X\in \fg.$
\smallskip
Let $\theta, \theta': \fh^*\longrightarrow N^* (\epsilon
(\fh^*))\simeq \fh^* \times \fh^*\times \fg$
be two sections expressed as $\theta (q)= (\lambda (q), 0, X(q))$ and
$\theta' (q)= (\lambda' (q), 0, X' (q))$.  We set
$${\overline {\theta}} (h,q,u)=
(T^*_h l_{h^{-1}} \lambda (q), 0, T^*_u l_{u^{-1}} X(q)),$$ and similarly
for ${\overline {\theta'}}.$ By Eqn. (3.1.2), it
suffices to calculate
$$[\theta, \theta'](q)= [{\overline {\theta}},{\overline { \theta'}}]
(1,q,1),$$
where the {\sl rhs} is given by Eqn. (3.1.1) with Hamiltonian operator
$$\Pi^\#_\Gamma (h,q,u) ({\overline \theta})= (0, -\lambda (q), -\lambda^+
(X(q)) (u)).$$
A  calculation making use of standard properties of the Lie derivative and
of the dressing 
field $\lambda^+$ shows that  $[(\lambda, X),(\lambda' , X')] (q)$
 coincides with that of the trivial Lie algebroid of Prop. 3.2.1.
Finally observe that the anchor map, which is the restriction of $-\Pi^\#$
to $N^* (\epsilon
(\fh^*))$, is given by $a(q, \lambda, X)= (q, \lambda).$ (Note that the
dressing field $\lambda^+ $ vanishes at $u= 1.$) This concludes the proof
of the theorem.   $\pf$
\bigskip
In the special case when $R=0$, we have $G^*= \fg^*$ equipped with the
Lie-Poisson structure. In this case $I^*= \iota^* : g^*\longrightarrow h^*$
and $\varphi^+_{h^{-1}} (A) = Ad^*_{h^{-1}} (A).$   Specializing Thm 5.1.4
to this situation, we
have
\bigskip
{\bf Corollary 5.1.5} (Dual Poisson groupoid for $R=0.$)
\smallskip
{\it Let $X=\fh^*\times G\times \fh^*$ be the coboundary dynamical Poisson
groupoid 
of Thm 2.1.4 with
$R=0.$ Then the Poisson groupoid dual $X^*$ of $X$ is the set $H\times
\fh^*\times \fg^*$ equipped with the Poisson bracket
 $$\eqalign {&\{f,g\}_{*} (h,p,A)= -<D'g, \delta_1 f> + <D'f, \delta_1
g>\cr
&\hskip 50 pt - <p,
[\delta_1 f, \delta_1 g]> + <A, [\delta f, \delta g]>,\cr}$$
 ($<\delta f, B>= {d\over dt}_{\mid_0} f(h,p,A+tB)$),
and the groupoid structure
$$\eqalign { & \alpha (h,p,A)= Ad^*_{h^{-1}} p + \iota^* A,\quad \beta
(h,p,A)= 
p,\quad \epsilon (p)= (1,p,0)\cr
&(h,\alpha (k,q,B), A)\cdot (k,q,B)= (hk, q, A+ Ad^*_{h^{-1}} B)\cr
&i (h,p,A)= (h^{-1}, \alpha (h,p,A), -Ad^*_h A).\quad \pf  \cr}$$}
\bigskip
We now describe the trivialization of the Lie groupoid $\Gamma=
H\times \fh^* \times G^*.$
\bigskip
Let $H^\perp$ be the connected and simply connected Lie subgroup of $G^*$
with $Lie (H^\perp)= (\fh^\perp, [\, ,\, ]_*)$ and
$$H\ltimes H^\perp \subset H\ltimes G^*$$
be the Lie subgroup with Lie algebra $\fh\ltimes \fh^\perp\subset \fh\ltimes
\fg^*$ (see the note in the proof of Prop. 5.1.3).
\bigskip
{\bf Proposition 5.1.6}  (Trivialization)
\smallskip
 {\it Equip $\fh^*\times (H\ltimes H^\perp)\times \fh^*$ with the
 trivial Lie groupoid structure over $\fh^*$, and let $\Gamma$ be the Lie
groupoid in
Theorem
 5.1.4. If $s$ is an arbitrary linear section $s:\fh^*\rightarrow
\fg^*$  of $\iota^*,$ the map
$$\eqalign {&\Sigma: \fh^*\times (H\ltimes H^\perp)\times
\fh^*\longrightarrow
\Gamma\cr
&(p, (k,u), q)\mapsto \big(k, q,  exp (s(p))\, u\, \varphi^+_{k^{-1}}
 (exp (-s(q)))\big)\cr}$$
is a Lie groupoid isomorphism.}
\bigskip
Proof. We use an elementary device (see \cite{M1}) according to which
if $\sigma: \fh^*\longrightarrow \Gamma$ is a global smooth section of the
restriction of $\alpha$ to the fiber $\beta^{-1} (0)$ and $G'\subset \Gamma$
is the 
isotropy subgroup at $0$, then the map
$$\eqalign {&\Sigma: \fh^*\times G'\times \fh^*\longrightarrow \Gamma\cr
&(p, x', q)\mapsto \sigma (p)\cdot x'\cdot i(\sigma (q))\cr}$$
is a Lie groupoid isomorphism.
\smallskip
Now $G'=\beta^{-1} (0)\bigcap \alpha^{-1} (0) = H\ltimes H^\perp$,
while $\alpha_{\mid_{\beta^{-1} (0)}}:\beta^{-1} (0)\rightarrow \fh^*:
(h,0,u)\mapsto I^* (u).$
Observe that for any linear section $s:\fh^* \rightarrow \fg^*$, the map
$$\sigma: \fh^*\rightarrow \Gamma: p\mapsto (1, 0, exp \, s(p))$$ is a
smooth section of  $\alpha_{\mid_{\beta^{-1} (0)}}$ since
$I^* (\sigma (p))= I^* (exp \,s(p))= exp_{h^*}\, \iota^*
(s(p))= p.$ Calculating $\sigma (p)\cdot (h,0, u)\cdot i(\sigma (p))$ in
$\Gamma$ immediately yields the claim. $\pf$
\bigskip
{\bf Caveat} Note that if $\fh\not= 0$ the group $G'$ is not in general
isomorphic to the
the Poisson Lie group dual $G^*.$ For example, if $R=0$, $G^*=
\fg^*$ but $G'= H\ltimes \fh^\perp$. So $G'$ and $G^*$ may differ even
topologically.
\bigskip
{\bf Remark 5.1.7} By Thm. 3.2.4, the Poisson bracket $\{\, , \, \}_*$
on $\fh^* \times 
(H\times H^\perp)\times \fh^*$ defined by $\{\Sigma^*, \Sigma^*
g\}_*:=\Sigma^* \{f,g\}_\Gamma$ has the form given by Thm. 2.2.5.
However, even for standard Poisson Lie groups, the explicit bracket
transport turns out to be rather
cumbersome.
\bigskip
We close this subsection with the following
\medskip
{\bf Example 5.1.8} Let $\fg= \fh\oplus {\frak n}_+ \oplus {\frak n}_- $ be
the root space decomposition
 of a complex simple Lie algebra $\fg$, as in section 4.
Let $R= \Pi_{{\frak n}_+} - \Pi_{{\frak n}_-}$ be the standard $r -$ matrix.
In what follows, we shall
scale the Poisson bracket by $1/2$ to match with standard conventions.
Let ${ N}_\pm$ be the (connected and simply connected) unipotent subgroups
of $G$ with Lie algebra ${\frak n}_\pm.$ Note that $H=\fh\simeq \fh^*.$
\smallskip
The dual group $G^*$ is the set $\fh\times { N}_+\times { N}_-$
with semi-direct group law
$$(Z, A, B)\cdot (Z', A', B')= (Z+Z', A e^{Z\over 2} A' e^{-{Z\over 2}},
e^{- {Z\over 2}}B' e^{Z\over 2} B).$$
The dressing action of $\fh$ on $G^*$ is given by
$$\varphi^+ _{-Y} (Z, A, B)= (Z, e^Y A e^{-Y}, e^Y B e^{-Y}).$$
The Poisson Lie group $H\ltimes G^*$ of Prop. 5.1.3 (b) is the set
$\fh\times 
\fh\times { N}_+ \times { N}_-$ with group law
$$
\align
&(Y, Z, A, B)\cdot (Y', Z', A', B') = (Y+Y', Z+Z', A e^{Z\over 2} e^Y A'
e^{-Y} e^{-{Z\over 2}},\\
&\hskip 160 pt e^{- {Z\over 2}} e^Y B' e^{-Y} e^{Z\over 2} B),\\
\endalign
$$
and hence the vertex subgroup $G'$ of the groupoid $\Gamma$ is the set
$\fh\times 
{ N}_+\times { N}_-$ with group law
$$(Y, A, B)\cdot (Y', A', B')= (Y+Y', A e^Y A' e^{-Y}, e^Y B' e^{-Y} B).$$
Finally, the map 
$$\eqalign { &\fh^*\times G'\times \fh^*\longrightarrow \Gamma= \fh^*\times
\fh\times { N}_+ \times { N}_-\times \fh^*\cr
& (p, Y, A, B, q)\mapsto (q, Y, e^{p\over 2} A e^{-{p\over 2}},
e^{-{p\over 2}} B e^{p\over 2} , p-q)\cr}$$
gives an explicit trivialization of the Lie groupoid $\Gamma.$
\bigskip
\bigskip

{\bf 5.2. Construction of the associated symplectic double groupoid.}
\bigskip

Our goal in this subsection is to construct, for the constant $r-$ matrix
case (taken to be $-R$), a
symplectic double groupoid having $\fh^*\times
G\times \fh^*$ and $H\times \fh^*\times G^*$ as its side Poisson
groupoids.
\smallskip
We begin by recalling the notion of double Lie groupoids \cite{E},
\cite{M3}, and
symplectic double groupoids \cite{W1}, \cite{LW2}, \cite{M2}.
\medskip
{\bf Definition 5.2.1}

{\it (a) A double Lie groupoid consists of a quadruple $({\Cal  S}; {\Cal
H}, {\Cal  V}, B)$ where ${\Cal  H}$ and ${\Cal  V}$ are Lie groupoids over
$B$, and 
${\Cal  S}$ is equipped with two Lie groupoid structures, a horizontal
structure with base ${\Cal  V}$, and a vertical structure with base ${\Cal
H}$, such that the structure maps (source, target, multiplication, unit
section and inversion) of each groupoid structure on ${\Cal  S}$ are morphisms
with respect 
to the other. We call ${\Cal  H}$ and ${\Cal  V}$ the side groupoids of
${\Cal  S}$, and $B$ the double base. $({\Cal  S}; {\Cal
H}, {\Cal  V}, B)$ is displayed as in Fig. 5.2.1 below.}
\smallskip
{\it (b)   A double Lie groupoid $({\Cal  S}; {\Cal
H}, {\Cal  V}, B)$ is called symplectic  if ${\Cal  S}$ is equipped with a
symplectic structure such that both ${\Cal  S} {\buildrel {\tilde
\alpha_{\Cal  H}, 
\tilde\beta_{\Cal  H}}\over \rightrightarrows} {\Cal
V}$ and ${\Cal  S} {\buildrel {\tilde \alpha_{\Cal  V},
\tilde\beta_{\Cal  V}}\over \rightrightarrows} {\Cal
V}$ are 
symplectic groupoids.}
\bigskip
$$\eqalign {&\hskip 30 pt {\Cal  S}\quad {\buildrel {\tilde \alpha_{\Cal
H},
\tilde\beta_{\Cal  H}}\over \rightrightarrows}\quad {\Cal
V}\cr
&{\tilde\alpha_{\Cal  V}}, {\tilde\beta_{\Cal  V} }\downdownarrows\hskip 45
pt
\downdownarrows \alpha_{\Cal  V}, \beta_{\Cal  V}\cr
&\hskip 30 pt {\Cal  H}\quad {\buildrel \alpha_{\Cal  H},\beta_{\Cal
H}\over
\rightrightarrows}\quad B\cr
&\cr
&\hskip 40 pt {\hbox { Fig. 5.2.1}}\cr }$$

\bigskip
We shall consider the case where the Poisson Lie group $G$ is complete. In
this
case, the Drinfeld double $D$ can be identified with $G\times G^*$
\cite{STS}, \cite{LW1} with multiplication
$$(g_1, u_1) \cdot (g_2, u_2)= ((\varphi^- _{u_2^{-1}} (g_1^{-1}))^{-1} g_2,
u_1 (\varphi^+_{g_1^{-1}} (u_2^{-1}))^{-1} ). \eqno (5.2.2)$$
As a first step in the construction, we show that $X= \fh^*\times G\times
\fh^*$ and $X^*= H\times \fh^*\times G^*$ form a matched pair
of Lie groupoids in the sense of the following
\medskip
{\bf Definition 5.2.2} \cite{M3} {\it Two Lie groupoids $\Cal  V$ and $\Cal
H$
over the same base $B$ are said to form a matched pair of Lie groupoids
iff the manifold
$${\Cal  V}\ast {\Cal  H}= \{ (v,h)\in {\Cal  V}\times {\Cal  H}\, \mid
\beta_{\Cal  V} (v)= \alpha_{\Cal  H} (h)\}$$
admits a Lie groupoid structure over $B$ such that

(a) the maps $h\mapsto {\overline h}= (\epsilon_{\Cal  V} (\alpha_{\Cal  H}
(h)), h)$ and $v\mapsto {\overline v}= (v, \epsilon_{\Cal  H} (\beta_{\Cal
V} (v)))$ are morphism of Lie groupoids from ${\Cal  H}$ and ${\Cal  V}$
to ${\Cal  V} \ast {\Cal  H}$ respectively ,
\smallskip
(b) the map ${\Cal  V}\ast {\Cal  H}\longrightarrow {\Cal  V} \ast {\Cal
H}: (v,h)\mapsto {\overline v} {\overline h}$ is a diffeomorphism.
\smallskip
In this case, the groupoid ${\Cal  V}\ast {\Cal  H}$ is called  the
matched product of ${\Cal  V} $ and ${\Cal  H}.$ }
\bigskip
{\bf Proposition 5.2.3} {\it The Lie groupoids $X$ and $X^*$ form a
matched pair with matched product given by the trivial groupoid
$\fh^*\times M \times \fh^*$ where the vertex
group $M = H\times (G\times G^*)$ is the direct product of $H$ with
the Drinfeld double (see (5.2.2)).}
\medskip
Proof. Clearly, $X\ast X^*$ may be identified with the manifold
$\fh^*\times M\times \fh^*$. Equip the latter with the trivial
groupoid structure. The groupoids $X$ and $X^*$ are embedded as wide
subgroupoids of $\fh^*\times M\times \fh^*$ through the morphisms
$(p,g,q)\mapsto (p,
1, g, 1, q)$, $(h,p,u)\mapsto (Ad^*_{h^{-1}} p + I^* (u) , h, h,
u, p).$ Finally, for $(p,h, g,u, q)\in \fh^*\times M\times \fh^*$, we
have the unique factorization
$$
\align
 & (p, h, g,u, q)= (p, 1, \varphi^- _u (gh^{-1}), 1, Ad^*_{h^{-1}} q +
I^* (\varphi^+ _{gh^{-1}} (u)))\\
&\hskip 50 pt\cdot (Ad^*_{h^{-1}} q +
I^*(\varphi^+ _{gh^{-1}} (u)), h, h, \varphi^+ _{gh^{-1}} (u), q).\\
\endalign
$$
Hence it follows that $(X, X^*)$ is a matched pair. $\pf$
\bigskip
We shall denote by ${\overline X}$ and ${\overline X}^*$ the images
of $X$ and $X^*$ under the morphisms in the proof above.
\smallskip
 Given $(h,p,u)\in X^*, (p, g,q)\in X$, the corresponding elements in
 $\fh^*\times M\times \fh^*$ are composable, and we have the unique
 factorization
$$\eqalign {& (Ad^*_{h^{-1}} p + I^* (u), h, h, u, q)
(p,1,g,1,q)\cr
&\hskip 25 pt = \big(Ad^*_{h^{-1}} p + I^* (u), 1, \varphi^- _u
(hgh^{-1}), 1, 
Ad^*_{h^{-1}} q+ I^* (\varphi^+ _{hgh^{-1}} (u))\big)\cr
&\hskip 20 pt \big(Ad^*_ {h^{-1}} q + I^* (\varphi^+ _{hgh^{-1}}
(u)), 
h, h, \varphi^+ _{hgh^{-1}} (u), q\big).\cr}\eqno (5.2.3) $$
We therefore obtain two maps
$$
\align
 & \psi^- : (H\times \fh^*\times G^*) \ast_{\alpha_X}
(\fh^*\times G\times \fh^*)\longrightarrow
 \fh^*\times G\times \fh^*\hskip 70 pt (5.2.4.a) \\
&((h,p,u),(p,g,q))\mapsto \big(Ad^*_{h^{-1}} p + I^* (u), \varphi^-
_u (hgh^{-1}), 
Ad^*_{h^{-1}} q+ I^* (\varphi^+ _{hgh^{-1}} (u))\big),\\
\endalign
$$
and 
$$\eqalign {& \psi^+ : (H\times \fh^*\times G^*)
\ast_{\beta_{X^*}} (\fh^* \times G\times \fh^*)\longrightarrow
H\times \fh^*\times G^*\cr
&\hskip 30 pt ((h,p,u),(p,g,q))\mapsto \big(h,q, \varphi^+_{hgh^{-1}}
(u)\big).\cr}\eqno 
(5.2.4. b)$$
\bigskip
{\bf Proposition 5.2.4} {\it $\,\, \psi^-$ is a left groupoid action of
$X^*$ on $X$ with moment map $\alpha_X$ and $\psi^+$ is a right
groupoid action of $X$ on $X^*$ with moment map $\beta_{X^*}.$
Furthermore, the following conditions are satisfied
\medskip
(a) $\beta_X (\psi^-_{g_-} (g_+)= \alpha_{X^*} (\psi^+ _{g_+}
(g_-)),\quad$
\medskip
 for all $g_+ \in X, g_-\in X^*$ with $\beta_{X^*} (g_-)= \alpha_X
(g_+)$,
\bigskip
(b) $\psi^- _{g_-} (g_+g_+') = \psi^- _{g_-} (g_+) \psi^- _{\psi^+_{g_+}
(g_-)} (g'_+),\quad$
\medskip
for all $g_-\in X^*, g_+, g'_+\in X$ with $\beta_{X^*} (g_-) =
\alpha_X (g_+)$, $\beta_X (g_+) = \alpha_X (g'_+)$,
\bigskip
(c) $\psi^+ _{g_+} (g_-g'_-) = \psi^+ _{\psi^- _{g'_-} (g_+)} (g_-) \psi^+
_{g_+} (g'_-),\quad$
\medskip
 for all $g_+\in X, g_-, g'_-\in X^*$ with  $\beta_{X^*} (g_-) =
\alpha_{X^*} (g'_-)$, $\beta_{X^*} (g'_-) = \alpha_X (g_+).$}
\bigskip
Proof. The proof consists of direct checking and we shall omit the details.
See, however, Prop. 5.2.9 below. $\pf$
\bigskip
From standard consideration \cite{M3}, the upshot of the above
proposition
is
that one can construct  a vacant double Lie groupoid $({\Cal  S}_{vac};
X^*, X, \fh^*)$ (vacant means that the double source map $\beta_+:
{\Cal  S}_{vac}\rightarrow X^*\ast_\beta X$ is a diffeomorphism) having
$X$ and $X^*$ as its side groupoids. Indeed the horizontal structure
of the vacant double is given by the left action groupoid $X^*
\ltimes X$ corresponding to $\psi^-$, while the vertical
structure is given by the right action groupoid $X^*
\rtimes X$ associated with $\psi^+$. However, ${\Cal  S}_{vac}$ 
is not the correct underlying double Lie groupoid of the symplectic double
groupoid which we are looking for, as is clear from dimension considerations. 
Nevertheless, as we
shall show in what follows, the sought for double Lie groupoid can be
constructed by extending the objects $X$ and $X^*.$ It turns out that
these extended objects act on the unextended ones through groupoid actions
which restrict to $\psi^\pm$, and the corresponding left/right action
groupoids then give the desired horizontal/vertical structures.
Before we carry out the details of this construction, let us make an
important remark.  As we know from Thm 5.1.4, the Poisson structure on
$X^*$ is the product of the standard symplectic structure on $H\times
\fh^*$ 
and the multiplicative structure on $G^*.$
 Since $H\times \fh^*$ is symplectic, the coarse  groupoid $(H\times
 \fh^*) \times (H\times \fh^*)^-\rightrightarrows H\times \fh^*$ is a
 symplectic groupoid. On the other hand, there is a symplectic groupoid
 $G\times G^*\rightrightarrows G^*$ \cite{Lu} with structure maps given as
 follows:
$$\eqalign { & \alpha (g,u)= \varphi^+ _{g^{-1}} (u^{-1})^{-1},\quad
\beta (g, u) = u\cr
&\hskip 20 pt (g, \alpha (g',u'))\cdot (g', u')= (gg', u'). \cr}\eqno
(5.2.5)$$
Therefore the product groupoid
$$(H\times \fh^*)\times (H\times \fh^*)^- \times (G\times
G^*)\rightrightarrows H\times \fh^*\times G^* \eqno (5.2.6)$$
is a symplectic groupoid over $X^*.$ It turns out that this product
groupoid is isomorphic to the right action groupoid alluded to above.
\medskip
We now introduce the extensions of $X$ and $X^*.$ Since $H$ is a
groupoid over a point, we have the product groupoid
$$X^*_e= X^* \times H\rightrightarrows \fh^* \eqno (5. 2.7).$$
Let $J_+ = \alpha_X$ and define $\Psi^-: X^*_e\ast _{J_+}
X\longrightarrow X$ by
$$\Psi^- _{(h,p,u,k)} (p, g, q)= \big(Ad^*_{h^{-1}}p + I^* (u),
\varphi^- 
_u (hgk^{-1}), Ad^*_{k^{-1}} q+ I^* (\varphi^+_{hgk^{-1}}
(u))\big).\eqno (5.2.8)$$
\bigskip
{\bf Proposition 5.2.5} {\it $\Psi^- $ is a left groupoid action of
$X^*_e$ on $X$ with moment map $J_+$ such that $\Psi^- _{(h,p,u,h) }
(p,g,q)= \psi^- _{(h,p, u)} (p,g,q)$ for all $(h,p,u)\in X^*,\,
(p,g,q)\in X.$}
\bigskip
Proof. Clearly, $({\tilde x}, x)= \big( ({\tilde h} , Ad^*_{h^{-1}} p+
I^* (u), {\tilde u}, {\tilde k}), (h, p, u,k)\big)$ is a composable
pair in $X^*_e$ and we have
$${\tilde x} x= ({\tilde h} h, p, {\tilde u} \varphi^+_{{\tilde h}^{-1} }
(u), {\tilde k} k).$$
Therefore,
$$\eqalign { &\Psi^- _{{\tilde x}x} (p,g,q)= \big(Ad^*_{({\tilde h}
h)^{-1}} p + 
I^* ({\tilde u} \varphi^+_{{\tilde h} ^{-1}} (u)),\, \varphi^- _{{\tilde
g} \varphi^+_{{\tilde h}^{-1}} (u)} ({\tilde h} h g
({\tilde k} k)^{-1}),\cr
&\hskip 30 pt Ad^*_{({\tilde k}k)^{-1}} q + I^*
(\varphi^+_{{\tilde h} h g ({\tilde k} k)^{-1}} ({\tilde u} \varphi^+
_{{\tilde h} ^{-1}} (u)))\big).\cr}$$
On the other hand,
$$\eqalign {&\Psi^-_{\tilde x} \Psi^-_x (p,g,q)= \big(Ad^*_{({\tilde h}
h)^{-1}} p+ Ad^*_{{\tilde h}^{-1}} I^* (u)+ I^* ({\tilde
u}),\, \varphi^- _{{\tilde u}} ({\tilde h} \varphi^-_u (hgk^{-1}) {\tilde
k} ^{-1},\cr
&\hskip 25 pt Ad^*_{({\tilde k} k)^{-1}} q + Ad^*_{{\tilde k} ^{-1}}
I^* 
(\varphi^+ _{hgk^{-1}} (u) + I^* (\varphi^+ _{{\tilde h} \varphi^-_u
(hgk^{-1}) {\tilde k}^{-1}} ({\tilde u})\big).\cr}$$
Since $I^*$ is an $H-$ equivariant homomorphism, the equality of the
first components is clear. Now using the same property of $I^*$, we have
$$\eqalign {& Ad^*_{{\tilde k} ^{-1}} I^* (\varphi^+ _{hgk^{-1}} (u))
) + I^* (\varphi^+ _{{\tilde h} \varphi^-_u
(hgk^{-1}) {\tilde k}^{-1}} ({\tilde u}))\cr
&= I^* (\varphi^+ _{{\tilde h} \varphi^-_u
(hgk^{-1}) {\tilde k}^{-1}} ({\tilde u}) \, \varphi^+ _{{\tilde h} h g
({\tilde k} k)^{-1}} (\varphi^+ _{{\tilde h} ^{-1}} (u))).\cr}$$
But from Eqn. (5.1.4) and the fact that $H$ is a trivial Poisson Lie subgroup
of $G$, we find
$$\eqalign {\varphi^+ _{{\tilde h} h g ({\tilde k} k)^{-1}} ({\tilde u}
\varphi^+ 
_{{\tilde h} ^{-1}} (u))&= \varphi^+ _{\varphi^- _{\varphi^+_{{\tilde
h}^{-1}}(u)} ({\tilde h} h g ({\tilde k} k)^{-1})} ({\tilde u})\,\,
\varphi^+ 
_{{\tilde h} h g ({\tilde k}k)^{-1}} (\varphi^+ _{{\tilde h} ^{-1}} (u))\cr
&= \varphi^+ _{{\tilde h} \varphi^- _u (hg ({\tilde k} k)^{-1}} ({\tilde
u}) \, \, \varphi^+ _{{\tilde h} h ({\tilde k} k)^{-1}} (\varphi^+ _{{\tilde
h} ^{-1}} (u))\cr
&= \varphi^+ _{{\tilde h} \varphi^- _u (hgk^{-1}) {\tilde k} ^{-1} }
({\tilde u})\,\, \varphi^+ _{{\tilde h} h g ({\tilde k} k)^{-1}} (\varphi^+
_{{\tilde h}^{-1}} (u)).\cr}$$
Hence we have equality of the third components. Finally it follows from
the above calculation that
$$\eqalign { \varphi^- _{\tilde u} ({\tilde h} \varphi^- _u (hgk^{-1})
{\tilde k}^{-1})&= \varphi^- _{{\tilde u}} \varphi^- _{\varphi^+ _{{\tilde
h}^{-1}} (u)} ({\tilde h} h g ({\tilde k} k)^{-1} )\cr
&=  \varphi^-_ {{\tilde u} \varphi^+_{{\tilde h}^{-1}} (u) }({\tilde h} h g
({\tilde k} k)^{-1}).\cr}$$
This completes the proof that $\Psi^- _{{\tilde x} x} = \Psi^- _{\tilde x}
\Psi^- _x.$ The assertion on the relationship between $\Psi^-$ and $\psi^-$
is clear . $\pf$
\bigskip
The following corollary is a direct consequence of the definition of an
action 

groupoid.
\medskip
{\bf Corollary 5.2.6} {\it The left action groupoid $X^*_e
\ltimes X\rightrightarrows X$
corresponding to $\Psi^-$ has structure maps given by}
$$\eqalign {& {\tilde \alpha}_{\Cal  H}:\, X^*_e
\ltimes X\rightarrow X: \big(
(h,p,u,k),(p,g,q)\big)\mapsto \Psi^- _{(h,p,u,k)} (p,g,q)\cr
&{\tilde \beta}_{\Cal  H}:\,X^*_e
\ltimes X\rightarrow X: \big( (h,p,u,k),
(p,g,q)\big)\mapsto (p,g,q)\cr
&{\tilde m}_{\Cal  H}:\,( X^*_e
\ltimes X)\ast (X^*_e
\ltimes X)\rightarrow X^*_e
\ltimes X\cr
&\hskip 5 pt  \big( (h_1, p_1, u_1, k_1), \Psi^- _{(h_2, p_2, u_2, k_2)}
(p_2, g_2, 
q_2)\big)\cdot ((h_2, p_2, u_2, k_2), (p_2,g_2,q_2)\big)\cr
&\hskip 15 pt= \big( (h_1h_2, p_2, u_1 \varphi^+ _{h_1^{-1}} (u_2), k_1k_2),
(p_2,g_2,q_2)\big)\cr
&{\tilde \epsilon}_{\Cal  H}: X\rightarrow  X^*_e
\ltimes X : (p,g,q)\mapsto \big(
(1,p,1,1), (p,g,q)\big)\cr
&{\tilde i}_{\Cal  H}: X^*_e
\ltimes X\rightarrow X^*_e
\ltimes X:\big( (h,p,u,k),
(p,g,q)\big)\mapsto \cr
&\hskip 5 pt\big( (h^{-1} , Ad^*_{h^{-1}} p+ I^* (u), \varphi^+ _h
(u^{-1}), 
k^{-1}),\cr
&\hskip 15 pt  (Ad^*_{h^{-1}} p + I^* (u), \varphi^- _u (hgk^{-1}),
Ad^*_{k^{-1}} p + I^* (\varphi^+ _{hgk^{-1}} (u)))\big). \quad\quad
\pf\cr}$$
\bigskip
For the extension of $X$, we consider the coarse groupoid $H\times
H\rightrightarrows H$ and let
$$X_e= (H\times H) \times X \rightrightarrows H\times \fh^*\eqno (5.2.9)$$
be the product groupoid. Introduce the map
$$J_-: X^*\longrightarrow H\times \fh^*: (h,p,u)\mapsto (h,p) \eqno
(5.2.10)$$
and define
$$\Psi^+ : X^*\ast_{J_-} X_e \longrightarrow X^*: \Psi^+
_{(h,k,p,g,q)} (h,p,u)= (k,q, \varphi^+ _{hgk^{-1}} (u)). \eqno (5.2.11)$$
\bigskip
{\bf Proposition 5.2.7} {\it $\Psi^+$ is a right groupoid action of $X_e$
on $X^*$ with moment map $J_-$ and we have $\Psi^+ _{(h,h,p,g,q)}
(h,p,u)= \Psi^+ _{(p,g,q)} (h,p,u)$ for all $(h,p,u)\in X^*,
(p,g,q)\in X.$}
\medskip
Proof. This is clear. $\pf$
\bigskip
{\bf Corollary 5.2.8} {\it The right action groupoid $ X^*
\rtimes X_e\rightrightarrows X^*$ corresponding to
$\Psi^+$ has 
structure maps given by}
$$\eqalign {&{\tilde \alpha}_{\Cal  V}:X^*
\rtimes X_e\rightarrow X^*: \big( (h,p,u),
(h,k,p,g,q)\big)\mapsto (h,p,u)\cr
&{\tilde \beta}_{\Cal  V} :X^*\rtimes X\rightarrow X^*: \big( (h,p,u),
(h,k,p,g,q)\big)\mapsto \Psi^+ _{(h,k,p,g,q)} (h,p,u)\cr
&{\tilde m}_{\Cal  V} : (X^* \rtimes X_e) \ast (X^* \rtimes X_e)\rightarrow
X^* \rtimes X_e:\cr
&\hskip 5 pt  \big( (h_1, p_1, u_1), (h_1, k_1, p_1, g_1, q_1)\big) \cdot
(\Psi^+ 
_{(h_1, k_1, p_1, g_1, q_1)} (h_1, p_1, u_1), (k_1, k_2, q_1, g_2,
q_2)\big)\cr
&\hskip 10 pt = \big( (h_1, p_1, u_1) , (h_1, k_2, p_1, g_1 g_2,
q_2)\big)\cr
&{\tilde \epsilon}_{\Cal  V}: X^* \rightarrow X^* \rtimes X_e:
(h,p,u)\mapsto \big(
(h,p,u),(h,h,p,1,p)\big)\cr
&{\tilde i}_{\Cal  V} : X^* \rtimes X_e \rightarrow X^* \rtimes X_e: \big(
(h,p,u),
(h,k,p,g,q)\big)\mapsto\cr
&\hskip 15 pt \big( (k,q, \varphi^+ _{(hgk^{-1})} (u)), (h,h, q, g^{-1} ,
p)\big).\quad\quad \pf\cr}$$
\bigskip
Let $Pr_1: X^*\longrightarrow H$ be the projection onto the first
factor of $X^*= H\times \fh^*\times G^*.$
\medskip
{\bf Proposition 5.2.9} {\it The groupoid actions $\Psi^\pm$ satisfy the
following properties:
\medskip
(a) $\beta_X (\Psi^- _ {(g_-, k)} (g_+))= \alpha_{X^*} (\Psi^+ _ {(Pr_1
(g_-), k, g_+)} (g_-))$
\medskip
for all $g_-\in X^*, g_+\in X$ with $\beta_{X^*} (g_-)= \alpha_X
(g_+)$ and all $k\in H$.
\bigskip 
(b) $\Psi^- _ {(g_-, k')} (g_+ g'_+) = \Psi^- _{(g_-, k)} (g_+)\, \Psi^- _
{\Psi^+ _ {(Pr_1 (g_-), k, g_+)} (g_-, k')} (g'_+)$
\medskip
for all $g_-\in X^*, g_+, g'_+\in X$ with $\beta_{X^*} (g_-) = \alpha_X
(g_+),\, \beta_X (g_+) = \alpha_X (g'_+)$ and for all $k,k'\in H.$
\bigskip
(c) $$
\align
&\Psi^+ _ {(Pr_1 (g_-) Pr_1 (g'_-), kk', g_+)} (g_- g'_-)= \\
&\hskip 60 pt \Psi^+
_{(Pr_1 (g_-), k, \Psi^- _ {(g'_-, k')} (g_+))} (g_-)\, \Psi^+ _ {(Pr_1
(g'_-), k', g_+)} (g'_-)\\
\endalign
$$
\medskip
for all $g_+\in X$, $g_-, g'_-\in X^*$ with $\beta_{X^*} (g_-) =
\alpha_{X^*} (g'_-)$, $\beta_{X^*} (g'_-)= \alpha_X (g_+)$, and
for all $k,k'\in H.$}
\bigskip
Proof. We shall check (b) and (c).
\smallskip
(b) Let $g_-= (h, p, u), g_+= (p,g_1,q), g'_+= (q, g_2, r).$ Then
$$\eqalign {&\Psi^- _ {(g_-, k')} (g_+, g'_+)\cr
&= \big(Ad^*_{h^{-1}} p + I^* (u), \varphi^-_u (hg_1g_2k^{\prime-1}),
Ad^*_{k^{\prime-1}} r + I^* (\varphi^+ _{hg_1g_2k^{\prime-1}} (u))\big),\cr}$$
and
$$\eqalign { \Psi^- _ {(g_-, k)} \, &\Psi^- _ {(\Psi^+ _ {(Pr_1 (g_-), k,
g_+)} (g_-), k')} (g'_+)= \big(Ad^*_{h^{-1}} p + I^* (u),\cr
& \varphi^- _{u} (hg_1 k^{-1}) \varphi^- _{\varphi^+ _{hg_1k^{-1}} (u)}
(kg_2 {k'}^{-1}), Ad^*_{{k'}^{-1}} r + I^* (\varphi^+ _ {hg_1g_2
{k'}^{-1}} (u))\big).\cr}$$
Hence the assertion follows from (5.1.4).
\smallskip
(c) Let $g_-= (h,p, u),\, g'_-= (h', p', u')$ satisfy
$\beta_{X^*} (g_-)= \alpha_{X^*} (g'_-),$ i.e. $p=Ad^*_{{h'}^{-1} }p'+
I^* (u')$, and let $g_+= (p', g, q)$ so that $\beta_{X^*}
(g'_-)= \alpha_X (g_+).$
We have
$$\eqalign {&\Psi^+ _ {(Pr_1 (g_-)  Pr_1 (g'_-), kk', g_+)} (g_- g'_-)\cr
&= (kk', q, \varphi^+ _ {hh'g(kk')^{-1}} (u \varphi^+ _ {h^{-1}} (u')))\cr
&= (kk', q, \varphi^+ _ {\varphi^- _{\varphi^+ _ {h^{-1}} (u')} (hh'g
(kk')^{-1})} (u)\, \varphi^+ _ {h'g (kk')^{-1}} (u')).\cr}$$
On the other hand,
$$\eqalign {&\Psi^+ _ {(Pr_1 (g_-), k, \Psi^- _ {(g'_-, k')} (g_+))}
(g_-)\, \Psi^+ _ {(Pr_1 (g'_-), k', g_+)} (g'_-)\cr
&= (kk', q, \varphi^+ _ {h\varphi^- _{u'} (h'g k^{\prime-1}) k^{-1}} (u) \,
\varphi^+ _ {h'g (kk')^{-1}} (u')).\cr}$$
The assertion then follows from the calculation in the proof of Prop.
5.2.5.$\quad \pf$
\bigskip
Let ${\Cal  S}= X^*_e \ast_{J_+} X.$ Then ${\Cal  S}$ supports both the
left action groupoid structure $X^*_e \ltimes X$ and the
right action groupoid structure $X^* \rtimes X_e.$
\bigskip
{\bf Theorem 5.2.10} {\it If the horizontal structure on ${\Cal  S}$ is
$X^*_e \ltimes X$
and the vertical structure on ${\Cal  S}$ is $X^* \rtimes
X_e$, then $({\Cal  S}; X^*, X, \fh^*)$ is a double Lie groupoid.}
\bigskip
Proof. We have to show that the structure maps of the horizontal (resp.
vertical) structure on ${\Cal  S}$ are morphisms with respect to the
vertical (resp. horizontal) structure. We shall illustrate the role played
by the properties in Prop. 5.2.9 by checking that ${\tilde \alpha}_{\Cal
H}$ and ${\tilde {i}} _{\Cal  V}$ are groupoid morphisms. The rest of the
proof will be left to the interested reader.
\smallskip
$(i)\, {\tilde\alpha}_{\Cal  H}$ is a groupoid morphism.
\smallskip
$$\eqalign {&{\tilde\alpha}_{\Cal  H} \big( (g_-, k, g_+) \cdot ( \Psi^+_
{(Pr_1 (g_-), k , g_+)} (g_-), k', g'_+)\big)\cr
&= {\tilde\alpha}_{\Cal  H} (g_-, k', g_+ g'_+)\cr
&= \Psi^-_{(g_-, k')} (g_+ g'_+)\cr
&= \Psi^-_ {(g_-, k)} (g_+) \Psi^- _ {\Psi^+ _ {(Pr_1 (g_-), k, g_+)} (g_-),
k')} (g'_+) \, {\hbox { by Prop. 5.2.9 (b) }}\cr
&= {\tilde\alpha}_{\Cal  H} (g_-, k, g_+)\, {\tilde\alpha}_{\Cal  H} (\Psi^+
_ {(Pr_1 (g_-), k, g_+)} (g_-), k', g'_+).\cr}$$
\medskip
$(ii)\, {\tilde {i}}_ {\Cal  V}$ is a groupoid morphism.
\smallskip
$$\eqalign { & {\tilde {i}}_{\Cal  V} \big( (g_-, k, \Psi^-_ {(g'_-, k')}
(g'_+))\cdot (g'_-, k', g'_+)\big)\cr
&= {\tilde {i}}_ {\Cal  V} (g_- g'_-, kk', g'_+)\cr
&= (\Psi^+ _ {(Pr_1 (g_-g'_-), kk', g'_+)} (g_-g'_-), Pr_1 (g_-g'_-),
{g'}_+ ^{-1})\cr
&= (\Psi^+ _ {(Pr_1 (g_-), k, \Psi^- _ {(g'_-, k')} (g'_+))} (g_-) \, \Psi^+
_ {(Pr_1 (g'_-), k', g_+)} (g'_-), Pr_1 (g_-) Pr_1 (g'_-), {g'}_-^{-1}),\cr
&\hskip 15 pt { \hbox { by Thm 5.1.4 and Prop. 5.2.9 (c)}}\cr
&= (\Psi^- _ {(Pr_1 (g_-), k, \Psi^-_ {(g'_-, k')} ({g'}_+))} (g_-), Pr_1
(g'_-), \Psi^- _ {(g',k')} ({g'}_-)^{-1})\cdot\cr
&\hskip 30 pt (\Psi^+ _ {(Pr_1 (g'_-), k', g'_+)} (g'_-), Pr_1 (g'_-) ,
{g'}_+^{-1})\cr
&= {\tilde {i}}_{\Cal  V} (g_-, k, \Psi^-_ {(g'_-, k')} (g'_+)) \,
{\tilde{i}}_{\Cal V} (g'_-, k', g'_+).\quad \pf\cr}$$
\bigskip
To clarify the relation between the vacant double Lie groupoid $({\Cal
S}_{vac}; 
X^*, X, \fh^*)$ and the double Lie groupoid $({\Cal S}; X^*, X, \fh^*)$ in
the 
above theorem, we introduce the following definition
\bigskip
{\bf Definition 5.2.11}\hskip 5 pt {\sl Let $({\Cal S}_1; {\Cal H}_1, {\Cal
V}_1, 
P_1)$ be a double Lie groupoid. A double Lie subgroupoid of  $({\Cal S}_1;
{\Cal H}_1, {\Cal V}_1,
P_1)$ is a double Lie groupoid $({\Cal S}_2; {\Cal H}_2, {\Cal V}_2,
P_2)$ such that the Lie groupoids ${\Cal S}_2\rightrightarrows {\Cal
H}_2,\, {\Cal 
S}_2\rightrightarrows {\Cal V}_2,\, {\Cal H}_2\rightrightarrows P_2,\,
{\Cal V}_2\rightrightarrows P_2$ are respectively Lie subgroupoids of
 ${\Cal S}_1\rightrightarrows {\Cal
H}_1,\, {\Cal 
S}_1\rightrightarrows {\Cal V}_1,\, {\Cal H}_1\rightrightarrows P_1,\,
{\Cal V}_1\rightrightarrows P_1.$}
\bigskip
{\bf Corollary 5.2.12} \hskip 5 pt {\sl The vacant double Lie groupoid
$({\Cal S}_{vac}; 
X^*, X, \fh^*)$ associated with the matched pair $(X,X^*)$ is a double Lie
subgroupoid of the double Lie groupoid $({\Cal S}; X^*, X, \fh^*)$ in Thm
5.2.10.}
\medskip
Proof.  As the side groupoids of the two double groupoids Lie groupoids
are identical, it suffices to show that the Lie groupoids
${\Cal S}_{vac} \rightrightarrows X, \, {\Cal S}_{vac} \rightrightarrows
X^*$ are respectively Lie subgroupoids of ${\Cal S}\rightrightarrows X, \,
{\Cal S}\rightrightarrows
X^*$. For the horizontal structures, it suffices to observe that
$X^*\ltimes X\rightarrow X^*_e \ltimes X: (g_-, g_+)\mapsto (g_-, Pr_1
(g_-), g_+)$ is an injective immersion. The other case is similar.
 $\pf$
\bigskip
\bigskip
We now turn to the description of the symplectic properties of the double
Lie groupoid $({\Cal S}; X^*, X, \fh^*)$ of Thm. 5.2.10. (Recall that
$X=\fh^*\times G\times \fh^*$ is the dynamical groupoid for the constant
$r-$ matrix $-R.$)
\smallskip
To begin with, a simple computation (using Eqn.(5.1.4)) shows that the map
$${\eqalign {&\rho: (H\times \fh^*) \times (H\times \fh^*) \times G\times
G^*\longrightarrow X^*\rtimes X_e \cr
&(h,p,k,q,g,u)\mapsto ((h,p,\varphi^+_{g^{-1}} (u^{-1})^{-1}),\, (h,k,p,
h^{-1} \varphi^- _{u^{-1}}(g^{-1})^{-1} k, q))\cr}}\eqno (5.2.12) $$
is an isomorphism of groupoids; here the domain is the product groupoid of
Eqn (5.2.6) and the range is
the right action groupoid of Cor. 5.2.8.
\bigskip
Recall that the Poisson bracket of the symplectic groupoid
${\Cal S}'= (H\times \fh^*)\times (H\times \fh^*)^- \times (G\times
G^*)\rightrightarrows X^*$ of Eqn. (2.5.6)
is explicitely given by
$${\eqalign {&\{F, F'\}_{{\Cal S}'} (h,p,k,q, g,u)\cr
&= - <D'_1 F', \delta_1 F> + <D'_1 F, \delta_1 F' > - <p, [\delta_1
F, \delta_1 F']>\cr
&  + <D'_2 F', \delta_2 F> - <D'_2 F, \delta_2 F' > +<q,
[\delta_2 F, \delta_2 F']>\cr
& - <\partial F, \lambda_- (D F') (g)> + <\partial_* F, \lambda_+ (D'_* F')
(u)>\cr
&- <D'F', D_*F>+<D'F, D_*F'>,\cr}}
\eqno (5.2.13)$$
where the indices $1$ and $2$ indicate partial derivatives and left/right
gradients w.r.t. the appropriate factor in the first and second copies of
$H\times \fh^*$ and the index $*$ indicates partial derivative w.r.t. $G^*.$
Using the bijection $\rho$ we may transport this Poisson bracket
to ${\Cal S}$ by setting
$$\{F_1, F_2\}_{\Cal S} \circ \rho= \{F_1\circ \rho,
F_2\circ \rho\}_{{\Cal S}'}.$$
Since $\rho$ is a Lie groupoid isomorphism,  $(X^*\rtimes X_e, \{\, , \,
\}_{\Cal 
S})$ is a symplectic groupoid.
\smallskip
We now come to the main result of this section.
\bigskip
{\bf Theorem 5.2.13} \hskip 5 pt {\sl The double Lie groupoid $({\Cal S};
X^*, X, \fh^*)$ where ${\Cal S}$ is equipped with the Poisson bracket
$\{\, , \, \}_{{\Cal S}}$ is a symplectic double groupoid.}
\bigskip
In order to prove the theorem, it remains to show that $(X^*_e\ltimes
X, \{\, , \, \}_{\Cal S})$ is a symplectic groupoid. For this purpose, we
shall use the isomorphic image ${\Cal S}'\rightrightarrows X$ of
$X^*_e\ltimes X\rightrightarrows X$ under the map $\rho$ and the bracket
$\{\, ,\, \}_{{\Cal
S}'}.$ By direct computation, ${\Cal S}'\rightrightarrows X$
has target and source maps
$${\eqalign { & \alpha (h,p,k,q,g,u)\cr
&= (Ad^*_{h^{-1}} p + I^* (\varphi^+_{g^{-1}} (u^{-1})^{-1}), g , Ad^*
_{k^{-1}} q + I^* (u)),\cr
&\beta (h,p,k,q,g,u)\cr
&= (p, h^{-1} \varphi^- _ {u^{-1}} (g^{-1})^{-1} k, q),\cr}}\eqno
(5.2.14.a)$$
multiplication map
$${\eqalign {& m ((h_1, p_1, k_1, q_1, g_1, u_1), (h_2, p_2, k_2, q_2, g_2,
u_2))\cr
&= (h_1h_2, p_2, k_1k_2, q_2, g_1, u_1 \varphi^+_{k_1^{-1}}
(u_2)),\cr}}\eqno 
(5.2.14.b)$$
where
$$(p_1, h_1^{-1} \varphi^- _{u_1^{-1}} (g_1^{-1})^{-1} k_1, q_1)=
(Ad^*_{h_2^{-1}} p_2+ I^*(\varphi^+_{g_2^{-1}} (u_2^{-1}) ^{-1}), g_2,
Ad^*_{k_2^{-1}} q_2 + I^*
(u_2)),$$
and unit section
$$\epsilon (p,g,q)= (1,p, 1, q, g,1).\eqno (5.2.14.c)$$
\bigskip
We shall verify the conditions of the following proposition
of Libermann in \cite{L}, and show that the unique Poisson structure induced
on
the base $X$ indeed coincides with that of Thm 2.1.5.
\bigskip
{\bf Proposition 5.2.14} \hskip 5 pt {\sl Let $\Gamma\rightrightarrows P$ be
a Lie 
groupoid equipped with a symplectic form $\Omega.$ If $\Gamma$
is $\beta -$ connected,  the $\alpha-$ foliations and the $\beta-$
foliations are 
symplectically orthogonal, and $\epsilon (P)\subset \Gamma$ is Lagrangian,
then $\Gamma$ is a symplectic groupoid over $P.$}
\bigskip
In our case, that the $\beta -$ fibers of ${\Cal
S}'\rightrightarrows X$
are connected and $\epsilon (\fh^*)$ is
Lagrangian are easy to check and we shall leave the details to the reader.
In order to establish the other condition, we begin with two Propositions
which allow us to identify the target and source maps of (5.2.14.a) with
canonical projections of natural group actions.
\smallskip
Let $H\times (H\ltimes G^*)$ be the product of $H$ with the Lie group
$H\ltimes
G^*$ of Prop. 5.1.3 (b)
\bigskip
{\bf Proposition 5.2.15} \hskip 5 pt {\sl The left action of
${\Cal S}'\rightrightarrows X$ on
itself induces a left action of the group $H\times (H\ltimes G^*)$ on ${\Cal
S}'$ given by 
$$
\align
&(h',k',u')\cdot (h,p,k,q,g,u)\\
&= (h'h, p, k'k, q, h'\varphi^- _ {\varphi^+_{k'} (u')} (g^{-1})^{-1}
{k'}^{-1}, u'\varphi^+_{{k'}^{-1}} (u)).\\
\endalign
$$
Moreover, the canonical projection ${\Cal S}'\longrightarrow H\times
(H\ltimes
G^*)\backslash {\Cal S}'\simeq X$ coincides with $\beta$.}
\bigskip
Proof. Clearly, 
$$
\align
&(h',k',u')\cdot (h,p,k,q,g,u)\\
&= m(( h', p', k', q', g', u'), (h, p, k, q, g, u))\\
\endalign
$$
for unique $(p', q', g')\in \fh^*\times\fh^*\times G$ and it is easy to
show that this defines a left $H\times (H\ltimes G^*)$ action on ${\Cal
S}'$.
Now, if we identify each $H\times (H\ltimes G^*)-$ orbit with its unique
intersection with $\{1\}\times \fh^*\times \{1\}\times \fh^*\times
G\times \{1\}$, then $H\times (H\ltimes G^*)\backslash {\Cal S}'\simeq X$
and an easy calculation shows the projection map coincides with
$\beta .\quad \pf$
\bigskip
In a similar way, we have
\bigskip
{\bf Proposition 5.2.16} {\sl The right action of ${\Cal
S}'\rightrightarrows X$ on
itself induces a right action of the group $H\times (H\ltimes G^*)$ on
${\Cal S}'$ given by
$$\eqalign {& (h,p,k,q,g,u)\cdot (h',k',u')\cr
&\hskip 30 pt = (hh', Ad^*_{h'} p + I^*(\varphi^+ _{k^{-1}
\varphi^-_{u^{-1}} (g^{-1}) hh'}
({u'}^{-1})), kk',\cr
&\hskip 80 pt Ad^*_{k'} q + I^* (\varphi^+_{k'} ({u'}^{-1})), g,
u\varphi^+_{k^{-1}}
(u')).\cr}$$
Moreover, the canonical projection ${\Cal S}'\longrightarrow {\Cal S}'/
H\times
(H\ltimes G^*)\simeq X$ coincides with $\alpha .\quad \pf$}
\bigskip
\bigskip
If $\phi\in C^\infty (X)$, then it follows frow the above propositions
that $\alpha^* \phi$ is right $H\times (H\ltimes G^*)-$ invariant and
$\beta^* \phi$ is left $H\times (H\ltimes G^*)-$ invariant.  Conversely,
it is clear that right/left $H\times (H\ltimes G^*)$ invariant functions
on ${\Cal S}'$ are of the above form.

\bigskip
{\bf Lemma 5.2.17}
 \hskip 5 pt {\sl For $\phi\in C^\infty (X),$ we have}
\medskip
(a.1)\hskip 5 pt $ D'_1 (\alpha^*\phi)= ad^*_{\delta_1 (\alpha^*\phi)} p$
\smallskip
(a.2)\hskip 5 pt $D'_2(\alpha^*\phi)= ad^*_{\delta_2 (\alpha^* \phi)} q$
\smallskip
(a.3) \hskip 5 pt $D'_* (\alpha^* \phi)= Ad_{\phi^- _{u^{-1}} (g^{-1}) h}
\delta_1 
(\alpha^*\phi) + Ad_k \,\delta_2 (\alpha^* \phi),$
\medskip
(b.1) \hskip 5 pt $D_1 (\beta^* \phi) + \iota^* D (\beta^* \phi)= 0$
\smallskip
 (b.2) \hskip 5 pt $D_2 (\beta^* \phi)-\iota^* D' (\beta^* \phi)+ \iota^*
T_u 
l_{u^{-1}} \lambda^+ \big(D'_* (\beta^*\phi)\big)(u)=0$
\smallskip
(b.3) \hskip 5 pt $D_* (\beta^* \phi)- T_g l_{g^{-1}} \lambda^- \big(D
(\beta^* 
\phi)\big)(g)=0$
\medskip
{\sl Here all partial derivatives and left/right gradients are evaluated at
$(h,p,k,q,g,u)\in {\Cal S}'.$}
\medskip 
Proof. These are the infinitesimal versions of the invariance properties
of $\alpha^* \phi$ and $\beta^* \phi$ which can be obtained by using the
following basic formulae:
${d\over dt}_{\mid_0} \varphi^+_ {e^{tX}} (u)= \lambda^+ (X) (u)$,
${d\over dt}_{\mid_0} \varphi^- _ {e^{t\alpha}} (g)= \lambda^- (\alpha)
(g)$,
${d\over dt}_{\mid _0} \varphi^+ _g (e^{t\alpha}) = Ad^*_g \alpha$,
${d\over dt}_{\mid_0}\varphi^-_u (e^{tX}) = Ad^*_{u^{-1}} X.$
\hskip 10 pt $\pf$
\bigskip
{\bf Lemma 5.2.18}
\smallskip
{\sl (a) For all $g\in G, u\in G^*$ and $Z\in \fh,$ we have}
$$ Ad^*_{u^{-1}} Ad_{\varphi^- _{u^{-1}} (g^{-1})} Z= Ad_{g^{-1}} Z + T_g
l_{g^{-1}} \lambda^- \big(Tr_{\varphi^+_{g^{-1}}(u^{-1})}
\lambda^+ (Z) (\varphi^+ _{\varphi^- _{u^{-1}} (g^{-1})} (u))\big) (g).$$
\medskip
{\sl (b)}\hskip 5 pt
$\varphi^+_{g^{-1}} \circ r_{u^{-1}}\circ \varphi^+
_{\varphi^-_{u^{-1}} (g^{-1})^{-1}}= r_{\varphi^+ _{g^{-1}} (u^{-1})}$
\smallskip
{\sl for all } $g\in G, \, u\in G^*.$
\medskip
Proof.
\smallskip
(a)  
$$\eqalign { & Ad^* _{u^{-1}} \, Ad_{\varphi^- _ {u^{-1}} (g^{-1})} Z\cr
&= {d\over dt}_{\mid_0} \varphi^- _u (e^{t Ad_{\varphi^- _ {u^{-1}}
(g^{-1})} 
Z})\cr
&= {d\over dt}_{\mid_0} g^{-1} e^{tZ} \varphi^- _ {\varphi^+_{e^{tZ}}
(\varphi^+ _{\varphi^-_{u^{-1}} (g^{-1})} (u)) \varphi^+_{g^{-1}}
(u^{-1})} (g),\cr
& {\hbox { by repeated application of Eqn (5.1.4) and the triviality }}\cr
&{\hbox { of the action of }} G^* {\hbox { on }} H\cr
&= Ad_{g^{-1}} Z + T_gl_{g^{-1}} \lambda^- \big( Tr_{\varphi^+_{g^{-1}}
(u^{-1})} \lambda^+ (Z) (\varphi^+_{\varphi^- _{u^{-1}} (g^{-1})} (u))\big)
(g)\cr}$$
where we have used the formulae in the proof of Lemma 5.2.17.
\medskip
(b) 
$$\eqalign {&\varphi^+_{g^{-1}} \circ r_{u^{-1}}\circ \varphi^+
_{\varphi^-_{u^{-1}} (g^{-1})^{-1}} (u')\cr
&= \varphi^+ _{g^{-1}} (\varphi^+ _{\varphi^-_ {u^{-1}} (g^{-1})^{-1}} (u')
\, u^{-1})\cr
&= \varphi^+ _{\varphi^- _{u^{-1}} (g^{-1})} (\varphi^+ _{\varphi^- _
{u^{-1}} (g^{-1})^{-1}} (u')) \varphi^+ _{g^{-1}} (u^{-1})\cr
& {\hbox { by Eqn (5.1.4)}}\cr
&= u' \varphi^+ _{g^{-1}} (u^{-1}).\cr}$$
Hence the assertion. $\pf$
\bigskip
{\bf Proposition 5.2.19} (Polarity condition) {\sl For all
$\phi, \psi\in C^\infty (X)$, we have}
$$\{\alpha^*\phi, \beta^* \psi\}_{{\Cal
S}'}= 0.$$
\medskip
Proof. Let $\hat\phi= \alpha^* \phi,\, \hat\psi= \beta^* \psi.$
By invoking the identities (a.1), (a.2) and (b.3) of Lemma 5.2.17, 
we have
$$\eqalign {&\{\hat\phi, \hat \psi\}_{{\Cal S}'}(h,p,k,q,g,u)\cr
&=-<D_1\hat\psi, Ad_h \delta_1 \hat\phi>+<D_2\hat\psi, Ad_k \delta_2
\hat\phi>\cr
&+<d_* \hat\phi, \lambda^+ (D'_* \hat\psi) (u)>-<D'\hat\psi, D_*
\hat\phi>.\cr}$$
Next,using (b.1), (b.2) and (b.3) of the same lemma successively gives
$$\eqalign { & \{\hat\phi, \hat\psi\}_{{\Cal S}'} (h,p,k,q,gu)\cr
&= <D'\hat\psi, Ad_{g^{-1}h} \delta_1\hat\phi + Ad_k \delta_2 \hat\phi -
D_* \hat\phi\cr
&+ T_g l_{g^{-1}} \lambda^- (Ad^*_{g^{-1}} T_u r_{u^{-1}} \lambda^+
(Ad_{\varphi^- _ {u^{-1}} (g^{-1} h) } \delta_1\hat\phi) (u)) (g)>\cr}$$
Now using (a.3) together with $Ad^*_{u^{-1}} Z= Z,\, \forall Z\in \fh$, we
obtain
$$\eqalign { &\{\hat\phi, \hat\psi\}_{{\Cal S}'} (h,p,k,q,g,u)\cr
&= <D'\hat\psi, Ad_{g^{-1}h} \delta_1 \hat\phi - Ad^*_{u^{-1}} Ad_{\varphi^-
_{u^{-1}} (g^{-1})h} \delta_1\hat\phi\cr
&+ T_g l_{g^{-1}} \lambda^- (Ad^*_{g^{-1}} T_u r_{u^{-1}} \lambda^+
(Ad_{\varphi^- _ {u^{-1}} (g^{-1})h}\delta_1\hat\phi) (u)) (g)>\cr
&= <D'\hat\psi, T_g l_{g^{-1}} \lambda^- (Ad^*_{g^{-1}} T_u r_{u^{-1}}
\lambda^+ (Ad_{\varphi^-_{u^{-1}} (g^{-1})h} \delta_1\hat\phi)(u))\cr
&-Tr_{\varphi^+_{g^{-1}} (u^{-1})}\lambda^+ (Ad_h \delta_1 \hat\phi)
(\varphi^+ _{\varphi^- _{u^{-1}} (g^{-1})}(u)))(g)>\cr}$$
where in the last equality we have used Lemma 5.2.18 (a).
\smallskip
But 
$$\eqalign {&\lambda_+ (Ad_{\varphi^- _{u^{-1}} (g^{-1})h}
\delta_1\hat\phi) (u)\cr
&= T\varphi^+_{\varphi^- _ {u^{-1}} (g^{-1})^{-1}} \lambda_+ (Ad_h
\delta_1\hat\phi) (\varphi^+ _{g^{-1}} (u^{-1})^{-1}),\cr}$$
hence the assertion that $\{\hat\phi, \hat\psi\}_{{\Cal S}'}=0$ now follows
from Lemma 5.18 (b). \hskip 5 pt $\pf$
\bigskip
\bigskip
{\bf Lemma 5.2.20} {\it For $\phi\in C^\infty (X)$, we have
\bigskip
(a) $\, \delta_1 (\alpha^* \phi)= Ad_{h^{-1}} \delta_1 \phi (\alpha (s)),$
\medskip
(b) $\, \delta_2 (\alpha^* \phi)= Ad_{k^{-1}}\delta_2  \phi (\alpha (s)),$
\medskip
(c) $\, D(\alpha^*\phi)= D\phi (\alpha (s))- Tl_{\varphi^+_{g^{-1}}
(u^{-1})^{-1}} \lambda^+ \big(\delta_1 \phi (\alpha (s))\big)
(\varphi^+_{g^{-1}}
(u^{-1})),$
\medskip
(d) $\, D_* (\alpha^* \phi)= \delta_2 \phi (\alpha (s))+ Ad^*_{u^{-1}}
Ad_{\varphi^- _ 
{u^{-1}} (g^{-1})} \delta_1 \phi (\alpha (s))$
\medskip
\hskip 60 pt = $\delta_2 \phi (\alpha (s))+ Ad_{g^{-1}} \delta_1 \phi
(\alpha (s))$
\smallskip
\hskip 25 pt $- T_g l_{g^{-1}} \lambda^- \big(Tl_{\varphi^+ _{g^{-1}}
(u^{-1})^{-1}} \lambda^+ (\delta_1 \phi (\alpha (s))) (\varphi^+_{g^{-1}}
(u^{-1}))\big) (g)\hskip 15 pt $( second form).
\medskip
where all the partial derivatives and right gradients of $\alpha^*\phi$
are evaluated at
$s=(h,p,k,q,g,u)\in {\Cal S}'.$}
\bigskip
Proof. We shall establish the formulae in (d), leaving the other parts to
the interested reader. For $\gamma\in \fg^*,$ we have
$$
\align
&<D'_*(\alpha^*\phi), \gamma>\\
&= {d\over dt}_{\mid_0} \, \phi \big(Ad^*_{h^{-1}} p -
I^*(\varphi^+_{g^{-1}}
(e^{-t\gamma}u^{-1})), g, Ad^*_{k^{-1}} q+ I^* (u e^{t\gamma})\big)\\
&=<\delta_1\phi (\alpha (s)), -TI^* {d\over dt}_{\mid_0}
\varphi^+_{g^{-1}} (e^{-t\gamma}u^{-1})> + <\delta_2 \phi, \iota^*
\gamma>.\\
\endalign
$$
Since
$${d\over dt}_{\mid_0} \varphi^+_{g^{-1}} (e^{-t\gamma} u^{-1})= - T_1
r_{\varphi^+_{g^{-1}} (u^{-1})}\, Ad^*_{\varphi^-_{u^{-1}}
(g^{-1})}\gamma,$$
and
$$T_1 (I^*\circ r_{\varphi^+_{g^{-1}} (u^{-1})})= T_1 I^*= \iota^*,$$
it follows that
$$D'_* (\alpha^* \phi)= \delta_2\phi (\alpha (s))+ Ad_{\varphi^- _{u^{-1}}
(g^{-1})} \delta_1\phi (\alpha (s)),$$
and this gives the first form of $D_* (\alpha^*\phi).$ To obtain the
second form
of $D_* (\alpha^* \phi)$ from the first one, simply apply Lemma 5.2.18
(a) and the fact that
$$
\align
&\lambda^+ (Z) (\varphi^+_{g^{-1}} (u^{-1})^{-1})\\
&= - T\big( l_{\varphi^+_{g^{-1}} (u^{-1})^{-1}}\circ r_{\varphi^+_{g^{-1}}
(u^{-1})^{-1}}\big)
\lambda^+ (Z) (\varphi^+_{g^{-1}} (u^{-1})),\quad Z\in \fh. \quad
\pf\\
\endalign
$$
\bigskip
\bigskip
{\bf Lemma 5.2.21}  {\it For all $g\in G, u\in G^*,$ and $Z\in \fh$,
we have
\bigskip
(a) $\, T_u I^* \lambda^+(Z) (u)= ad^*_Z I^* (u),$
\medskip
(b) $Ad^*_g\, Tl_{\varphi^+_{g^{-1}} (u^{-1})^{-1}} \lambda^+ (Z)
(\varphi^+_{g^{-1}} (u^{-1}))$
\medskip
\hskip 15 pt $= -T_u r_{u^{-1}} \lambda^+ (Ad_{\varphi^-_{u^{-1}}
(g^{-1})}Z) (u),$
\medskip
(c) $Tl_{\varphi^+_{g^{-1}} (u^{-1})^{-1}} \lambda^+ (Ad_g Z)
(\varphi^+_{g^{-1}} (u^{-1}))$
\medskip
\hskip 15 pt $= - Ad_{\varphi^+_{g^{-1}} (u^{-1})^{-1}}
Ad^*_{\varphi^-_{u^{-1}} (g^{-1})} T_u l_{u^{-1}} \lambda^+ (Z) (u).$}
\bigskip
\bigskip
Proof. (a)  This is the infinitesimal version of the $H-$ equivariance of
the map $I^*.$
\medskip
(b) 
$$
\align
&\lambda^+ (Ad_ {\varphi^-_{u^{-1}} (g^{-1})} Z)(u)\\
&= {d\over dt}_{\mid_0} \varphi^+_ {e^{tAd_{\varphi^-_{u^{-1}} (g^{-1})}Z}}
(u)\\
&=T\varphi^+_ {\varphi^-_ {u^{-1}} (g^{-1})^{-1}}\, \lambda^+ (Z)
(\varphi^+_{g^{-1}} (u^{-1})^{-1})\\
&= - T \varphi^+_ {\varphi^-_ {u^{-1}} (g^{-1})^{-1}} Tr_{\varphi^+_{g^{-1}}
(u^{-1})^{-1}} Tl_{\varphi^+_{g^{-1}} (u^{-1})^{-1}} \lambda^+ (Z)
(\varphi^+_{g^{-1}} (u^{-1}))\\
& ({\hbox {where we have used }} Ad^*_u Z= Z {\hbox { for }} u\in G^*, Z\in
\fh).\\
\endalign
$$
Now $Ad^*_g= T_1 \varphi^+_g$ and by a straightforward calculation using
Eqn.(5.1.4), we find 
$$\varphi^+_g \circ r_{\varphi^+_{g^{-1}}
(u^{-1}) }\circ \varphi^+_ {\varphi^-_ {u^{-1}} (g^{-1})}=
r_{u^{-1}}.$$
Hence the assertion follows.
\medskip
(c) We have
$$
\align
&\lambda^+ (Ad_g Z) (\varphi^+_{g^{-1}} (u^{-1}))\\
&={d\over dt}_{\mid_0} \varphi^+_ {e^{t Ad_g Z}} (\varphi^+_{g^{-1}}
(u^{-1}))\\
&= {d\over dt}_{\mid_0} \varphi^+_ {g^{-1}} \big( (\varphi^+_ {\varphi^-
_ {u^{-1}} (e^{tZ})} (u))^{-1}\big)\\
&= {d\over dt}_{\mid_0} \varphi^+_{g^{-1}} (\varphi^+_ {e^{tZ}} (u)^{-1})\\
& {\hbox { by the triviality of the action of }} G^* {\hbox { on }} H\\
&= - T_{u^{-1}} \varphi^+_ {g^{-1}} T_1 r_{u^{-1}} T_u l_{u^{-1}} \lambda^+
(Z) 
(u).\\
\endalign
$$
On the other hand, for ${\tilde u}\in G^*$, we find
$$l_{\varphi^+_{g^{-1}} (u^{-1})^{-1}} \circ \varphi^+_{g^{-1}}\circ
r_{u^{-1}} ({\tilde u})= \varphi^+_ {g^{-1}} (u^{-1})^{-1} \varphi^+_
{\varphi^- _ {u^{-1}} (g^{-1})} ({\tilde u}) \varphi^+_ {g^{-1}} (u^{-1})$$
from which we deduce that
$$T_1\big( l_{\varphi^+_{g^{-1}} (u^{-1})^{-1}} \circ
\varphi^+_{g^{-1}}\circ
r_{u^{-1}}\big)= Ad_{\varphi^+_{g^{-1}} (u^{-1})^{-1}} Ad^*_{\varphi^- _
{u^{-1}} (g^{-1})}.$$
The assertion is now clear.\hskip 5 pt $\pf$
\bigskip
\bigskip
{\bf Proposition 5.2.22} {\it The map $\alpha: ({\Cal S}', \{\,,\,\}_{{\Cal
S}'})\longrightarrow (X, \{\, , \, \}_X)$ is a Poisson map.}
\bigskip
Proof. We want to show
$$\{\alpha^* \phi, \alpha^* \psi\}_{{\Cal S}'}=
\alpha^* \{\phi, \psi\}_X$$
 for all $\phi, \psi\in C^\infty (X).$ We shall
evaluate $\{\alpha^* \phi, \alpha^* \psi\}_{{\Cal S}'}$ at
$s=(h,p,k,q,g,u)\in {\Cal S}'$ and set $x=\alpha (s).$
\smallskip
By using Lemma 5.2.17 (a.1), (a.2) and Lemma 5.2.20 (a), (b), we have
$$
\align
&\{\alpha^* \phi, \alpha^* \psi\}_ {{\Cal S}'} (s)\\
&=<Ad^*_ {h^{-1}}p, [\delta_1\phi (x), \delta_1 \psi (x)]> - <Ad^*_
{k^{-1}}q, [\delta_2\phi (x), \delta_2 \psi (x)]>\\
&-<D' (\alpha^* \phi), T_g l_{g^{-1}} \lambda^- (D(\alpha^*\psi))(g) -
D_* (\alpha^*\psi)>\\
&+ <D_* (\alpha^* \phi), T_ur_{u^{-1}} \lambda^+ (D'_* (\alpha^*\psi)) (u)-
D'(\alpha^*\psi)>\\
\endalign
$$
Now, it follows from Lemma 5.2.20 (c) and (d) (second form) that
$$
\align
&T_g l_{g^{-1}} \lambda^- (D (\alpha^* \psi))(g)- D_* (\alpha^*\psi)\\
&= T_g l_{g^{-1}} \lambda^- (D\psi (x)) (g)- \delta_2 \psi (x) - Ad_{g^{-1}}
\delta_1 \psi (x).\\
\endalign
$$
On the other hand, by using Lemma 5.2.20 (c), (d) (first form) and Lemma
5.2.21 (b), we obtain
$$
\align
& T_u r_{u^{-1}} \lambda^+ (D'_* (\alpha^* \psi))(u)- D'(\alpha^*\psi)\\
&= T_u r_{u^{-1}} \lambda^+ (\delta_2 \psi (x)) (u) - D'\psi (x).\\
\endalign
$$
Consequently,
$$
\align
&-<D' (\alpha^* \phi), T_g l_{g^{-1}} \lambda^- (D(\alpha^*\psi))(g) -
D_* (\alpha^*\psi)>\\
&+ <D_* (\alpha^* \phi), T_ur_{u^{-1}} \lambda^+ (D'_* (\alpha^*\psi)) (u)-
D'(\alpha^*\psi)>\\
&= <D'\phi (x), \delta_2\psi (x)>+ <D\phi (x), \delta_1 \psi (x)>\\
&- <D'\psi (x), \delta_2 \phi (x)> - <D\psi (x), \delta_1 \phi (x)>\\
&-<D'\phi (x), T_gl_{g^{-1}} \lambda^- (D\psi (x))(g)>\\
&+T_1+T_2+T_3+T_4,\\
\endalign
$$
where
$$
\align
&T_1 =<Ad^*_g\, Tl_ {\varphi^+_ {g^{-1}} (u^{-1})^{-1}} \lambda^+ (\delta_1
\phi (x)) ({\varphi^+_ {g^{-1}} (u^{-1})}), T_g l_{g^{-1}} \lambda^-
(D\psi (x)) (g)>\\
&+< T_g l_{g^{-1}} \lambda^- ( Tl_{\varphi^+_ {g^{-1}} (u^{-1})^{-1}}
\lambda^+ (\delta_1 \phi (x)) ({\varphi^+_ {g^{-1}} (u^{-1})}) (g), D'\psi
(x)>,\\
&\\
&T_2= -<\iota^* Tl_{\varphi^+_ {g^{-1}} (u^{-1})^{-1}} \lambda^+ (\delta_1
\phi (x)) ({\varphi^+_ {g^{-1}} (u^{-1})}), \delta_1 \psi (x)>,\\
&\\
&T_3= <\delta_1 \phi (x), \iota^* T_u r_{u^{-1}} \lambda^+ (\delta_2 \psi
(x)) (u)>,\\
&T_4= <Ad^*_ {u^{-1}} Ad_ {\varphi^- _ {u^{-1}} (g^{-1})} \delta_1 \phi
(x), T_u r_{u^{-1}} \lambda^+ (\delta_2 \psi (x)) (u)>\\
&-<Ad^*_g\, Tl_{\varphi^+_ {g^{-1}} (u^{-1})^{-1}} \lambda^+ (\delta_1 \phi
(x))({\varphi^+_ {g^{-1}} (u^{-1})}), \delta_2 \psi (x)>.\\
\endalign
$$
Using the relation
$$<\lambda^- (\gamma) (g), v>= - <\gamma, T_g r_{g^{-1}} \lambda^- (T^*_1
r_g v) (g)>,$$
it is immediate that $T_1=0.$ For the term $T_2,$ note that
$$\iota^* Tl_{\varphi^+_ {g^{-1}} (u^{-1})^{-1}}= T_{\varphi^+_ {g^{-1}}
(u^{-1})} I^*.$$ Hence it follows from Lemma 5.2.21 (a) that
$$
\align
& T_2= - <ad^*_ {\delta_1 \phi (x)} I^* (\varphi^+_ {g^{-1}} (u^{-1})),
\delta_1 \psi (x)>\\
&= - <I^*({\varphi^+_ {g^{-1}} (u^{-1})}), [\delta_1 \phi (x), \delta_1
\psi (x)]>.\\
\endalign
$$
Similarly we have
$$T_3= - <I^* (u), [\delta_2 \phi (x), \delta_2 \psi (x)]>.$$
Assembling the calculations, we find
$$\{\alpha^*\phi, \alpha^* \psi\}_ {{\Cal S}'} (s)= \{\phi, \psi\}_X (x) +
T_4.$$
Hence it remains to show that $T_4=0.$ To do so, we
invoke the relation
$$<\lambda^+ (X) (u), \omega>= - <X, T_u l_{u^{-1}} \lambda^+ (T_1^* l_u
\omega) (u)>$$
and Lemma 5.2.21 (c) to rewrite $T_4$ as
$$
\align
&T_4= <\delta_1 \phi (x), -\iota^* Ad_{\varphi^+_ {g^{-1}} (u^{-1})^{-1}}
Ad^*_ {\varphi^-_ {u^{-1}} (g^{-1})} T_u l_{u^{-1}} \lambda^+ (\delta_2
\psi (x)) (u)\\
&\hskip 20 pt + \iota^* Ad^*_ {\varphi^- _ {u^{-1}} (g^{-1})} T_u
l_{u^{-1}} \lambda^+ (\delta_2 \psi (x)) (u)>.\\
\endalign
$$
But from the triviality of the $Ad^*_{G^*}$ action on $\fh$, it follows that
$$\iota^* Ad^*_ {\varphi^- _{u^{-1}} (g^{-1})}= \iota^* Ad_{\varphi^+_
{g^{-1}} (u^{-1})^{-1}}\, Ad^*_ {\varphi^- _{u^{-1}} (g^{-1})}.$$
Therefore, $ T_4=0.\quad \pf$
\bigskip
\bigskip
Combining the above Proposition with Proposition 5.2.16, we have
\bigskip
{\bf Corollary 5.2.23} {\it The right action of $H\times (H\ltimes G^*)$
on ${\Cal S}'$ in Proposition 5.2.16 is admissible, i.e. functions in
$C^\infty ({\Cal S}')$ invariant under the action form a Lie subalgebra of
$C^\infty ({\Cal S}').$ Furthermore, the quotient Poisson structure on
$X\simeq {\Cal S}'/ H\times (H\ltimes G^*) $ coincides with $\{\, , \,
\}_X.\quad \pf$}
\bigskip
We shall skip the proof of the next two lemmas.
\bigskip
{\bf Lemma 5.2.24} {\it For $\phi\in C^\infty (X)$, we have
\medskip
(a) $\, \delta_1 (\beta^*\phi)= \delta_1 \phi (\beta (s)),$
\medskip
(b) $\, \delta_2 (\beta^* \phi)= \delta_2\phi (\beta (s)),$
\medskip
(c) $\, D'_1 (\beta^*\phi)= -\iota^* D\phi (\beta (s)),$
\medskip
(d) $\, D'_2 (\beta^* \phi)= \iota^* D'\phi (\beta (s)),$
\medskip
(e)$\, D (\beta^* \phi)= Ad_ {\varphi^+_ {g^{-1}} (u^{-1})^{-1}}
Ad^*_{h^{-1}} D\phi (\beta (s)),$
\medskip
(f)$\, D'_* (\beta^* \phi)= - T r_{\varphi^-_ {u^{-1}} (g^{-1})^{-1}}
\lambda^- (Ad^* _{k^{-1}} D'\phi (\beta (s))) (\varphi^-_ {u^{-1}}
(g^{-1})),$
\medskip
where all the partial derivatives and left/right gradients of $\beta^*\phi$
are evaluated at
$s= (h,p,k,q, g,u)\in {\Cal S}'.\quad \pf$}
\bigskip
{\bf Lemma 5.2.25} {\it For all $g\in G, u\in G^*,$ and $\gamma\in \fg^*,$
we have
$$
\align
&Ad^*_g\, Ad_ {\varphi^+_ {g^{-1}} (u^{-1})^{-1}} \gamma= Ad_u\, Ad^*_
{\varphi^- _ {u^{-1}} (g^{-1})^{-1}} \gamma\\
\hskip 5 pt &-T_1 r_{u^{-1}} \lambda^+ (Tr_ {\varphi^-_ {u^{-1}}
(g^{-1})^{-1}} 
\lambda^- (Ad^*_{\varphi^-_ {u^{-1}} (g^{-1})^{-1}}
\gamma) (\varphi^-_ {u^{-1}} (g^{-1}))) (u).\, \pf\\
\endalign
$$}
\bigskip
\bigskip
{\bf Proposition 5.2.26} {\it The map $\beta: ({\Cal S}', \{\,,\,\}_{{\Cal
S}'})\longrightarrow (X, \{\, , \, \}_X)$ is an anti-Poisson map.}
\bigskip
Proof.  Let $\phi, \psi\in C^\infty (X), \, s= (h,p,k,q,g,u)\in {\Cal S}'$
and denote $\beta (s)$ by $x.$ From Lemma 5.2.24 (a) - (d) and Lemma 5.2.17
(b.1), we immediately have
$$
\align
&\{\beta^* \phi, \beta^* \psi\}_ {{\Cal S}'} (s)\\
&= <\iota^* D\psi (x), \delta_1 \phi (x)>- <\iota^* D\phi (x), \delta_1
\psi (x)> -<p, [\delta_1 \phi (x), \delta_1 \psi (x)]>\\
&+<\iota^* D'\psi (x), \delta_2 \phi (x)>- <\iota^* D'\phi (x), \delta_2
\psi (x)> + <q, [\delta_2 \phi (x), \delta_2 \psi (x)]>\\
&+ <D'_* (\beta^* \phi), T_ul_{u^{-1}} \lambda^+ (D'_* (\beta^* \psi))
(u) - Ad_ {u^{-1}} D' (\beta^* \psi)>.\\
\endalign
$$
Hence we have to show that
$$ 
\align
& <D'_* (\beta^* \phi), T_ul_{u^{-1}} \lambda^+ (D'_* (\beta^* \psi))
(u) - Ad_ {u^{-1}} D' (\beta^* \psi)>\\
&= <d_G\phi (x), \Pi^\#_G (h^{-1} \varphi^-_ {u^{-1}} (g^{-1})^{-1} k) d_G
\psi (x)>.\\
\endalign
$$
To do so, we apply Lemma 5.2.24 (e), (f) and Lemma 5.2.25, this yields
$$
\align
& T_u l_ {u^{-1}} \lambda^+ (D'_* (\beta^* \psi)) (u) - Ad_{u^{-1}}
D'(\beta^* \psi)\\
&= - Ad^* _ {\varphi^-_ {u^{-1}} (g^{-1})^{-1}} Ad^* _ {h^{-1}} D\psi
(x).\\
\endalign
$$
Consequently,
$$
\align
& <D'_* (\beta^* \phi), T_ul_{u^{-1}} \lambda^+ (D'_* (\beta^* \psi))
(u) - Ad_ {u^{-1}} D' (\beta^* \psi)>\\
&= <Tr_{\varphi^- _ {u^{-1}} (g^{-1})} \lambda^- (Ad^*_{k^{-1}} D'\phi
(x))({\varphi^- _ {u^{-1}} (g^{-1})}), Ad^*_{\varphi^- _ {u^{-1}}
(g^{-1})^{-1}}  Ad^*_ {h^{-1}} D\psi (x)>\\
&= <d_G \psi (x), T_{\varphi^- _ {u^{-1}} (g^{-1})}  (l_{h^{-1} {\varphi^- _
{u^{-1}} (g^{-1})^{-1}} }\circ r_{{\varphi^- _ {u^{-1}} (g^{-1})^{-1}}
k})\, \Pi^\#_G ({\varphi^- _ {u^{-1}} (g^{-1})})\\
&\hskip 30 pt T^*_{\varphi^- _ {u^{-1}} (g^{-1})} (l_{h^{-1} {\varphi^- _
{u^{-1}} (g^{-1})^{-1}} }\circ r_{{\varphi^- _ {u^{-1}} (g^{-1})^{-1}}k})
\, d_G \phi (x)>\\
&= - <d_G \psi (x), \Pi^\#_G (h^{-1} \varphi^-_ {u^{-1}} (g^{-1})^{-1}k)
d_G \phi (x)>,\\
\endalign
$$
where in the last step, we have used the fact that $\Pi_G$ is
multiplicative
and that $\Pi_G$ vanishes on $H.$  This completes the proof. $\quad \pf$
 \bigskip
As a consequence of the Prop. 5.2.26 and Prop. 5.2.15, we obtain
\bigskip
{\bf Corollary 5.2.27} {\it The left action of $H\times (H\ltimes G^*)$
on ${\Cal S}'$ in Proposition 5.2.15 is admissible (see Cor.5.2.23).
Furthermore, the quotient Poisson structure on $X\simeq H\times (H\ltimes
G^*)\backslash {\Cal S}'$ coincides with $-\{\, , \, \}_X.\quad \pf$
\bigskip}
This completes the proof that ${\Cal S}'\rightrightarrows X$ is a
symplectic groupoid.
\bigskip
We now turn to the description of the symplectic foliation of $X.$
\bigskip
Equip $H\times (H\ltimes G^*)$ with the product of the trivial Poisson Lie
group structure on $H$ and the Poisson Lie group structure of Prop. 5.1.3
on $H\ltimes G^*.$
It is easy to see that the groupoid action $\Psi^-$  of Eqn. (5.2.8)
restricts to a
group action
$$\tilde\Psi^-: H\times (H\ltimes G^*)\times X\longrightarrow X,$$
given by
$$((k,(h,u)), (p,g,q))\mapsto \Psi^- _{(h,p,u,k)} (p,g,q).$$
\bigskip
{\bf Theorem 5.2.28}\hskip 5 pt
\smallskip
{\it (a) \hskip 5 pt $\tilde\Psi^-$ is a left Poisson Lie group action.
\medskip
(b) The symplectic leaf ${\Cal L} _{(p,g,q)}$ in $(X, \{\, , \, \}_X)$
passing through the point
$(p,g,q)$ is the orbit of $(p,g,q)$ under the action $\tilde\Psi^-$ i.e.
$$\eqalign {{\Cal L}_{(p,g,q)}= \{ (Ad^* _{h^{-1}} p + I^* (u), &\varphi^-
_u 
(hgk^{-1}), Ad^*_{k^{-1}} q + I^* (\varphi^+_{hgk^{-1}} (u)))\, \cr
&\mid (k, (h,u))\in H\times (H\ltimes G^*)\}.\cr}$$}
\bigskip
Proof.\hskip 5 pt (a) Equip $X$ with the Poisson bracket of Thm. 2.1.5 and
$H\times (H\ltimes G^*)$ with the Poisson Lie bracket
$$\{\phi, \psi\} (k, (h,u)) = \partial_* \phi (u) \, \Pi_* (u) \partial_*
\psi (u)$$
of Prop. 5.1.3.   We have to show that the action

$$
\align
&{\tilde \Psi}^-:\hskip 10 pt  H\times (H\ltimes G^*)\hskip 5 pt \times
\hskip 5 pt  X\hskip 10 pt
\longrightarrow\hskip 10 pt  X\\
&((k, (h,u)), (p,g,q))\mapsto (Ad^* _{h^{-1}} p + I^* (u), \varphi^-
_u (hgk^{-1}),  Ad^*_{k^{-1}} q + I^* (\varphi^+_{hgk^{-1}} (u)))\\
\endalign
$$
satisfies
$$\{f\circ {\tilde \Psi}^-, f'\circ {\tilde \Psi}^-\}_ {H\times (H\ltimes
G^*)\times X} = \{f, f'\}_X \circ {\tilde \Psi}^-.\eqno (P)$$
To begin with, a direct calculation making use of the equivariance of
$I^*$, Eqn. (5.1.4), and the triviality of the action of $G^*$ on $H$ yields
the following expressions for the partial derivatives of $f\circ {\tilde
\Psi}^-$ 
at $((k, (h,u), p,g,q)\in H\times (H\ltimes G^*) \times X:$
$$
\align
&\partial_* (f\circ {\tilde \Psi}^-) = T^*_u r_{u^{-1}}\big( \iota \delta_1
f - \Pi^r (\varphi^-
_u (hgk^{-1})) Df +  Ad_{\varphi^-_u (hgk^{-1})} \iota \delta _2 f\big),\\
& \delta_1 (f\circ {\tilde \Psi}^-)= Ad_{h^{-1}} \delta_1 f,\\
&\delta_2 (f\circ {\tilde \Psi}^-)= Ad_{k^{-1}} \delta_2 f,\\
&\partial\, (f\circ {\tilde \Psi}^-)= T^*_g l_{g^{-1}} \big(
Ad_{\varphi^+_{hg} 
(u)^{-1}} Ad^* _k D'f + \Pi_*^l (\varphi^+_{hg} (u)) \iota Ad_{k^{-1}}
\delta_2 f\big),\\
\endalign
$$
where $\Pi^r$ (resp. $\Pi_*^l$) stands for the Poisson tensor of $G$
(resp. $G^*$) in the right (resp. left) invariant
frame. 
\smallskip
We shall restrict ourselves to an outline of the main steps of the
calculation of $lhs (P) - rhs (P)$. We use the shorthand notation $d_{ij} =
(lhs)_{ij} (P) - (rhs)_{ij}
(P),\, i,j\in \{1,2,\star\}$,
where 1 (resp. 2) stands for $f= p_1^*\phi$ (resp. $f=p_2^* \phi$),
$\,\phi\in C^\infty (\fh^*)$, and $\star$ stands for $f= p_G^* \psi, \,
\psi\in C^\infty (G).$ Thus for example
$$d_{1*}= \{p_1^* \phi\circ {\tilde
\Psi}^-, p_G^* \psi\circ {\tilde \Psi}^-\}_ {H\times (H\ltimes
G^*)\times X} - \{p_1^*\phi, p_G^* \psi\}_X \circ {\tilde \Psi}^-.$$
\smallskip
We now compute $d_{ij}$ for the various cases.
\smallskip
(1) $d_{11}= <\iota\delta\phi, \Pi_*^r (u) \iota \delta\psi> - <I^*
(u) , [\delta\phi, \delta\psi]>= 0$ by Lemma 5.2.21 (a).
\medskip
(2) 
$$
\align
 d_{22} &= -2 <\iota Ad_{k^{-1}} \delta \phi, \Pi_*^l
(\varphi^+_{hg}(u)) \iota Ad_{k^{-1}} \delta \psi>+ <\iota \delta \phi,
\Pi_*^l (\varphi^+ _{hgk^{-1}} (u))\iota \delta \psi>\\
& + < Ad_{\varphi^- _u (hgk^{-1})} \iota \delta \phi, \Pi^r_* (u)
Ad_{\varphi^-_u (hgk^{-1})} \iota\delta \psi>\\
&+ <\iota Ad_{k^{-1}} \delta \phi , \Pi^l _* (\varphi^+_{hg} (u)) \Pi^l (g)
\Pi^l _* (\varphi^+_{hg} (u)) \iota Ad_{k^{-1}} \delta
\psi>.\\
\endalign
$$
Using the fact that the action $H\times G^*\longrightarrow G^*: (k,v)\mapsto
\varphi^+_{k^{-1}} (v) $ is Hamiltonian, we have
$$\Pi_*^l (\varphi^+_{hgk^{-1}}(u))=\Pi_*^l (\varphi^+_{k^{-1}}
(\varphi^+_{hg}(u)))=Ad^*_ {k^{-1}} \Pi_*^l (\varphi^+_{hg} (u) )
Ad_{k^{-1}}\eqno (E.1)$$
while, since $\Pi$ is multiplicative and vanishes on $H\subset G$, we obtain
$$\Pi^l (hgk^{-1})= Ad_k \Pi^l (g) Ad^*_k.\eqno (E.2)$$
Therefore,
$$
\align
 d_{22} &= - <\iota \delta \phi,
\Pi_*^l (\varphi^+ _{hgk^{-1}} (u)) \iota \delta \psi>\\
& + < Ad_{\varphi^- _u (hgk^{-1})} \iota \delta \phi, \Pi^r_* (u)
Ad_{\varphi^-_u (hgk^{-1})} \iota\delta \psi>\\
&+ <\iota  \delta \phi , \Pi^l _* (\varphi^+_{hgk^{-1}} (u)) \Pi^l
(hgk^{-1}) 
\Pi^l _* (\varphi^+_{hgk^{-1}} (u)) \iota\delta
\psi>.\\
\endalign
$$
Next, observe that the Poisson property of $\varphi^+ : G^* \times
{\overline G} \longrightarrow G^*$
may be written as 
$$
\align
&Ad^* _{\varphi^-_u (x)} \Pi^r_* (u) Ad_{\varphi^-_u (x)}\\
&= Ad_{\varphi^+_x (u)} \big(-\Pi^l_* (\varphi^+_x (u)) \Pi^l (x) \Pi^l_*
(\varphi^+_x (u)) + \Pi^l_* (\varphi
^+_x (u))\big) Ad^*_{\varphi^+_x (u)}.\\
\endalign
$$
Thus, $d_{22}= 0$ follows from $Ad^*_v \iota Z= \iota Z,\, Z\in \fh, \, v\in
G^*.$
\medskip
(3)
$$
\align
d_{12}= &<\iota\delta\phi, \big(\Pi_*^r (u) Ad_{\varphi^-_u (hgk^{-1})} -
Ad^*_{g^{-1} h^{-1}} \Pi_*^l (\varphi_{hg}^+ (u)) Ad_{k^{-1}}\big)\iota
\delta\psi>\\
=& <\iota \delta\phi, \big(\Pi_*^r (u) Ad_{\varphi^-_u (hgk^{-1})} -
Ad^*_{(hgk^{-1})^{-1}} \Pi_*^l (\varphi_{hgk^{-1}}^+ (u)\big) \iota
\delta\psi>.\\
\endalign
$$
But it follows from Lemma 5.2.21 (b) that
$$Ad^*_{x^{-1}}\, \Pi^l_* (\varphi^+_x (u))\iota Z= Ad_{u^{-1}} \Pi^r_*
(u) Ad_{\varphi^-_u (x)} \iota Z.$$
Thus, $d_{12} =0$ follows from $Ad^*_u \iota Z= \iota Z.$
\medskip
(4) 
$$
\align
d_{1\star} = & -<\iota \delta\phi, \Pi^r_* (u) \Pi^r (\varphi^-_u
(hgk^{-1}) D\psi>\\
&-<\iota \delta\phi, Ad^*_ {g^{-1} h^{-1}} Ad_{\varphi^+_{hg} (u)^{-1}}
Ad^*_k Ad^*_ {\varphi^-_u (hgk^{-1})}D\psi> + <\iota \delta \phi, D\psi>.\\
\endalign
$$
which, upon using $\varphi^-_u (hgk^{-1})= \varphi^-_u (hg) k^{-1}$ and
the multiplicativity
$$\Pi^r (\varphi^-_u (hg) k^{-1})= \Pi^r (\varphi^-_u
(hg)),$$ becomes
$$d_{1\star}= <\iota \delta \phi, \big( -\Pi^r_* (u) \Pi^r (\varphi^-_u
(hg))  -Ad^*_ {(hg)^{-1}} Ad_{\varphi^+_{hg} (u)^{-1}} Ad^*_
{\varphi^-_u (hg)} + Id \big) D\psi>.$$
But by taking the derivative at $t=0$ in the identity
$$\varphi^-_u (e^{tX} x)= \varphi^-_u (e^{tX}) \varphi^-_
{\varphi^+_{e^{tX}} (u)} (x)$$
and dualizing immediately yields
$$Id= Ad_u Ad^*_{x^{-1}} Ad_ {\varphi_x^+ (u) ^{-1}} Ad^*_ {\varphi^-_u
(x)}+ \Pi^r_* (u) \Pi^r (\varphi^-_u (x)).\eqno (E.3)$$
Thus $d_{1\star}= 0$ again follows from the triviality of $Ad^*_{G^*}$ on
$\fh.$
\medskip
(5)
$$
\align
d_{2\star}=& <\iota \delta \phi,\big(- Ad^*_ {\varphi^-_u (hgk^{-1})}
\Pi^r_* 
(u) \Pi^r (\varphi^-_u (hgk^{-1}))Ad^*_ {\varphi^-_u (hgk^{-1})^{-1}}\\
&\hskip 10 pt + Ad^*_ {k^{-1}} \Pi^l_* (\varphi^+_{hg} (u)) \Pi^l (g)
Ad_{\varphi^+_ {hg} (u) ^{-1}}
Ad^*_k \big) D'\psi>.\\
\endalign
$$
Now, by taking the derivative at $t=0$ in the identity
$$k\varphi^-_v (e^{tX}) k^{-1}= \varphi^-_ {\varphi^+_ {k^{-1}} (v)}
(ke^{tX} k^{-1})$$
and dualizing, we find
$$Ad_ {v^{-1}} Ad^*_k= Ad^*_k Ad_ {(\varphi^+_{k^{-1}}v)^{-1}}.\eqno (E.4)$$
Therefore, 
$$
\align
d_{2\star}=& <\iota \delta \phi,\big(- Ad^*_ {\varphi^-_u (hgk^{-1})}
\Pi^r_* 
(u) \Pi^r (\varphi^-_u (hgk^{-1}))Ad^*_ {\varphi^-_u (hgk^{-1})^{-1}}\\
&\hskip 10 pt +  \Pi^l_* (\varphi^+_{hgk^{-1}} (u)) \Pi^l (hgk^{-1})
Ad_{\varphi^+_ {hgk^{-1}} (u) ^{-1}}  \big) D'\psi>,\\
\endalign
$$
where we have also used the identities (E.1) and (E.2) above.
\smallskip
Next, by taking the derivative at $t=0$ of the expression
$$\varphi^+_x (ue^{t\Lambda})= \varphi^+_ {\varphi^-_{e^{t\Lambda}} (x)}
(u) \varphi^+_x (e^{t\Lambda}),$$
we obtain
$$Id= Ad_{\varphi^+_x (u)^{-1}} Ad^* _{\varphi^-_u (x)} Ad_u
Ad^*_{x^{-1}}+ \Pi_*^l (\varphi^+_x (u)) \Pi^l (x).\eqno (E.5)$$
Combining $(E.3)$ with $(E.5)$ then yields
$$Ad^*_ {\varphi^-_u (x)} \Pi^r_* (u) \Pi^r (\varphi^-_u (x)) Ad^*_
{\varphi^-_u (x) ^{-1}}= Ad_{\varphi^+_x (u)} \Pi^l_* (\varphi^+_x (u))
\Pi^l (x) Ad_{\varphi^+_x (u)^{-1}}.$$
Thus $Ad^*_u Z= Z$ implies $d_{2\star}=0.$
\medskip
(6) 
$$
\align
&d_{\star\star} =\\
&-<D'\phi, Ad_{\varphi^-_u (hgk^{-1})^{-1}} \Pi^r
(\varphi^-_u (hgk^{-1})) \Pi^r_* (u) \Pi^r (\varphi^-_u (hgk^{-1}))
Ad^*_ {\varphi^-_u (hgk^{-1})^{-1}} D'\psi>\\
&+<D'\phi, \big( - Ad_k Ad^*_ {\varphi^+_{hg} (u) ^{-1}} \Pi^l (g) Ad_
{\varphi^+_{hg} (u) ^{-1}}Ad^*_k + \Pi^l (\varphi^-_u (hgk^{-1})) \big)
D'\psi>.\\
\endalign
$$ 
Here we use the Poisson property of $\varphi^-: {\overline {G^*}}\times
G\longrightarrow G$
which may be expressed in the form
$$
\align
&Ad_{\varphi^-_u (x)^{-1}} \Pi^r (\varphi^- _u (x)) \Pi^r_* (u)
\Pi^r(\varphi^- _u (x)) Ad^*_ {\varphi^-_u (x)^{-1}} \\
&= - Ad^*_ {\varphi^+_x (u)^{-1}} \Pi^l (x) Ad_{\varphi^+_x (u)^{-1}}  +
\Pi^l (\varphi^-_u (x)).\\
\endalign
$$
Therefore,
$$
\align
d_{\star\star}=& <D'\phi,  Ad^*_ {\varphi^+_{hgk^{-1}} (u)^{-1}} \Pi^l
(hgk^{-1}) Ad_{\varphi^+_{hgk^{-1}} (u)^{-1}} D'\psi>\\
&\hskip 10 pt   -<D'\phi,  Ad_k Ad^*_ {\varphi^+_{hg} (u) ^{-1}} \Pi^l (g)
Ad_ 
{\varphi^+_{hg} (u) ^{-1}}Ad^*_k D'\psi>.\\
\endalign
$$ 
Thus, $d_{\star\star}=0$ follows from the identities (E.2) and (E.4).
\smallskip
This concludes the verification that $lhs (P)- rhs (P)=0.$
\bigskip
(b) This follows from a general result of Weinstein according to which the
symplectic leaf passing through $(p,g,q)$ of the base $X$ in the full
symplectic realization
$$\eqalign { &\hskip 40 pt {\Cal S}\cr
&{\tilde \beta}_{\Cal H}\hskip 10 pt \swarrow\hskip 20 pt \searrow\hskip 10
pt {\tilde \alpha}_{\Cal
H}\cr
&\hskip 10 pt X\hskip 50 pt X\cr}$$
of Theorem 5.2.13 is given by ${\tilde \alpha}_{\Cal H} \big({\tilde
\beta}_{\Cal H}^{-1} (p,g,q)\big).\quad \pf$
\bigskip
To conclude the paper we give two corollaries of Thm 5.2.28.
\smallskip
First, for the special case when
$R=0$, we have $G^*=\fg^*$,
$\,I^*=\iota^*:\fg^*\rightarrow \fh^*,$
$\varphi^+_{g}=Ad^*_{g},$ and $\varphi^-_u= Id.$
Therefore, the Poisson Lie group structure on $H\times (H\times \fg^*)$ is
given by
$$
\align
&(k,(h,A))\cdot (k', (h',A'))=(kk', (hh',A+Ad^*_{h^{-1}} A')) \\
&\{f,g\}_{H\times (H\times \fg^*)} (k,(h,A))= <A, [\delta f, \delta g] >,\\
\endalign
$$
and the group action becomes
$$
\align
&{\tilde \psi}^-_0:H\times (H\ltimes g^*)\times X\longrightarrow X\\
&((k, (h,A)), (p,g,q))\mapsto (Ad^*_{h^{-1}} A + \iota^* (A), hgk^{-1},
Ad^*_{k^{-1}} q + \iota^* (Ad^*_{hgk^{-1}} A)).\\
\endalign
$$
\medskip
{\bf Corollary 5.2.29} (Symplectic leaves for $R=0$)
\smallskip
{\it (a) ${\tilde \psi}^-_0$ is a left Poisson Lie group action
\smallskip
(b) The symplectic leaf ${\Cal L}_{(p,g,q)}$ in $(X, \{\, , \, \}_X)$
passing through the point
$(p,g,q)$ is the orbit of $(p,g,q)$ under the action ${\tilde
\psi}^-_0.\quad \pf$}
\bigskip
Next, we consider the symplectic foliation of a Poisson quotient which we
now introduce. Recall from Theorem 2.1.4 (b) that $X$ has a pair
of Hamiltonian $H-$ actions with $\alpha$ and $\beta$ as momentum maps
respectively. Combining the two actions, we obtain the action
$$h\cdot (p,g,q)= (Ad^*_{h^{-1}} p, hgh^{-1}, Ad^*_{h^{-1}} q)\eqno
(5.2.15)$$
which is also Hamiltonian and its equivariant momentum map is given
by $J= \alpha - \beta.$ Now, $0\in \fh^*$ is clearly a regular value of $J$
and the corresponding isotropy subgroup is $H.$
Hence it follows from Poisson reduction \cite{MR} that $J^{-1}(0)/H$
inherits a Poisson structure $\{\, , \, \}_{J^{-1}(0)/H}$ satisfying
$$\{f_1, f_2\}_{J^{-1}(0)/H} \circ \pi= \{{\tilde {f_1}}, {\tilde {f_2}}\}_X
\circ i.\eqno (5.2.16)$$
Here, $i: J^{-1} (0)\rightarrow X$ is the inclusion map, $\pi:
J^{-1}(0)\rightarrow J^{-1}(0)/H$ is the canonical projection, $f_1,
f_2\in C^\infty (J^{-1}(0)/H),$ and ${\tilde {f_1}}, {\tilde {f_2}}$ are
(locally defined) smooth extensions of $\pi^* f_1, \, \pi^* f_2$ with
differentials vanishing on the tangent spaces of the $H-$ orbits.
\bigskip
{\bf Corollary 5.2.30} {\it The symplectic leaves of $(J^{-1}
(0)/H, \{\, , \, \}_{J^{-1}(0)/H})$ are given by the connected components
of ${\Cal L}_{(p,g,q)} \bigcap J^{-1} (0)/H,\, (p,g,q)\in X.$}
\medskip
Proof. This is a consequence of the theorem and a result in \cite{MR}, as
the triple $(X, J^{-1} (0), E)$ is Poisson reducible, where $E$ is the
tangent space to the $H-$ orbits of the action in (5.2.15).$\quad \pf$
\medskip
Clearly, the symplectic leaves of the quotient $G/H\times U$ in
Proposition 3.2.5 can also be obtained in a similar way.
\vfill
\eject

{\bf Appendix.}
\bigskip
{\bf A1. Proof of Proposition 2.2.3.}
\smallskip
The most general bivector field on $X$ is of the form
$$\eqalign {&\Pi (df,dg) (p,x,q)= K_1 (\delta_1 f, \delta_1 g)+ K_2
(\delta_2 f, \delta_2 g) +R(\delta_1 f , \delta_2 g) \cr
&\hskip 50 pt -R(\delta_1g, \delta_2 f) +P_1(\delta_1 f,
\partial g) -P_1(\delta_1g, \partial f)\cr
&\hskip 40 pt +P_2(\delta_2 f, \partial g) - P_2(\delta_2 g,
\partial f)+ P_G (\partial f, \partial g),\cr}$$
where $K_i, R, P_i , P_G$ are evaluated at $(p,x,q).$
\smallskip
Set 
$$\Omega_{(\omega, Z_1, Z_2, Z_3)}= \big( (Z_1, \omega, Z_2), (-Z_2,
T^*_y (r_{y^{-1}} \circ l_x) \omega, Z_3), (-Z_1, -T^*_{xy}
r_{y^{-1}} \omega, -Z_3)\big),$$
and denote $\big(\Pi\oplus \Pi\oplus -\Pi\big) \big( (p,x,q), (q,y,r), (p,
xy, r)\big) $ by $\Pi_m.$  Fix a reference point $q_0\in U.$
We have 
\smallskip
$$\eqalign {&\Pi_m (\Omega_{(0,Z,0,0)}, \Omega_ {(0, 0, Z',0)})=0
\Leftrightarrow R=0\cr
& \Pi_m (\Omega_{(0, Z, 0, 0)}, \Omega_{(0, Z', 0, 0)})=0
\Leftrightarrow K_1 (p,x,q)= K_1 (p,1,q_0)=: K(p)\cr
& \Pi_m (\Omega_{(0, 0, Z, 0)}, \Omega_{(0, 0, Z', 0)})=0
\Leftrightarrow K_2 (p,x,q)= -K(q).\cr}$$
Now,
$$ \Pi_m (\Omega_{(\omega, 0,0,0)}, \Omega_{(0, Z', 0,0)})= 0
\Leftrightarrow P_1(p, x,q) (Z', \omega)= P_1(p,xy, r) (Z', T^*_{xy}
r_{y^{-1}} \omega).$$
Setting successively $y=1, r=q_0$, and  $x=1, \omega= T^*_1 r_y \omega'$
in the latter equality yields
$$P_1(p,y,r) (Z', \omega')=
P_1(p,1,q_0) (Z',T^*_1 r_y \omega')=: <A_1(p) Z', T^*_1 r_y
\omega'>.$$
Similarly
$$\eqalign {&\Pi_m (\Omega_{(\omega, 0,0,0)},
\Omega_{(0,0,0,Z')})=0\Leftrightarrow\cr
& P_2 (p,x,r) (Z', \omega)= P_2
(q_0,1,r) (Z', T^*_1 l_x \omega)=: <A_2 (r) Z', T^*_1 l_x
\omega>.\cr}$$   
Moreover,
$$\Pi_m (\Omega_{(\omega, 0,0,0)},\Omega_{ (0,0,Z',0)})=0\Leftrightarrow
A_1 (p)= A_2 (p).$$
It only remains to demand that $\Pi_m (\Omega_{(\omega, 0,0,0)},
\Omega_{(\omega', 0,0,0)})=0.$
But working in the right invariant frame $P_G (\omega, \omega')= <T_1^*
r_x \omega, P (T_1^* r_x \omega')>$, the latter condition is equivalent
to the cocycle 
property
$$P(p,xy,r)= P(p,x,q) + Ad_x \, P(q,y,r) Ad^*_x.$$
Hence the assertion. $\pf$
\bigskip
\bigskip
{\bf A2. Proof of Thm 2.2.5 (b).}
\smallskip
We have to check the Jacobi identity for the bracket
$$\eqalign {&\{f,g\}_X (p,x,q)= <p,[\delta_1 f, \delta_1 g]> -
<q,[\delta_2f,\delta_2g]>\cr
&\hskip 30 pt -<A_\chi (p)\delta_1f, Dg> -<A_\chi (q) \delta_2 f, D'g>\cr
&\hskip 30 pt +<A_\chi (p) \delta_1g,
Df>+<A_\chi (q)\delta_2 g, D'f>\cr
&\hskip 30 pt+ <Df, P(p,x,q) Dg>.\cr}$$
We shall use (up to sign) the same notation $J_{ijk}$ as in the text. If
$a\in \fh$, 
we define (as in [EV]) the functions $a_1, a_2\in C^\infty (U\times
G\times 
U)$ by $a_1 (p,h,q)= <p,a>$ and $a_2 (p,h,q)= <q,a>.$ Finally, for
$Y\in \fg$, the left (resp. right) invariant vector field on $G$ whose value
at $1$ is $Y$ will be denoted by $Y^l$ (resp. $Y^r $).
\smallskip
We now compute $J_{ijk}$ for the various cases.
\smallskip
First of all, it is clear that $J_{ijk}=0,\, i,j,k\in \{1,2\}$.
\smallskip
On the other hand, we have
 $$\eqalign {J_{* 1 2}&= \{\{p^*_G f, a_1\}, b_2\}+\{\{a_1, b_2\},
p^*_G 
f\}+ \{ \{b_2, p^*_G f\},
a_1\}\cr
&= (A_\chi (q) (b) )^l \, (A_\chi (p) (a))^r (f) (x)- (A_\chi (p) (a))^r
(A_\chi (q) (b))^l (f) (x) =0.\cr
while
&\cr
J_{*11} &= \{\{p^*_G f, a_1\}, b_1\} + \{\{ a_1,
b_1\}, p^*_G f\} + \{\{b_1, p^*_G f\}, a_1\}\cr
 &= - < dA_\chi (p) \cdot ad^*_b p\cdot a , Df (x)>+ (A_\chi (p) b)^l
 (A_\chi (p) a)^l (f) (x)\cr
& - <A_\chi (p) [a,b] , Df (x)> + <dA_\chi (p)\cdot  ad^*_a p \cdot b ,
Df (x)>\cr
&\hskip 30 pt - (A_\chi (p) (a))^l (A_\chi (p) b)^l (f) (x)\cr
&= <dA_\chi (p)\cdot  ad^*_a p \cdot b - dA_\chi (p)\cdot  ad^*_b p
\cdot a\cr
&\hskip 20 pt + [A_\chi (p) a , A_\chi (p) b]
- A_\chi (p) ([a,b]), Df (x)>.\cr}$$
Similarly,
$$\eqalign {J_{* 2 2}&= <-dA_\chi (p)\cdot  ad^*_a p \cdot b +
dA_\chi (p)\cdot  ad^*_b p
\cdot a\cr
&\hskip 20 pt - [A_\chi (p) a , A_\chi (p) b]
+ A_\chi (p) ([a,b]), D'f (x)>.\cr}$$
So $J_{* i j}= 0,\, i,j\in \{1,2\} \Leftrightarrow A_\chi:
\fh^*\times \fh\longrightarrow \fg$ is a morphism of Lie algebroids.
\smallskip
Now,
$$
\align
&J_{1**}= \{a_1, \{p^*_G f, p^*_G g\}\}+
\{p^*_G f, 
\{p^*_G g, a_1\}\} + \{p^*_G g, \{ a_1, p^*_G f\}\}\\
&= <Df, \delta_1 P \cdot ad^*_a p \cdot Dg>- (A_\chi (p) a)^r (P
Dg)^r (f) \\
& - <Df, D P \cdot A_\chi (p) (a)\cdot Dg>+ (A_\chi (p) a)^r (P Df)^r (g)\\
 &+ <dA_\chi (p) \cdot (A_\chi (p)^* \,
Df) \cdot a , Dg> - (P Df)^r (A_\chi (p) a)^r (g)\\
&- <dA_\chi (p) \cdot (A_\chi (p)^* Dg) \cdot a , Df> + (P Dg)^r (A_\chi
(p) a)^r  (f)\\
&= <Df, \delta_1P \cdot ad^*_a p\cdot Dg + ad_{A_\chi (p) a} (P (Dg))>\\
&+<Df, P(ad^*_{A_\chi (p) a} Dg) - DP \cdot (A_\chi (p) a)\cdot Dg>\cr
&+<dA_\chi (p) \cdot (A_\chi (p)^* \, Df)\cdot a, Dg> - <dA_\chi (p)
\cdot (A_\chi (p)^* \, Dg)\cdot a, Df>.\\
\endalign
$$
Similarly,
$$\eqalign {&J_{2* *}= <Df, -\delta_2 P \cdot  ad^*_a q \cdot
Dg - D'P \cdot (A_\chi (q) (a))\cdot Dg>\cr
& + <dA_\chi (q) \cdot (A_\chi (q)^* D'f) \cdot a , D'g>- <dA_\chi
(q)\cdot (A_\chi (q)^* D'g)\cdot a, D'f>.\cr}$$
Writing the groupoid $1$- cocycle as
$$P(p,x,q)= -l(p) + \pi (x) + Ad_x l(q) Ad_x^*,$$
we have
$$\eqalign {&\delta_1 P \cdot \Lambda= - dl (p) \cdot \Lambda\cr
& DP \cdot X= d\pi (1) (X) + ad_X \, \pi (x) + \pi (x)\,
ad^*_X \cr
&\hskip 20 pt + ad_X\, Ad_x \, l(q)\, Ad^*_x+ Ad_x \, l(q) Ad^*_x \,
ad^*_X\cr
&\delta_2 P \cdot \Lambda = Ad_x \, dl(q)\cdot \Lambda \, Ad^*_x\cr
& D'P\cdot X= Ad_x\, d\pi (1) \cdot X \, Ad_x^* + Ad_x\, ad_X\,
l(q)\, Ad^*_x \cr
&\hskip 20 pt + Ad_x \, l(q) ad^*_X\, Ad^*_x\cr}$$
Inserting the latter into $J_{1**}$ and $J_{2**}$ yields
$$
\align
 &\hskip 40 pt J_{* * 1}= 0\Leftrightarrow J_{**
2}= 0 
\Leftrightarrow\\
& <\alpha, (dl(p) ad^*_a p  + ad _{A_\chi (p) a}
l(p)  + l(p) ad^*_{A_\chi (p)a} + d\pi (1) A_\chi(p) a) \beta>\\
& =  + <dA_\chi(p) \cdot (A_\chi^* (p)
\alpha)\cdot a, \beta>  - <dA_\chi (p)\cdot (A_\chi(p)^* \beta) \cdot a,
\alpha>,\\
&\hskip 60 pt \forall \alpha,\beta\in \fg^*, a\in \fh.\\
\endalign
$$
\smallskip
Finally,
$$
\align
 J_{***}=& \{p^*_G f, \{p^*_G g, p^*_G
h\}\}+ c.p.(f,g,h)\\
&= <Dg, \delta_1 P \cdot (A_\chi (p)^*\, Df ) \cdot Dh+ \delta_2 P
\cdot (A_\chi (p)^* D'f)\cdot Dh>\\
&\hskip 20 pt - (P Df)^r (P Dh)^r (g) + (P Df)^r (P Dg)^r (h) \\
&\hskip 10 pt  - <Dg, DP\cdot (P Df)\cdot Dh> + c.p. (f,g,h),\\
\endalign
$$
which may easily be brought to the form stated in Thm 2.2.5 (b).
$\pf$
 \bigskip
\bigskip
{\bf A3. Proof of Theorem 5.1.4.}
\smallskip
We have to show that the graph of the multiplication
$$Gr (m) \subset \Gamma\times \Gamma\times {\overline \Gamma}$$
is a coisotropic submanifold.  We use (as in (5.1.6)) the notation $\phi^l
_h (v)= \phi^l _v (h) = \varphi^+ _{h^{-1}} (v),\, h\in H, v\in G^*.$ We
have
$$
\align
&Gr (m)=\\
& \{ \big( (h, Ad^*_{k^{-1}} q + I^* (v), u), (k,q,v),
(hk, q, u\phi^l _h (v))\big) \, \mid \, h,k\in H, u,v \in G^*, q\in
\fh^*\},\\
\endalign
$$
therefore 
$$\eqalign {  T_\ast (Gr (m))&= \{ \big( (Z_1, -Ad^*_{k^{-1}} \,
ad^* _{(T_k 
l_{k^{-1}} Z_2)} q + Ad^*_{k^{-1}} \lambda + T_v (I^*) V, U),\cr
(Z_2, \lambda, V), &\, (T_k l_h Z_2 + T_h r_k Z_1, \lambda, T_u r_{\phi^l _h
(v)} U + T_{\phi^l _h (v)} l_u (T_h \phi^l _v Z_1+ T_v
\phi^l _h V)\big)\cr
&\mid\,  Z_1\in T_h H , Z_2\in T_k H, \lambda \in \fh^*, U\in T_u G^*,
V\in 
T_v G^*\}.\cr}$$
Hence $\Omega\in \big( T_\ast Gr (m)\big)^\perp$ if and only if
$$\eqalign {\Omega&= \big( (-T^*_h r_k \mu - T^*_h (l_u\circ \phi^l
_v) A, 
z_1, -T^*_u r_{\phi^l _h (v)} A),\cr
&( - T^*_k l_h \mu - T^*_k l_{k^{-1}} ad^*_ {(Ad_{k^{-1}} z_1)}
q, z_2, - T_v ^* I^* z_1 - T^*_v (l_u\circ \phi^l _h ) A),\cr
& (\mu, - Ad_{k^{-1}} z_1 -z_2, A)\big),\cr}$$
for some $\mu\in T^*_{hk} H, z_1, z_2\in \fh, A\in T^*_{u\phi^l_h
(v)} G^*.$
For $\Omega, \Omega'\in \big( T_\ast Gr (m) \big)^\perp$, we have
$$
\align
& (\Pi \oplus \Pi \oplus -\Pi) (\Omega, \Omega')= < T^*_h
r_k 
\mu'+ T^*_h (l_u\circ \phi^l _v )A', T_1 l_h z_1>\\
&- < T^*_h r_k \mu+ T^*_h (l_u\circ \phi^l _v )A, T_1 l_h
z'_1>- <Ad^*_{k^{-1}} q + I^* (v), [z_1, z'_1]>\\
& -< T^*_u r_{\phi^l _h (v)} A, \lambda^+ (T_1^* l_u (T^*_u
r_{\phi^l _h (v)} A')) (u)>\\
&+ < T^*_k l_h \mu' + T^*_k l_{k^{-1}} ad^*_{(Ad_{k^{-1}}
z'_1)} q , T_1 l_k z_2>\\
& - < T^*_k l_h \mu + T^*_k l_{k^{-1}}
ad^*_{(Ad_{k^{-1}}
z_1)} q , T_1 l_k z'_2>- <q, [z_1, z'_1]>\\
& - < T^*_v I^* z_1 + T^*_v (l_u \circ
\phi^l _h) A, \lambda^+ (T^*_1 l_v (T^*_v I^* z'_1 +
T^*_v (l_u\circ \phi^l _h ) A'))(v)>\\
& - <\mu', T_1 l_{hk} (Ad_{k^{-1}} z_1 + z_2)> +  <\mu, T_1 l_{hk}
(Ad_{k^{-1}} z'_1 + z'_2)>\\
& + <q, [Ad_{k^{-1} } z_1 + z_2, Ad_{k^{-1}} z_1' + z_2']> + <A,
\lambda^+ (T^*_1 l_{u\phi^l _h (v)} A') (u\phi^l _h (v))>.\\
\endalign
$$
We now treat separately the three types of terms which do not obviously
cancel out:
\smallskip 
$(1):= -<I^* (v), [z_1, z'_1]> - <T^*_v I^* z_1, \lambda^+
(T^*_1 l_v \, T^*_v I^* z'_1) (v)>.$
\smallskip
Since $I^*$ is a morphism of groups, $I^*\circ l_v (w)= I^* (u)+
I^* (w)$ therefore $T_1
(I^* \circ l_v) = \iota^*$, and hence
$$\eqalign { T_v I^* \lambda^+
(T^*_1 (I^*\circ l_v ) z'_1) (v)&= T_v I^* \lambda^+ (\iota
(z'_1)) (v)\cr
&= {d\over dt}_{\mid_0} I^* (\phi^l_{e^{t\iota (z'_1)}} (v)\cr
&= {d\over dt}_{\mid_0} Ad^*_{e^{t\iota (z'_1)}} (I^* (v))=
ad^*_{z'_1} I^* (v).\cr}$$ Thus $(1)=0.$
\smallskip
$$
\align
(2):=& -< T^*_u r_{\phi^l _h (v)} A,\lambda^+ (T_1^* l_u
(T^*_u 
r_{\phi^l _h (v)} A')) (u)>\\
&- <T^*_v (l_u \circ
\phi^l _h) A,\lambda^+ (T^*_1 l_v (
T^*_v (l_u\circ \phi^l _h ) A'))(v)>\\
&+ <A, \lambda^+ (T^*_1
l_{u\phi^l _h (v)} A') (u\phi^l _h (v))>.\\
\endalign
$$
Let $\Pi_\ast$ be the Poisson tensor of $G^*.$ We have
$$\eqalign {& (2)=\Pi_\ast (u) (T^*_u r_{\phi^l _h (v)} A, \, T^*_u
r_{\phi^l _h 
(v)} A')\cr
& + \Pi_\ast (v) (T^*_v (l_u\circ \phi^l _h) A, \,T^*_v (l_u\circ
\phi^l _h) A') - \Pi_\ast (u\phi^l _h (v)) (A, A')\cr
&= + \Pi_\ast (v) (T^*_v (l_u\circ \phi^l _h) A, \,T^*_v (l_u\circ
\phi^l _h) A') - \Pi_\ast (\phi^l_h (v))  (T^*_{\phi^l _h (v)} l_u A,
\,T^*_{\phi^l _h (v)}l_u A'),\cr}$$
where in the last equality we have used the multiplicativity of
$\Pi_\ast$.
 Thus $(2)=0$ follows from the Hamiltonian property of $\phi^l _h :
 G^*\longrightarrow G^*$.
\smallskip
$(3):= <T^*_h (l_u\circ \phi^l _v )A', T_1 l_h z_1>  - < T^*_v
I^* z_1 
, \lambda^+ (T^*_1 l_v (T^*_v (l_u\circ \phi^l _h ) A'))(v)>.$
\smallskip
We have 
$$\eqalign {<T^*_h (l_u\circ \phi^l _v )A', T_1 l_h z_1> &= <A', T_1
(l_u\circ \phi^l _v \circ l_h) z_1>\cr
&= <A', T_1 (l_u\circ \phi^l_h\circ
\phi^l _v)z_1>\cr
&= <T^* _v (l_u\circ \phi^l _h) A', T_1 \phi^l _v
z_1>,\cr}$$
and $T^*_v I^* = T^*_v l_{v^{-1}} \circ \iota.$
Therefore
$$\eqalign {(3)=& <T^* _v (l_u\circ \phi^l _h) A', T_1 \phi^l _v
(z_1)>+ <T^* _v l_{v^{-1}} \iota (z_1), \Pi_\ast (v) (T^*_v
(l_u\circ \phi^l _h) A')>\cr
=& <T^* _v (l_u\circ \phi^l _h) A',\, T_1 \phi^l _v
(z_1) -\Pi_* (v) (T^*_v l_{v^{-1}} \iota (z_1))>= 0.\cr}$$
Hence the proof. $\pf$


\bigskip
\bigskip

\Refs

\ref\key AM
\by Almeida, R. and Molino, P.
\paper Suites d'Atiyah et feuilletages transversalement complets
\jour C. R. Acad. Sci. Paris, Serie I, t. 300
\yr 1985
\pages 13-15
\endref


\ref\key BD
\by Belavin, A.A.  and Drinfel'd, V.
\paper Triangle equations for simple Lie algebras
\jour Mathematical Physics Reviews (Ed. Novikov et al.)
 Harwood, New York
\yr 1984
\pages 93-165\endref

\ref\key BH
\by Brown, R. and Higgins, P. J.
\paper On the connection between the second relative
homotopy groups of some related spaces
\jour Proc. London Math. Soc.
\vol 36
\yr 1978 \pages 193-212
\endref

\ref\key BKS
\by  Bangoura, M. and Kosmann-Schwarzbach, Y.
\paper Equations de Yang-Baxter dynamique classique et alg\'ebroides de Lie
\jour C. R. Acad. Sc. Paris, Serie I
\vol 327 
\yr 1998 \pages 541-546\endref

\ref\key CDW
\by Coste, A., Dazord, P. and Weinstein, A.
\paper Groupoides symplectiques
\jour Publications du D\'epartment de Math\'ematiques de l'Universit\'e de
Lyon
\vol 2/A \yr 1987 \pages1-65
\endref

\ref\key{CF}
\by Crainic, M. and Fernandes, R.
\paper Integrability of Lie brackets
\jour LANL e-print archive
\newline  math.DG/0105033; http://xxx.lanl.gov/
\endref

\ref\key D
\by  Drinfeld, V.G.
\paper Hamiltonian structures on Lie groups, Lie bialgebras, and the
geometric
meaning of the classical Yang-Baxter equations
\jour Soviet Math. Dokl.
\vol 27 \yr 1983\pages 667--671
\endref 

\ref\key {DSW}
\by Cannas da Silva, A. and Weinstein, A.
\paper Geometric models for noncommutative algebras
\jour Berkeley Mathematics
Lecture Notes 10. Amer. Math. Soc., Providence, RI \yr 1999
\endref

\ref\key E
\by Ehresmann, C.
\paper Cat\'egories structur\'ees
\jour Ann. Sci. Ecole Norm. Sup.
\vol 80 \yr 1963\pages 349-426
\endref


\ref\key{EV}
\by Etingof, P. and Varchenko, A.
\paper Geometry and classification of solutions of the classical dynamical
Yang-Baxter equation
\jour Commun. Math. Phys.\vol 192\yr 1998 \pages77-120
\endref


\ref\key{F}
\by Felder, G.,
\paper Conformal field theory and integrable systems associated to elliptic
curves
\jour  Proc. ICM Zurich, Birkh\"auser, Basel\yr 1994\pages1247--1255
\endref

\ref\key HM
\by Hurtubise, J., Markman, E.
\paper Elliptic Sklyanin integrable systems for arbitrary reductive groups
\jour LANL e-print Archive math.AG/0203031
\endref

\ref\key K
\by Karasev, M.
\paper Analogues of objects of the theory of Lie groups for nonlinear
Poisson brackets
\jour Math. USSR Izvestiya\vol 28\yr 1987
\pages 497-527
\endref

\ref\key{KW}
\by Kinyon,M. and  Weinstein, A.,
\paper Leibniz algebras, Courant algebroids, and multiplications on
reductive 
homogeneous spaces.
\jour Amer. J. math. \vol 123\yr 2001 \pages 525-550
\endref

\ref\key L
\by Libermann, P.
\paper On symplectic and contact groupoids
\jour Differential geometry and its applications (Opava 1992), Math. Publ.1,
Silesian Univ. Opava, Opava\yr 1993
\pages 29-45
\endref

\ref\key{Lu} 
\by Lu, J.H.
\book Ph.D. Thesis, Berkeley\yr 1990
\endref


\ref\key {LW1} 
\by Lu, J.-H., Weinstein, A.
\paper Poisson Lie groups, dressing transformations, and Bruhat
decompositions.
\jour J. Diff. Geom. 31\yr 1990 \pages 501--526
\endref

\ref\key {LW2}
\by Lu, J.-H., Weinstein, A.
\paper Groupoides symplectiques doubles des groupes de Lie-Poisson
\jour C. R. Acad. Sci. Paris, Ser. I, t. 309
\yr 1989\pages 951-954
\endref


\ref\key {LWX}
\by Liu, Z.-J., Weinstein, A. , and Xu,P.
\paper Manin triples for Lie bialgebroids
\jour J. Diff. Geom.\vol 45 \yr 1997 \pages 547-574
\endref


\ref\key {LX1}
\by  Li, L. C. and Xu, P.
\paper Spin Calogero-Moser
systems associated with simple Lie algebras
\jour C. R. Acad. Sci.
Paris, t.331, S\'erie I\yr 2000\pages 55--60
\endref

\ref\key{LX2}
\by Li, L. C. and Xu, P.
\paper A class of integrable spin Calogero-Moser systems
\jour to appear in Commun. Math. Phys.
\endref

\ref\key{M1}
\by  MacKenzie, K.
\book Lie groupoids and Lie algebroids in differential geometry. LMS Lecture
Notes Series 124
\publ Cambridge University Press\yr 1987
\endref

\ref\key{M2}
\bysame
\paper On symplectic double groupoids
and the duality of Poisson
groupoids
\jour Internat. J. Math. \vol 10 \yr 1999\pages 435-456
\endref

\ref\key {M3}
\bysame 
\paper Double Lie algebroids and second-order geometry I.
\jour Adv. Math. \vol 94 \yr 1992 \pages 180-239
\endref

\ref\key{MR}
\by Marsden J.E. and Ratiu T.
 \paper Reduction of Poisson manifolds
\jour Letters in Math. Phys.\vol 11\yr 1986\pages 161-169
\endref

\ref\key{MW}
\by Mikami K. and Weinstein A.
\paper Moments and reduction for symplectic groupoids
\jour Publ. RIMS, Kyoto University\vol 24\yr 1988\pages 121-140
\endref

\ref\key{MX1} 
\by MacKenzie, K. and Xu, P.
\paper Lie bialgebroids and Poisson groupoids
\jour Duke Math. J.\vol 73\yr 1994\pages 415--452
\endref

\ref\key{MX2}
\bysame  
\paper Integration of Lie bialgebroids
\jour Topology \vol 39\yr 2000 \pages 445-467
\endref 

\ref\key P
\by Pradines, J.
\paper G\'eometrie differentielle au-dessus d'un groupoide
\jour C. R. Acad. Sc. Paris, Ser. I, t. 266\yr 1968
\pages 1194-1196
\endref

\ref\key{STS}
\by  Semenov-Tian-Shansky, M.
\paper Dressing transformations and Poisson Lie group actions
\jour Publ. RIMS, Kyoto University \vol 21\yr 1985 \pages1237-1260
\endref

\ref\key{W1}
\by Weinstein, A.
\paper Coisotropic calculus and Poisson groupoids
\jour J. Math. Soc. Japan 4 no. 40\yr 1988 \pages 705--727
\endref  

\ref\key{W2}
\by Weinstein, A.
\paper Symplectic groupoids and Poisson manifolds
\jour Bull. Amer. Math. Soc.
\vol 16\yr 1987\pages 101-104
\endref

\ref\key{X}
\by  Xu, P.
\paper On Poisson groupoids.
\jour Internat. J. Math. \vol 6
\yr 1995\pages 101-124
\endref

\endRefs
\enddocument